\documentclass[multiarticle,nochapter,download]{msjmemoirs}
\usepackage{amsmath,amsthm,amssymb}
\usepackage{bbm}

\usepackage{makeidx}
\makeindex

\usepackage[nottoc]{tocbibind}

\usepackage{etoolbox}
\patchcmd{\theindex}{\thispagestyle{plain}}{}{}{}

\usepackage{epic,eepic}
\usepackage[all]{xy}

\allowdisplaybreaks

\newtheorem{thm}{Theorem}[section]
\newtheorem{prop}[thm]{Proposition}

\newtheorem{cor}[thm]{Corollary}
\newtheorem{lem}[thm]{Lemma}

\theoremstyle{definition}

\numberwithin{equation}{section}
\newtheorem{rem}[thm]{\bf Remark}
\newtheorem{ex}[thm]{\bf Example}

\newtheorem{defn}[thm]{\bf Definition}

\newtheorem{obs}[thm]{\bf Observation}

\def\bbA{\mathbb{A}}
\def\bbC{\mathbb{C}}
\def\bbP{\mathbb{P}}
\def\bbQ{\mathbb{Q}}
\def\bbR{\mathbb{R}}
\def\bbT{\mathbb{T}}

\def\bbZ{\mathbb{Z}}
\def\bfr{\mathbf{S}}
\def\bbone{\mathbbm{1}}

\def\rmsf{\mathrm{sf}}

\def\bfB{\mathbf{B}}
\def\bfC{\mathbf{C}}
\def\bfF{\mathbf{F}}
\def\bfG{\mathbf{G}}
\def\bfP{\mathbf{P}}

\def\bfX{\mathbf{X}}
\def\bfY{\mathbf{Y}}
\def\bfSigma{\mathbf{\Sigma}}
\def\bfa{\mathbf{a}}
\def\bfb{\mathbf{b}}
\def\bfc{\mathbf{c}}
\def\bfd{\mathbf{d}}
\def\bfe{\mathbf{e}}
\def\bfg{\mathbf{g}}

\def\bfr{\mathbf{r}}

\def\bfu{\mathbf{u}}
\def\bfv{\mathbf{v}}
\def\bfx{\mathbf{x}}
\def\bfy{\mathbf{y}}
\def\bfz{\mathbf{z}}
\def\bfzero{\mathbf{0}}
\def\rmGr{\mathrm{Gr}}

\def\calA{\mathcal{A}}

\def\calF{\mathcal{F}}
\def\calL{\mathcal{L}}

\def\calU{\mathcal{U}}
\def\haty{\hat{y}}

\usepackage{etoolbox}
\makeatletter
\pretocmd{\@startsection}{\gdef\thesectiontype{#1}}{}{}
\pretocmd{\@sect}{\@namedef{the\thesectiontype title}{#8}}{}{}
\pretocmd{\@ssect}{\@namedef{the\thesectiontype title}{#5}}{}{}
\makeatother

\usepackage{fancyhdr}
\begin{document}
\bibliographystyle{amsalpha}
\fontsize{11pt}{12.5pt}\selectfont

\mainmatter[page=1]
%\title[Part III.\ \thesection.\   \thesectiontitle]{Cluster Algebras and Cluster Patterns\\
\title[\ ]{Cluster Algebras and Scattering Diagrams\\
\ \\
Part I\\
Basics in Cluster Algebras\thanks{
This is the final manuscript of
 Part I of the monograph ``Cluster Algebras and Scattering Diagrams'',
 MSJ Mem. 41 (2023) by the author. 
%This is a preliminary draft of
% Part I of the forthcoming monograph ``Cluster Algebras and Scattering Diagrams''
%by the author. Any comments are welcome.
}}
%\subtitle{subtitle}
%\author[\thesubsection.\ \thesubsectiontitle]{Tomoki Nakanishi\\
\author[\ ]{Tomoki Nakanishi\\
\small Graduate School of Mathematics, Nagoya University}
%\dedicatory{}
%\keywords{keywords-list}
%\subjclass{Primary 13F60}
\maketitle

\begin{abstract}
This  is a first step guide to
the theory of cluster algebras.
We especially focus  on
basic notions, techniques, and results concerning seeds,
cluster patterns, and cluster algebras.
\end{abstract}

\newpage
\tableofcontents[depth=2]

\newpage

% last updated on 2022/12/8 by T.N.

%%%%%%%%%%%%%%%%%%%%%%%%%%%%%%%%%%%%%
%%%%%% part1.tex %%%%%%%%%%%%%%%%%%%%%%%%
\setcounter{section}{-1}
\pagestyle{fancy}
% there is no subsection here; so we temporarily reset it
\fancyhead[CO]{\thesection.\ \thesectiontitle}

\section{Introduction to Part I}

This  is a first step guide to
the theory of cluster algebras.

Cluster algebra theory was founded by the series of
papers entitled as ``Cluster Algebras I--IV'' (commonly
referred to as CA1--4)
by Fomin and Zelevinsky, together with Berenstein
\cite{Fomin02,Fomin03a,Berenstein05,Fomin07}.
These fundamental papers are still the best sources to learn
the basic notions, formulations, techniques, important examples, and the perspective
on cluster algebras.
However, even for someone  who seriously wishes to learn  cluster algebra theory,
 reading through all of them, even selectively, at the beginning is not an easy task.
 This is partly because each of CA1--4 has multiple goals,
 and some of them are rather advanced and/or specific
 for the beginners.

Fortunately,  several excellent texts/monographs/reviews (including preliminary drafts) are already available
to fill the gap.
(e.g., \cite{Carter06}, \cite{Keller08}, \cite{Gekhtman10}, \cite{Keller12}, \cite{Marsh13}, \cite{Williams12}, \cite{Glick17}, \cite{Fomin16},
\cite{Fomin17}, \cite{Fomin20}, \cite{Fomin21}).
As another text of the same kind, 
we especially focus  on \emph{basic aspects\/} in cluster algebra theory.
To be more explicit,
we paid attention
in the following points while preparing the manuscript.

\begin{itemize}
\item
Most importantly, the text is aimed for anyone
who starts to learn cluster algebras seriously  without any 
preliminary knowledge.
\item
The text should be read in a self-contained way.
Also, it should be reasonably concise.
 \item
Consequently,
we focus on (most but not all of) the basic and fundamental aspects of cluster algebra theory,
namely, basic notions,  techniques, and results  concerning seeds,
cluster patterns, and cluster algebras in CA1--4.
\item
We especially employ the formulation of seeds and mutations in CA4 throughout the text,
 because the structure of mutations
is most transparent in our point of view.
By the same reason, 
in some part of proofs,
we use the notion of \emph{free coefficients\/}
(coefficients in a universal semifield)
that replaces the role of coefficients of geometric type
in CA1--4.
\item
We try  to present
elementary manipulations of mutations
in the proofs of basic results explicitly,
even if many of them are easy or straightforward.
This is based on our belief that
the essence of cluster algebras is found  
in these details of mutation mechanism.
\item
Naturally, we do not try to go  deep into a particular topic.
In particular, we  have to omit the proof
of the finite type classification of cluster algebras in CA2,
because it requires deep analysis in conjunction with root systems of finite type.
\item
We add a short 
section on generalized cluster algebras
as extra material beyond CA1--4.
\end{itemize}

All results presented in this text are  well-known,
and, except for the last section,
they are taken from  CA1--4 or easily derived from the results therein 
unless otherwise mentioned.

Bon voyage!

\newpage
% we  reset it again
\fancyhead[CO]{\thesubsection.\ \thesubsectiontitle}

\section{Getting Started}
\subsection{The very first step: pentagon periodicity}

\label{1sec:first1}

Let us play with a prototypical example of the cluster algebraic structure
without explaining  the background.

Let $x_1, x_2$ be formal variables.
We consider the following recursion, or discrete dynamical system with discrete time $t=0,1,2,\dots$:

\begin{itemize}
\item
initial condition:  $x_1(0)=x_1$, \ $x_2(0)=x_2$.
\item
time development: For even $t$,
\begin{align}
\label{1eq:xA21}
\begin{cases}
\displaystyle
x_1(t+1)=\frac{1}{x_1(t)}(x_2(t)+1),\\
x_2(t+1)=x_2(t),
\end{cases}
\end{align}
and, for odd $t$,
\begin{align}
\label{1eq:xA22}
\begin{cases}
x_1(t+1)=x_1(t),\\
\displaystyle
x_2(t+1)=\frac{1}{x_2(t)}(x_1(t)+1).
\rule{0pt}{18pt}\hskip-1pt
\end{cases}
\end{align}
\end{itemize}

Let us calculate $x_1(t)$ and $x_2(t)$ up to $t=5$.
Anyone who wishes to learn cluster algebras seriously has to work on the
following calculations by oneself.
\begin{align}
&
\begin{cases}
x_1(0)=x_1,\\
x_2(0)=x_2.
\end{cases}\\
&
\begin{cases}
\displaystyle
x_1(1)=\frac{x_2+1}{x_1},\\
x_2(1)=x_2.
\end{cases}
\\
&
\begin{cases}
\displaystyle
x_1(2)=\frac{x_2+1}{x_1},\\
\displaystyle
x_2(2)
%&=\frac{1}{x_2(1)}(x_1(1)+1)\\
=\frac{1}{x_2}\frac{x_1+x_2+1}{x_1}
=\frac{x_1+x_2+1}{x_1 x_2}
.
\rule{0pt}{17pt}\hskip-1pt
\end{cases}
\end{align}
So far, nothing special happened.
Let us continue.
\begin{align}
\label{1eq:x13}
\begin{cases}
\displaystyle
x_1(3)
%&=\frac{1}{x_1(2)}(x_2(2)+1),\\
=\frac{x_1}{x_2+1}
\frac{x_1x_2+x_1+x_2+1}{x_1x_2}
\buildrel  ! \over{=}
\frac{x_1+1}{x_2},
\\
\displaystyle
x_2(3)
=\frac{x_1+x_2+1}{x_1 x_2}.
\rule{0pt}{18pt}\hskip-1pt
\end{cases}
\end{align}
Now, at the equality with symbol ``$!$'' in \eqref{1eq:x13}, the first nontrivial reduction happens.
Let us continue, and we  mark the equality 
with symbol ``$!$''
whenever a similar reduction
occurs.

\begin{align}
&
\begin{cases}
\displaystyle
x_1(4)=\frac{x_1+1}{x_2},\\
\displaystyle
x_2(4)
%&=\frac{1}{x_2(3)}(x_1(3)+1)\\
=\frac{x_1x_2}{x_1+x_2+1}\frac{x_1+x_2+1}{x_2}
\buildrel  ! \over{=} x_1,
\rule{0pt}{18pt}\hskip-1pt
\end{cases}
\\
&
\begin{cases}
\displaystyle
x_1(5)
%&=\frac{1}{x_1(4)}(x_2(4)+1),\\
=\frac{x_2}{x_1+1} (x_1+1)
\buildrel  ! \over{=}
x_2,
\\
\displaystyle
x_2(5)
=x_1.
\end{cases}
\end{align}
After  $t=3$,  reductions occur  \emph{systematically},
so that at $t=5$ the result is  extremely reduced.

Let us observe the above result more closely.

\begin{obs}
\label{1obs:x1}
(a). \emph{Periodicity.} 
We have a \emph{half\/} periodicity at $t=5$.
If we continue calculation, we have a \emph{full\/} periodicity at $t=10$.
\index{pentagon periodicity}
\par
(b). \emph{Laurent phenomenon.}
Thanks to the reduction,
each $x_i(t)$ is expressed, not as a general rational function,
but as a \emph{Laurent polynomial\/} in the initial variables
$x_1$, $x_2$ with integer coefficients.
\par
(c). \emph{Laurent positivity.} 
Moreover, every nonzero coefficient of 
the above Laurent polynomial is \emph{positive}.
\end{obs}

\begin{rem}
The property (c) is nontrivial even though the  rules
\eqref{1eq:xA21} and \eqref{1eq:xA22} do not contain any negative coefficients.
For example, we have
\begin{align}
\frac{x^3+1}{x+1}=x^2-x +1.
\end{align}
\end{rem}

Next, let us consider a similar but different system.

Let $y_1, y_2$ be another formal variables.
We consider the following recursion, or discrete dynamical system with discrete time $t=0,1,2,\dots$:
\begin{itemize}
\item
initial condition:  $y_1(0)=y_1$, \ $y_2(0)=y_2$.
\item
time development: For even $t$,
\begin{align}
\label{1eq:yA2}
\begin{cases}
y_1(t+1)=y_1(t)^{-1},\\
y_2(t+1)=y_2(t) (1+y_1(t)),
\end{cases}
\end{align}
and, for odd $t$,
\begin{align}
\begin{cases}
y_1(t+1)=y_1(t)(1+y_2(t)),\\
y_2(t+1)=y_2(t)^{-1}.
\end{cases}
\end{align}
\end{itemize}

We calculate up to $t=5$.
As before, we mark the symbol  ``$!$'' whenever some
reduction occurs.
Again, it is important to do it oneself.

\begin{align}
&
\begin{cases}
y_1(0)=y_1,\\
y_2(0)=y_2,
\end{cases}
\\
&
\begin{cases}
y_1(1)=y_1^{-1},\\
y_2(1)=y_2(1+y_1),
\end{cases}
\\
&
\begin{cases}
\displaystyle
y_1(2)=\frac{1+y_2+y_1y_2}{y_1},
\\
\displaystyle
y_2(2)=\frac{1}{y_2(1+y_1)}
,
\rule{0pt}{20pt}\hskip-1pt
\end{cases}
\\
&
\begin{cases}
\displaystyle
y_1(3)
=\frac{y_1}{1+y_2 + y_1y_2},\\
\displaystyle
y_2(3)=\frac{1}{y_2(1+y_1)}
\frac{1+y_1+y_2+y_1y_2}{y_1}
\buildrel  ! \over{=}
\frac{1+y_2}{y_1y_2},
\rule{0pt}{18pt}\hskip-1pt
\end{cases}
\\
&
\begin{cases}
\displaystyle
y_1(4)=\frac{y_1}{1+y_2 + y_1y_2}
\frac{1+y_2+y_1y_2}{y_1y_2}
\buildrel  ! \over{=} \frac{1}{y_2},
\\
\displaystyle
y_2(4)=\frac{y_1y_2}{1+y_2},
\rule{0pt}{18pt}\hskip-1pt
\end{cases}
\\
&
\begin{cases}
y_1(5)
=y_2,\\
\displaystyle
y_2(5)=\frac{y_1y_2}{1+y_2}
\frac{1+y_2}{y_2}
\buildrel  ! \over{=}
y_1.
\end{cases}
\end{align}

Unlike the previous case, the variables $y_i(t)$ are not necessarily 
Laurent polynomials in the initial variables $y_1$, $y_2$.
Nevertheless, similar reductions systematically occur after $t=3$,
and the same periodicity is obtained.
It is natural to speculate that there is a close relation between two systems.

Indeed this is the simplest example of mutations of \emph{cluster variables\/} 
 for $x_i(t)$
and \emph{coefficients\/} 
 for $y_i(t)$,  which we are going to study.
The periodicity we observed is the celebrated \emph{pentagon periodicity\/}
for a cluster algebra of type $A_2$.

\subsection{Semifields}
\label{1sec:semifields1}

We introduce the notion of a semifield.
There are some variations of the definition depending on the purpose, and we especially use the one
in \cite{Fomin03a}.

\begin{defn}[Semifield]\index{semifield}
A multiplicative abelian group $\mathbb{P}$
equipped with 
 a binary operation
$\oplus$
 is called a \emph{semifield\/} if 
the following properties hold:
For any $a,b,c,\in \mathbb{P}$,
\begin{align}
\label{1eq:sf1}
a\oplus b &= b\oplus a,\\
(a\oplus b)\oplus c &= a\oplus (b\oplus c),\\
\label{1eq:sf3}
(a\oplus b) c &= ac\oplus bc.
\end{align}
The operation $\oplus$ is called the \emph{addition\/} in $\mathbb{P}$.
Note that there is no \emph{subtraction\/} in $\mathbb{P}$.
\end{defn}

\begin{ex}
The set of all positive rational numbers $\mathbb{Q}_+$
is a semifield by the usual multiplication and addition.
On the other hand, 
the set of all nonnegative rational numbers $\mathbb{Q}_{\geq 0}$
is \emph{not\/} a semifield by the usual multiplication and addition,
because it is not a multiplicative group due to the presence of $0$.
Similarly, 
the set of all positive real numbers $\mathbb{R}_+$
is a semifield by the usual multiplication and addition.
\end{ex}

The following three examples are especially important for cluster algebras.

\begin{ex}
\label{1ex:sf1}
(a). \emph{Universal semifield} 
$\mathbb{Q}_{\mathrm{sf}}(\bfu)$.
\index{semifield!universal} 
Let $\bfu=(u_1,\dots, u_n)$ be an $n$-tuple of 
formal variables.
Let $\bbQ(\bfu)$ be the  rational function field of $\bfu$.
We say that a 
rational function $f(\bfu)\in \bbQ(\bfu)$ 
has a \emph{subtraction-free expression\/}
if it is expressed as
$f(\bfu)=p(\bfu)/q(\bfu)$, where both $p(\bfu)$ and $q(\bfu)$
are \emph{nonzero\/} polynomials in $\bfu$ whose coefficients
are \emph{nonnegative\/} integers.
\index{subtraction-free expression} 
For example, $f(\bfu)=u_1^2-u_1+1=(u_1^3+1)/(u_1+1)$ has a subtraction-free expression.
Let   $\mathbb{Q}_{\mathrm{sf}}(\bfu)$
be the set of all rational functions in $\bfu$
 having subtraction-free expressions.
Then, $\mathbb{Q}_{\mathrm{sf}}(\bfu)$
is a semifield by the usual 
multiplication and addition in $\bbQ(\bfu)$.
\par
(b). \emph{Tropical semifield}
 $\mathrm{Trop}(\bfu)$.
\index{semifield!tropical}
Let $\bfu=(u_1,\dots, u_n)$ be an $n$-tuple of 
formal variables.
Let $\mathrm{Trop}(\bfu)$ be the set of all
Laurent monomials in $\bfu$ with coefficient 1,
which is a multiplicative abelian group by the usual
multiplication.
We define the addition $\oplus$ by
\begin{align}
\label{1eq:ts1}
\prod_{i=1}^n u_i^{a_i} \oplus
\prod_{i=1}^n u_i^{b_i}
:=
\prod_{i=1}^n u_i^{\min(a_i,b_i)} .
\end{align}
Then,  $\mathrm{Trop}(\bfu)$ becomes a semifield.
The addition $\oplus$ is called the \emph{tropical sum}.
\index{tropical sum}

\par
(c). \emph{Trivial semifield\/}
 $\bbone$. \index{semifield!trivial}
Let $\bbone=\{1\}$ be the trivial multiplicative group.
We define the addition by
$1\oplus 1 =1$.
Then, $\bbone$ becomes a semifield.
\end{ex}

Let $\bbP$ be any semifield.
For any $a\in \bbP$ and any positive integer $m$,
 we write
\begin{align}
\label{1eq:ev1}
ma:=\underbrace{a\oplus \cdots \oplus a}_{\text{$m$ times}}.
\end{align}
Also, any positive integer $m$ is identified with an element of $\bbP$
as $m=m1\in \bbP$.
For example, in the trivial semifield $\bbone$,
we have $2=1\oplus 1=1$,
which is a little confusing.
(See Remark \ref{1rem:conf1}
for another problem of this notation.)

\begin{defn}[Semifield homomorphism]
\index{semifield!homomorphism}
For any semifields $\mathbb{P}$ and $\mathbb{P}'$,
a map $\varphi:\mathbb{P}\rightarrow \mathbb{P}'$ is called a
\emph{semifield homomorphism\/}
 if it preserves the multiplication and the addition.
 \end{defn}

\begin{ex}[Trivial homomorphisms] \index{semifield!trivial homomorphism}
For any semifield $\bbP$,
the map 
\begin{align}
\varphi_{\mathrm{triv}}:\bbP\rightarrow \bbone,
\
a\mapsto 1
\end{align}
is a semifield homomorphism.
\end{ex}

The following fact
justifies the name ``the universal semifield'' for $\mathbb{Q}_{\mathrm{sf}}(\bfu)$.

\begin{prop}
\label{1prop:uni1}
Let $\mathbb{Q}_{\mathrm{sf}}(\bfu)$ be the universal semifield
with $\bfu=(u_1,\dots, u_n)$.
Let $\bbP$ be any semifield, and let $\bfa=(a_1,\dots,a_n)$ be any 
$n$-tuple of elements in
$\bbP$.
Then, there is a unique semifield homomorphism $\pi:
\mathbb{Q}_{\mathrm{sf}}(\bfu)\rightarrow  \bbP$
such that $\pi(u_i)=a_i$ for any $i=1,\dots,n$.
\end{prop}
\begin{proof}
First,
for any nonzero polynomial $p(\bfu)$ with nonnegative integer coefficients,
the image $\pi(p(\bfu))$ is uniquely determined 
by replacing $u_i$ with $a_i$ and $+$ with $\oplus$ in $\bbP$ in the polynomial $p(\bfu)$,
where the notation \eqref{1eq:ev1} is taken into account. 
The equality $\pi(p(\bfu)q(\bfu))=\pi(p(\bfu))\pi(q(\bfu))$ is guaranteed by the 
axiom \eqref{1eq:sf1}--\eqref{1eq:sf3}.
Next, for a given  $f(\bfu)\in 
\mathbb{Q}_{\mathrm{sf}}(\bfu)$, 
take any subtraction-free expression  $f(\bfu)=p(\bfu)/q(\bfu)$.
Then, the image $\pi (f(\bfu))$ is defined by $\pi(p(\bfu))/\pi(q(\bfu))$.
To see that it is well-defined,
let us take another  subtraction-free expression $f(\bfu)=p'(\bfu)/q'(\bfu)$.
Then, we have $p(\bfu)q'(\bfu)=p'(\bfu)q(\bfu)$. Thus, $\pi(p(\bfu))\pi(q'(\bfu))=\pi(p'(\bfu))\pi(q(\bfu))$.
Therefore, we have $\pi(p(\bfu))/\pi(q(\bfu))=\pi(p'(\bfu))/\pi(q'(\bfu))$.
It is easy to show that $\pi$ is a semifield homomorphism.
\end{proof}

\begin{defn}[Specialization]
\index{semifield!specialization}
\label{1defn:sp1}
Let $ \mathbb{Q}_{\mathrm{sf}}(\bfu)$,
$\bfa\in \bbP$, and $\pi$ be
the ones in Proposition \ref{1prop:uni1}.
For $f(\bfu)\in \mathbb{Q}_{\mathrm{sf}}(\bfu)$,
the image $\pi(f(\bfu))$ 
 is called the \emph{specialization\/}
of $f(\bfu)$ at $\bfa$ in $\bbP$,
and denoted by $f\vert_{\bbP}(\bfa)$.
\end{defn}

\begin{ex}[Tropicalization homomorphism]
\index{tropicalization!homomorphism}\index{semifield!tropicalization homomorphism}
\label{1ex:trop1}
Let us
consider  the semifields $\mathbb{Q}_{\mathrm{sf}}(\bfu)$ and
$\mathrm{Trop}(\bfu)$ with common formal variables
$\bfu=(u_1,\dots, u_n)$.
By applying Proposition \ref{1prop:uni1} with 
$\bbP=\mathrm{Trop}(\bfu)$ and $\bfa=\bfu$,
we have a unique semifield homomorphism
\begin{align}
\label{1eq:trophom1}
\pi_{\mathrm{trop}}: \mathbb{Q}_{\mathrm{sf}}(\bfu)
\rightarrow \mathrm{Trop}(\bfu)
\end{align}
such that $\pi_{\mathrm{trop}}(u_i)=u_i$ for any $i=1,\dots,n$.
We call it the \emph{tropicalization homomorphism}.
For example, for $\bfu=(u_1,u_2,u_3)$,
\begin{align}
\pi_{\mathrm{trop}}
\left(
\frac{3u_1u_2^2u_3^2 + 2u_1^2u_2u_3}{3u_2^2 + u_1^2 u_2^2+
u_1u_2^3u_3}
\right)
=
\frac{u_1u_2u_3}{u_2^2}=u_1u_2^{-1}u_3.
\end{align}
Roughly speaking, it extracts the ``leading Laurent monomial'' of 
a function $f(\bfu)$ in $ \mathbb{Q}_{\mathrm{sf}}(\bfu)$.
\end{ex}

Since a semifield $\bbP$ is a multiplicative abelian group,
we can construct the group ring $\bbZ\bbP$ of $\bbP$.
The addition in $\bbZ\bbP$,
denoted by $+$ as usual,
should be distinguished from 
the addition $\oplus$ in $\bbP$.

To construct the field of fractions of the ring $\bbZ\bbP$,
 the ring  $\bbZ\bbP$ should be a domain.
This is automatically guaranteed.

\begin{prop}
\label{1prop:dom1}
For any semifield $\bbP$, the following facts hold.
\par
(a).
There is no torsion element in $\bbP$ other than 1.
Namely, for any $p\in \bbP$,
if $p^n= 1$ for some positive integer $n$,
then $p=1$.
\par
(b).
The group ring $\bbZ\bbP$ is a domain.
Namely, there is no zero divisor other than 0.
\end{prop}

\begin{proof}
(a). Suppose that $p^n=1$ for some positive integer $n$.
Then, we have
\begin{align}
p=\frac{p\oplus p^2 \oplus \cdots \oplus p^{n}}
{1\oplus p \oplus \cdots \oplus p^{n-1}}
=\frac{p\oplus p^2 \oplus \cdots \oplus 1}
{1\oplus p \oplus \cdots \oplus p^{n-1}}
=1.
\end{align}
\par
(b). This follows from the  known fact that, for any  multiplicative abelian group $G$
with no torsion element other than $1$,
its group ring $\bbZ G$ is a domain (e.g., \cite{May69}).
For completeness, let us give a proof.
Take any $a,b\in\bbZ G$.
They are written as finite sums $a=\sum_i m_i g_i$, $b=\sum_i m'_i g'_i$
($g_i, g'_i \in G$) with nonzero coefficients
$m_i, m'_i\in \bbZ$.
Let $H$ be the subgroup of $G$ generated by
all $g_i, g'_i$.
Then, $H$ is finitely generated,
and we have $a,b, ab\in \bbZ H$.
By the assumption, $H$  has  no torsion element other than $1$.
Therefore, thanks to the fundamental theorem of
finitely generated abelian groups,
we have $H\simeq \bbZ^n$ for some $n$.
Thus, $\bbZ H$ is isomorphic to the Laurent polynomial ring
$\bbZ[x_1^{\pm1}, \dots, x_n^{\pm1}]$, which is a domain. Therefore, if $ab=0$, then $a=0$ or $b=0$.
\end{proof}

Thanks to Proposition \ref{1prop:dom1},
the field of fractions of  $\bbZ\bbP$ is well-defined,
and it is denoted by $\bbQ\bbP$.

\begin{rem}
\label{1rem:conf1}
There is some notational conflict between $\bbP$ and  $\bbZ\bbP$ (and also  $\bbQ\bbP$).
Namely, for $a\in \bbP$,
$2a$, for example, may stand for $a\oplus a$ in $ \bbP$, or $a +a$
in $\bbZ\bbP$ depending on the context, and they  are distinct.
Fortunately, we do not have serious difficulty from this conflict,
because usually this notation appears together with $\oplus$ or $+$,
which clarifies the context.
Anyway, we have to be careful.
\end{rem}

\subsection{Matrix and quiver mutations}
\label{1sec:matrix1}
Below we fix a positive integer $n$,
which is called the \emph{rank} of the forthcoming seeds, cluster patterns,
cluster algebras, etc.
\index{rank} 

\begin{defn}[Skew-symmetrizable matrix]
\index{matrix!skew-symmetrizable}
\label{1defn:skew1}
An $n\times n$ integer  matrix $B=(b_{ij})_{i,j=1}^n$
is said to be \emph{skew-symmetrizable\/}
if there is a diagonal matrix $D=(d_i\delta_{ij})_{i,j=1}^n$
whose diagonal entries $d_i$ are positive rational numbers
such that $DB$ is skew-symmetric,
i.e,
\begin{align}
\label{1eq:ss1}
d_i b_{ij}= - d_jb_{ji}.
\end{align}
The matrix $D$ is called a \emph{(left) skew-symmetrizer\/} 
of $B$,
which is not unique to $B$.
\index{skew-symmetrizer}
In particular, any skew-symmetric matrix is skew-symmetrizable
with a skew-symmetrizer $D=I$.
\end{defn}

The condition \eqref{1eq:ss1} can be rephrased in the matrix notation as
\begin{align}
DB=-B^T D,
\
\text{or }
DBD^{-1}=-B^T,
\end{align}
where $B^T$ is the transpose of $B$.
By \eqref{1eq:ss1}, we have
\begin{gather}
\label{1eq:bii1}
b_{ii}=0,\\
b_{ij}=0 \Longleftrightarrow b_{ji}=0,\\
b_{ij}> 0 \Longleftrightarrow b_{ji}< 0.
\end{gather}

\begin{ex}
The following matrices exhaust all $2\times 2$ skew-symmetrizable matrices:
\begin{align}
\begin{pmatrix}
0 & 0\\
0 & 0
\end{pmatrix}
,
\quad
\begin{pmatrix}
0 & \mp b\\
\pm a & 0
\end{pmatrix}
\quad
 (a,b\in \bbZ_{> 0}).
\end{align}
In the former case, any diagonal matrix
whose diagonals are positive integers
 is
a skew-symmetrizer.
In the latter case, a skew-symmetrizer is given by
\begin{align}
D=\left(
\begin{matrix}
a & 0\\
0 & b
\end{matrix}
\right).
\end{align}
\end{ex}

For any integer $a$,
we define 
\begin{align}
[a]_+:=\max(a,0).
\end{align}
We have the following useful equalities:
\begin{align}
\label{1eq:pos1}
a&=[a]_+ - [-a]_+,
\\
\label{1eq:pos2}
a[b]_+ + [-a]_+b
&=
a[-b]_+ + [a]_+b,
\end{align}
where \eqref{1eq:pos2} follows from \eqref{1eq:pos1}.

\begin{defn}[Matrix mutation]
\index{mutation!matrix}
\label{1defn:bmut1}
For any $n\times n$ skew-symmetrizable matrix $B=(b_{ij})$
and $k\in \{1,\dots,n\}$,
we define a new matrix $B'=\mu_k(B)=(b'_{ij})$
by the following rule:
\begin{align}
\label{1eq:bmut1}
b'_{ij}&=
\begin{cases}
-b_{ij}
&
\text{$i=k$ or $j=k$,}
\\
b_{ij}
+
b_{ik} [b_{kj}]_+
+
[-b_{ik}]_+b_{kj}
&
i,j\neq k.
\end{cases}
\end{align}
The matrix $\mu_k(B)$ is called the \emph{mutation of $B$ in direction  $k$}.
\end{defn}

\begin{rem}
\label{1rem:bmutsign1}
By the identity \eqref{1eq:pos2}, the second case of \eqref{1eq:bmut1}
is also written as
\begin{align}
\label{1eq:bmutsign1}
b'_{ij}=
b_{ij}
+
b_{ik} [-b_{kj}]_+
+
[b_{ik}]_+b_{kj}
\quad
(i,j\neq k).
\end{align}
\end{rem}

\begin{rem}
\label{1rem:bmut1}
The formula \eqref{1eq:bmut1} is verbally rephrased
that $B'$ is obtained from $B$ by the following elementary transformations of matrices:
\begin{itemize}
\item
For each $j\neq k$, add 
the $k$th column  multiplied by $[b_{kj}]_+$
to the $j$th column.
Also, for each $i \neq k$,
add the $k$th row  multiplied by $[-b_{ik}]_+$
to the $i$th row.
(By Remark \ref{1rem:bmutsign1}, one can simultaneous change the sign of
$b_{kj}$ and $b_{ik}$ above.)
\item Then, multiply $-1$ to the  $k$th column and the $k$th row.
(This is well-defined because  $b_{kk}=0$.)
\end{itemize}
\end{rem}

\begin{prop}
\label{1prop:Bmut1}
 Let $B$ a skew-symmetrizable matrix,
and let $B'=\mu_k(B)$.
Then, the following facts hold:
\par
(a). Any skew-symmetrizer of $D$ is also a skew-symmetrizer of $B'$.
Therefore, $B'$ is also skew-symmetrizable.
\par
(b). We have $\det B=\det B'$.  In particular, if $B$ is nonsingular, $B'$ is also nonsingular.
\par
(c). We have $B=\mu_k(B')$. Namely, each mutation $\mu_k$ is involutive.
\end{prop}

\begin{proof}
(a). For $i=k$ or $j=k$,
we have
\begin{align}
d_i b'_{ij} = -d_i b_{ij} =d_j b_{ji}=
 - d_j b'_{ji}.
\end{align}
For $i,j\neq k$, we have
\begin{align}
\begin{split}
d_i b'_{ij}&=
d_i(
 b_{ij}
+
b_{ik} [b_{kj}]_+
+
[-b_{ik}]_+b_{kj})\\
&=
-d_j b_{ji}
-d_k b_{ki} [b_{kj}]_+
+
[d_kb_{ki}]_+b_{kj}\\
&=
-d_j(
 b_{ji}
+
b_{ki} [-b_{jk}]_+
+
[b_{ki}]_+b_{jk})
= -d_j b'_{ji}.
\end{split}
\end{align}
\par
(b). 
The verbal version of the mutation in Remark \ref{1rem:bmut1}
 preserves the determinant.
\par
(c). Let $B''=\mu_k(B')$. 
For $i=k$ or $j=k$, we have
\begin{align}
b''_{ij}=-b'_{ij}=b_{ij}.
\end{align}
For $i,j\neq k$,
\begin{align}
\begin{split}
b''_{ij}&=
b'_{ij}
+
b'_{ik} [b'_{kj}]_+
+
[-b'_{ik}]_+b'_{kj}\\
&=
(b_{ij}
+
b_{ik} [b_{kj}]_+
+
[-b_{ik}]_+b_{kj}
)
-
b_{ik} [-b_{kj}]_+
-
[b_{ik}]_+ b_{kj}
= b_{ij},
\end{split}
\end{align}
where in the last equality we used \eqref{1eq:pos1} or \eqref{1eq:pos2}.
\end{proof}

\begin{ex}
Here is an example of a matrix mutation.
\begin{align}
B=
\left(
\begin{matrix}
0 & 6 & -3\\
-12 & 0 & 6\\
2 & -2 & 0
\end{matrix}
\right),
\quad
B'=\mu_1(B)=
\left(
\begin{matrix}
0 & -6 & 3\\
12 & 0 & -30\\
-2 & 10 & 0
\end{matrix}
\right),
\
\end{align}
where
a common skew-symmetrizer of $B$ and $B'$ is given by
\begin{align}
D=
\left(
\begin{matrix}
2 & 0 & 0\\
0 & 1 & 0\\
0 & 0 & 3
\end{matrix}
\right).
\end{align}
\end{ex}

The matrix mutation is compatible with the matrix decomposition.

\begin{prop}
\label{1prop:Bdecom1}
Suppose that an $n\times n$ skew-symmetrizable matrix $B$  is decomposed as
\begin{align}
\label{1eq:Bdecom1}
B
=
\begin{pmatrix}
B_1 & O\\
O & B_2
\end{pmatrix},
\end{align}
where   $B_1$ and $B_2$ are $n_1\times n_1$ and $n_2\times n_2$ matrices, respectively.
Then, we have
\begin{align}
\label{1eq:Bdecom2}
\mu_k(B)
=
\begin{cases}
\begin{pmatrix}
\mu_k(B_1) & O\\
O & B_2
\end{pmatrix}
& k=1,\dots, n_1,\\ 
\rule{0pt}{25pt}\hskip-2pt
\begin{pmatrix}
B_1 & O\\
O & \mu_k(B_2)
\end{pmatrix}
& k=n_1+1,\dots, n.
\end{cases}
\end{align}
\end{prop}
\begin{proof}
This is clear from the verbal version of the mutation
in Remark \ref{1rem:bmut1}.
\end{proof}

Any \emph{skew-symmetric\/} matrix can be represented by a \emph{quiver}.

\begin{defn}[Quiver]
\index{quiver}
A (finite) \emph{quiver\/} 
 is a finite directed graph. Namely, it consists of a finite set of vertices
and a finite set of arrows between the vertices.
The following arrows are called a \emph{loop\/} \index{quiver!loop} and a \emph{2-cycle}, \index{quiver!2-cycle} respectively.
$$
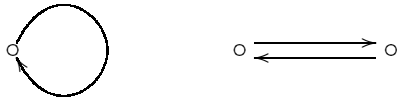
\begin{xy}
(0,0)*\cir<2pt>{},
(30,0)*\cir<2pt>{},
(50,0)*\cir<2pt>{},
(0.5,1); (0.5,-1) **\crv{(5,10) &(20,0)&(5,-10)};
\ar@{->} (0.62,-1.2);(0.5,-1)
\ar@{->} (32,1);(48,1)
\ar@{<-} (32,-1);(48,-1)
\end{xy}
$$

\end{defn}

\begin{ex} Here is an example of a quiver with a loop and a 2-cycle.
$$
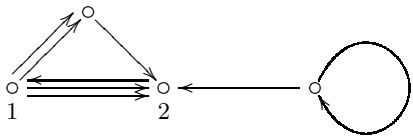
\begin{xy}
(0,0)*\cir<2pt>{},
(10,10)*\cir<2pt>{},
(20,0)*\cir<2pt>{},
(40,0)*\cir<2pt>{},
(0,-3)*{\text{\small1}},
(20,-3)*{\text{\small2}},
(40.5,1); (40.5,-1) **\crv{(45,10) &(60,0)&(45,-10)};
%**\crv{(41,1) & (60,0) & (41,-1)};
\ar@{->} (2,0);(18,0)
\ar@{->} (2,-1);(18,-1)
\ar@{->} (18,1);(2,1)
\ar@{<-} (22,0);(38,0)
\ar@{->} (1,1);(9,9)
\ar@{->} (0,1.8);(8,9.8)
\ar@{->} (11,9);(19,1)
\ar@{->} (40.62,-1.2);(40.5,-1)
\end{xy}
$$
To be more precise, there are two ways to count $2$-cycles in this quiver.
If we do not distinguish two arrows from the vertex $1$ to the vertex $2$,
there are only \emph{one} $2$-cycle.
On the other hand, if we distinguish them, there are \emph{two} $2$-cycles.
Here, we employ the former viewpoint.
Namely,  we only care about the multiplicity of arrows.
\end{ex}

For a quiver with $n$ vertices, we assume that its vertices are labeled  with 1, \dots, $n$
without duplication.
For a skew-symmetric matrix $B=(b_{ij})$, one can associate a quiver $Q(B)$ without loops and 2-cycles by the following rule:
\begin{itemize}
\item 
For each matrix element $b_{ij}>0$, we assign $b_{ij}$ arrows from the vertex $i$ to the vertex $j$.
\end{itemize}
Since $b_{ii}=0$, there are no loops. Also,  there are no 2-cycles,
because $b_{ji}=-b_{ij}<0$ if $b_{ij}>0$.

Conversely, one can  recover a skew-symmetric matrix $B$ from a quiver without loops and 2-cycles by applying the above rule in the opposite direction.
It is clear that this correspondence is one-to-one.

One can translate the matrix mutation into the quiver mutation as follows.

\begin{defn}[Quiver mutation]
\index{mutation!quiver}
\label{1defn:qmut1}
For any quiver $Q$ without loops and 2-cycles and $k=1,\dots,n$, we define a new quiver
$Q'=\mu_k(Q)$ by the following operation:
\begin{itemize}
\item 
For each  pair $i$, $j$ ($i\neq j$) such that $i,j\neq k$,
if
there are $p>0$ arrows from the vertex $i$ to the vertex $k$,
 and $q>0$ arrows from the vertex $k$ to the vertex $j$,
 then  add $pq$ arrows from the vertex $i$ to the vertex $j$.
\item Remove  the resulting 2-cycles as many as possible.
\item Invert all arrows into and out of the vertex $k$.
\end{itemize}
The quiver $Q'$ is called the \emph{mutation of $Q$ at the vertex $k$}.
\end{defn}

To see the equivalence to the matrix mutation \eqref{1eq:bmut1},
we rewrite the second case of \eqref{1eq:bmut1}
as follows: For $i,j\neq k$,
\begin{align}
\label{1eq:bqmut1}
b'_{ij}&=
\begin{cases}
b_{ij}+b_{ik}b_{kj}
&
b_{ik}, b_{kj}>0,
\\
b_{ij}-b_{ik}b_{kj}
&
b_{ik}, b_{kj}<0,
\\
b_{ij}
&
\mbox{otherwise}.
\end{cases}
\end{align}
It is easy to see that
the first two operations in Definition \ref{1defn:qmut1}
correspond to \eqref{1eq:bqmut1},
while the last operation corresponds to the first case
of \eqref{1eq:bmut1}.
\begin{ex}
\label{ex:QB1}
Here is an example of a skew-symmetric matrix $B$
and the corresponding quiver $Q(B)$,
where the number attached to each arrow represents
the multiplicity of the arrows of the same kind:
$$
B=
\left(
\begin{matrix}
0 & 3 & -2 &2\\
-3 & 0 & 4 & 0\\
2& -4 & 0 & 1\\
-2 & 0 & -1& 0
\end{matrix}
\right),
\quad
Q(B)=
\
\lower12pt
\hbox{
\begin{xy}
(0,0)*\cir<2pt>{},
(0,10)*\cir<2pt>{},
(10, 0)*\cir<2pt>{},
(10,10)*\cir<2pt>{},
(-3,13)*{\text{\small1}},
(13,13)*{\text{\small2}},
(13,-3)*{\text{\small3}},
(-3,-3)*{\text{\small4}},
(5,12)*{\text{\small3}},
(12,5)*{\text{\small4}},
(-2,5)*{\text{\small2}},
(3.5,4)*{\text{\small2}},
\ar@{->} (2,10);(8,10) %12
\ar@{->} (10,8);(10,2) %23
\ar@{->} (8,0);(2,0) %34
\ar@{->} (0,8);(0,2) %14
\ar@{->} (8,2);(2,8) %31
\end{xy}
}.
$$
The mutation of $Q(B)$ at the vertex 1 is done as follows:
$$
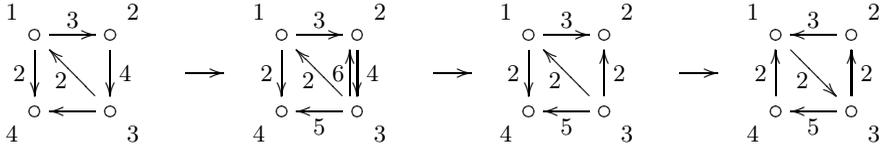
\begin{xy}
(0,0)*\cir<2pt>{},
(0,10)*\cir<2pt>{},
(10, 0)*\cir<2pt>{},
(10,10)*\cir<2pt>{},
(-3,13)*{\text{\small1}},
(13,13)*{\text{\small2}},
(13,-3)*{\text{\small3}},
(-3,-3)*{\text{\small4}},
(5,12)*{\text{\small3}},
(12,5)*{\text{\small4}},
(-2,5)*{\text{\small2}},
(3.5,4)*{\text{\small2}},
\ar@{->} (2,10);(8,10) %12
\ar@{->} (10,8);(10,2) %23
\ar@{->} (8,0);(2,0) %34
\ar@{->} (0,8);(0,2) %14
\ar@{->} (8,2);(2,8) %31
\ar@{->} (20,5);(25,5) %31
\end{xy}
\quad
\begin{xy}
(0,0)*\cir<2pt>{},
(0,10)*\cir<2pt>{},
(10, 0)*\cir<2pt>{},
(10,10)*\cir<2pt>{},
(-3,13)*{\text{\small1}},
(13,13)*{\text{\small2}},
(13,-3)*{\text{\small3}},
(-3,-3)*{\text{\small4}},
(5,12)*{\text{\small3}},
(7.5,5)*{\text{\small6}},
(12,5)*{\text{\small4}},
(-2,5)*{\text{\small2}},
(3.5,4)*{\text{\small2}},
(5,-2)*{\text{\small5}},
\ar@{->} (2,10);(8,10) %12
\ar@{->} (10,8);(10,2) %23
\ar@{<-} (9,8);(9,2) %32
\ar@{->} (8,-0);(2,-0) %34
\ar@{->} (0,8);(0,2) %14
\ar@{->} (8,2);(2,8) %31
\ar@{->} (20,5);(25,5) %31
\end{xy}
\quad
\begin{xy}
(0,0)*\cir<2pt>{},
(0,10)*\cir<2pt>{},
(10, 0)*\cir<2pt>{},
(10,10)*\cir<2pt>{},
(-3,13)*{\text{\small1}},
(13,13)*{\text{\small2}},
(13,-3)*{\text{\small3}},
(-3,-3)*{\text{\small4}},
(5,12)*{\text{\small3}},
(12,5)*{\text{\small2}},
(-2,5)*{\text{\small2}},
(3.5,4)*{\text{\small2}},
(5,-2)*{\text{\small5}},
\ar@{->} (2,10);(8,10) %12
\ar@{<-} (10,8);(10,2) %23
\ar@{->} (8,0);(2,0) %34
\ar@{->} (0,8);(0,2) %14
\ar@{->} (8,2);(2,8) %31
\ar@{->} (20,5);(25,5) %31
\end{xy}
\quad
\begin{xy}
(0,0)*\cir<2pt>{},
(0,10)*\cir<2pt>{},
(10, 0)*\cir<2pt>{},
(10,10)*\cir<2pt>{},
(-3,13)*{\text{\small1}},
(13,13)*{\text{\small2}},
(13,-3)*{\text{\small3}},
(-3,-3)*{\text{\small4}},
(5,12)*{\text{\small3}},
(12,5)*{\text{\small2}},
(-2,5)*{\text{\small2}},
(3.5,4)*{\text{\small2}},
(5,-2)*{\text{\small5}},
\ar@{<-} (2,10);(8,10) %12
\ar@{<-} (10,8);(10,2) %23
%\ar@{<-} (19,8);(19,-8) %32
\ar@{->} (8,0);(2,0) %34
\ar@{<-} (0,8);(0,2) %14
\ar@{<-} (8,2);(2,8) %31
\end{xy}
$$
\smallskip
Putting it back to the matrix form, we obtain
\begin{align}
\mu_1(B)=
\left(
\begin{matrix}
0 & -3 & 2 &-2\\
3 & 0 & -2 & 0\\
-2& 2 & 0 & 5\\
2 & 0 & -5& 0
\end{matrix}
\right).
\end{align}
\end{ex}

The above correspondence between skew-symmetric matrices and quivers
is the basis of the connection between cluster algebras and quiver representations.
See \cite{Plamondon16} for a concise survey on the subject.

\newpage

\section{Basic Notions}

In this section we introduce basic notions for cluster algebras.

\subsection{Seeds and mutations}

Let us introduce the most fundamental notions in cluster algebra theory,
namely, \emph{seeds\/} and \emph{mutations}.

Recall that, for a given semifield $\bbP$, $\bbQ\bbP$
is the  field of fractions of the group algebra $\bbZ\bbP$
in Section \ref{1sec:semifields1}.

\begin{defn}[Seed/Cluster variable/Coefficient]
\label{1defn:seed1}
Let $n$ be any positive integer.
Let $\bbP$ be any semifield,
and let
$\calF$  be a field that is isomorphic to
the rational function field of $n$-variables
with coefficients in the field $\bbQ\bbP$.
\begin{itemize}
\item
A \emph{(labeled) seed 
\index{seed}\index{seed!labeled} 
with coefficients in $\bbP$
(or a seed in $\calF$) of rank $n$}
is a triplet
$\Sigma=(\bfx,\bfy,B)$
such that
$\bfx=(x_1,\dots,x_n)$ is an
$n$-tuple of algebraically independent and generating elements
in $\calF$
 (i.e., a transcendence basis of $\calF$),
$\bfy=(y_1,\dots,y_n)$ is an
$n$-tuple of any  elements
in $\bbP$,
and $B=(b_{ij})_{i,j=1}^n$
is an $n\times n$ skew-symmetrizable (integer) matrix.
\item
We call $\bfx$, $\bfy$, and $B$,
respectively, 
the \emph{cluster}, \index{cluster}
the \emph{coefficient tuple}, \index{coefficient!tuple}
and the \emph{exchange matrix\/} \index{matrix!exchange} of $\Sigma$. 
The elements $x_i$ and $y_i$, respectively, are called
the \emph{cluster variables\/} \index{cluster variable} and the \emph{coefficients}. \index{coefficient}
\item
We  call $\bbP$ and $\calF$, respectively,  the \emph{coefficient semifield\/} \index{semifield!coefficient}
and  the \emph{ambient field\/} \index{ambient field} of a seed $\Sigma$,
and also, of the forthcoming cluster patterns, cluster algebras, etc.
\end{itemize}
\end{defn}

\begin{rem}
Sometimes,  cluster variables $x_i$ and coefficients $y_i$ are
casually  called \emph{$x$-variables} and \emph{$y$-variables}, respectively.
Or, we may simply say \emph{variable $x_i$}, \emph{variable $y_i$}, etc.
Alternatively, they are also denoted by $A_i$ and $X_i$, and called \emph{$A$-coordinates ($A$-variables)} and \emph{$X$-coordinates ($X$-variables)}, respectively,
following the convention of another pioneering works on cluster algebra theory by Fock and Goncharov \cite{Fock03,Fock07}. 
\end{rem}

There are some  related notions to Definition 
\ref{1defn:seed1}.
\begin{defn}[Seed without coefficients/$Y$-seed]
\ % necessary for noindent of Def
\begin{itemize}
\item
When the coefficient semifield $\bbP$ is taken to be the trivial semifield $\bbone$,
all coefficients $y_i$ are 1.
Then, we can reduce a triplet $(\bfx,\bfy,B)$
in  Definition 
\ref{1defn:seed1}
to  a pair $(\bfx,B)$, which is called
 a \emph{seed without coefficients}. \index{seed!without coefficients}
\item 
For any coefficient semifield $\bbP$,
we have an option to ignore cluster variables
and concentrate on a pair $Y=(\bfy,B)$
in Definition 
\ref{1defn:seed1},
which is called a \emph{$Y$-seed\/}  in $\bbP$. \index{seed!$Y$-}\index{$Y$-seed}
\end{itemize}
\end{defn}
For any seed $\Sigma=(\bfx,\bfy,B)$,
we attach an $n$-tuple 
$\hat\bfy=(\hat{y}_1,\dots,\hat{y}_n)$
of elements in $\calF$,
\begin{align}
\label{1eq:yhat1}
\hat{y}_i
=y_i \prod_{j=1}^n x_j^{b_{ji}}.
\end{align}
They play an important role in the cluster algebra theory.
We call them \emph{$\hat{y}$-variables}. \index{$\hat{y}$-variable}

\begin{defn}[Seed mutation]
\index{mutation!seed}
\label{1defn:smut1}
For any seed $\Sigma=(\bfx,\bfy,B)$ in $\calF$
and $k\in \{1,\dots,n\}$,
we define a new seed 
$\Sigma'=(\bfx',\bfy',B')$ in $\calF$
by the following rule:
\begin{align}
\label{1eq:xmut1}
x'_i
&=
\begin{cases}
\displaystyle
x_k^{-1}
\Biggl(\, \prod_{j=1}^n x_j^{[-b_{jk}]_+}
\Biggr)
\frac{ 1+\hat{y}_k}{ 1\oplus y_k}
& i=k,
\\
x_i
&i\neq k,
\end{cases}
\\
\label{1eq:ymut1}
y'_i
&=
\begin{cases}
\displaystyle
y_k^{-1}
& i=k,
\\
y_i y_k^{[b_{ki}]_+} (1\oplus y_k)^{-b_{ki}}
&i\neq k,
\end{cases}
\\
\label{1eq:bmut2}
b'_{ij}&=
\begin{cases}
-b_{ij}
&
\text{$i=k$ or $j=k$,}
\\
b_{ij}
+
b_{ik} [b_{kj}]_+
+
[-b_{ik}]_+b_{kj}
&
i,j\neq k,
\end{cases}
\end{align}
where $\hat{y}_k$ in \eqref{1eq:xmut1}
is defined by \eqref{1eq:yhat1}.
The seed $\Sigma'$ is called the \emph{mutation of $\Sigma$ in direction $k$},
and denoted by $\mu_k(\Sigma)=\mu_k(\bfx,\bfy,B)$.
The mutations of a seed without coefficients $(\bfx,B)$
and a $Y$-seed $(\bfx,B)$ are also defined by the same formulas.
\end{defn}

We will soon show that $\Sigma'$ in the above is indeed a seed in $\calF$.

The mutation of $B$ in \eqref{1eq:bmut2} is the matrix mutation already introduced
in Definition \ref{1defn:bmut1}.
Thanks to the identity \eqref{1eq:pos1},
the first case of \eqref{1eq:xmut1} is also written as the following
more standard form:
\begin{align}
\label{1eq:xmutstandard1}
x'_k =
x_k^{-1}
\Biggl(
\frac{ y_k}{ 1\oplus y_k}
\prod_{j=1}^n
x_j^{[b_{jk}]_+}
+
\frac{ 1}{ 1\oplus y_k}
\prod_{j=1}^n
x_j^{[-b_{jk}]_+}
\Biggr).
\end{align}

Related with the mutation \eqref{1eq:xmut1},
we note the following useful identity:
\begin{align}
\label{1eq:yy1}
\frac{ 1+\haty_k^{-1}}{ 1\oplus y_k^{-1}}
=
\frac{1+ \haty_k}{ 1\oplus y_k}
\prod_{j=1}^n x_{j}^{-b_{jk}}.
\end{align}

The following fact,  together with Proposition \ref{1prop:Bmut1},
 ensures that $(\bfx',\bfy',B')$ is indeed a seed in  $\calF$.

\begin{prop}
\label{1prop:smut1}
\par
(a). The mutation $\mu_k$  is involutive.
Namely, $\mu_k(\Sigma')=\Sigma$ for $\Sigma'=\mu_k(\Sigma)$.
\par
(b). The elements in $\bfx'$ are algebraically independent
and generating elements in $\calF$.
\end{prop}
\begin{proof}
(a). Let $(\bfx'',\bfy'',B'')=\mu_k(\bfx',\bfy',B')$ for
$(\bfx',\bfy',B')$ in Definition \ref{1defn:smut1}.
$B''=B$ was already shown in Proposition \ref{1prop:Bmut1} (b).
To show $\bfx''=\bfx$, it is enough to show that $x''_k=x_k$.
We first note that
\begin{align}
\label{1eq:yhatinv1}
\hat y'_k =
 y_k' \prod_{j=1}^n x'_{j}{}^{b'_{jk}}
=
 y_k^{-1} \prod_{j=1}^n x_{j}{}^{-b_{jk}}
 =
 \hat y_k^{-1},
\end{align}
where we used the fact $b_{kk}=0$.
Then, we have
\begin{align}
\begin{split}
x''_k &=
x'_k{}^{-1}
\Biggl(\, \prod_{j=1}^n x_j{^{[-b'_{jk}]_+}}
\Biggr)
\frac{ 1+\hat{y}'_k}{ 1\oplus y'_k}
=
x'_k{}^{-1}
\Biggl(\, \prod_{j=1}^n x_j{^{[b_{jk}]_+}}
\Biggr)
\frac{ 1+\hat{y}_k^{-1}}{ 1\oplus y_k^{-1}}
\\
&=
x'_k{}^{-1}
\Biggl(\, \prod_{j=1}^n x_j{^{[-b_{jk}]_+}}
\Biggr)
\frac{ 1+\hat{y}_k}{ 1\oplus y_k}=x_k,
\end{split}
\end{align}
where we used \eqref{1eq:yy1}
and  \eqref{1eq:pos1}.
Let us show $\bfy''=\bfy$.
We have $y''_k=y'_k{}^{-1}=y_k$.
For $i\neq k$,
\begin{align}
\begin{split}
y''_i &=
y'_i y'_k{}^{[b'_{ki}]_+} (1\oplus y'_k)^{-b'_{ki}}
\\
&=
(y_i y_k^{[b_{ki}]_+} (1\oplus y_k)^{-b_{ki}})
 y_k^{-[-b_{ki}]_+} (1\oplus y_k^{-1})^{b_{ki}}
=
y_i.
\end{split}
\end{align}
\par
(b). By (a), $x_1$, \dots, $x_n$ are expressed as rational functions of
$x'_1$, \dots, $x'_n$ over $\bbQ\bbP$.
It follows that $x'_1$, \dots, $x'_n$ generate $\calF$; moreover,
they are algebraically independent, because, 
if not, the transcendence degree of $\calF$
over $\bbQ\bbP$ becomes less than $n$,
which is a contradiction.
\end{proof}

The following fact is the first manifestation of the close relationship (duality)
between the mutations of $\bfx$ and $\bfy$.
\begin{prop}
\label{1prop:yhat1}
The $\hat{y}$-variables in \eqref{1eq:yhat1} mutate
in the ambient field $\calF$ as
\begin{align}
\hat{y}'_i
=
\begin{cases}
\displaystyle
\hat{y}_k^{-1}
& i=k,
\\
\hat{y}_i \hat{y}_k^{[b_{ki}]_+} (1+ \hat{y}_k)^{-b_{ki}}
&i\neq k,
\end{cases}
\end{align}
which is the same rule for the coefficients in \eqref{1eq:ymut1}.
\end{prop}
\begin{proof}
Below (and elsewhere) we use the fact $b_{kk}=0$ effectively.
For $i=k$,
it was already shown in \eqref{1eq:yhatinv1}.
For $i\neq k$,
\begin{align}
\begin{split}
\hat{y}'_i
&=y'_i
\prod_{j=1}^n x'_j{}^{b'_{ji}}\\
&=
y_i y_k^{[b_{ki}]_+} (1\oplus y_k)^{-b_{ki}}
\Biggl(
\prod_{\scriptstyle j=1\atop \scriptstyle  j\neq k}^n x_j^{b_{ji}+
b_{jk} [b_{ki}]_+
+
[-b_{jk}]_+b_{ki}}
\Biggr)\\
&\qquad
\times\Biggl(
x_k^{-1}
\Biggl(
\prod_{j=1}^n x_j^{[-b_{jk}]_+}
\Biggr)
\frac{ 1+\hat{y}_k}{ 1\oplus y_k}
\Biggr)^{-b_{ki}}
\\
&=\hat{y}_i\hat{y}_k^{[b_{ki}]_+}
(1+\hat{y}_k)^{-b_{ki}}.
\end{split}
\end{align}
\end{proof}

There is some flexibility to write the mutation formulas
in   Definition \ref{1defn:smut1}
and Proposition \ref{1prop:yhat1}.

\begin{prop}[$\varepsilon$-expressions]
\index{$\varepsilon$-expression}
\label{1prop:epsilon1}
The right hand sides of the following formulas do not
depend on the choice of $\varepsilon\in \{1, -1\}$.
\begin{align}
\label{1eq:xmut2}
x'_i
&=
\begin{cases}
\displaystyle
x_k^{-1}\Biggl(\, \prod_{j=1}^n x_j^{[-\varepsilon b_{jk}]_+}
\Biggr)
\frac{ 1+\hat{y}_k^{\varepsilon}}{ 1\oplus y_k^{\varepsilon}}
& i=k,
\\
x_i
&i\neq k,
\end{cases}
\\
\label{1eq:ymut2}
y'_i
&=
\begin{cases}
\displaystyle
y_k^{-1}
& i=k,
\\
y_i y_k^{[{\varepsilon}b_{ki}]_+} (1\oplus y_k^{\varepsilon})^{-b_{ki}}
&i\neq k,
\end{cases}
\\
\label{1eq:bmut3}
b'_{ij}&=
\begin{cases}
-b_{ij}
&
\text{$i=k$ {\rm or} $j=k$,}
\\
b_{ij}
+
b_{ik} [\varepsilon b_{kj}]_+
+
[-\varepsilon b_{ik}]_+b_{kj}
&
i,j\neq k,
\end{cases}
\\
\label{1eq:yhatmut2}
\haty'_i
&=
\begin{cases}
\displaystyle
\haty_k^{-1}
& i=k,
\\
\haty_i \haty_k^{[{\varepsilon}b_{ki}]_+} (1+ \haty_k^{\varepsilon})^{-b_{ki}}
&i\neq k.
\end{cases}
\end{align}
\end{prop}
\begin{proof}
For the matrix mutation, this was already pointed out
in Remark \ref{1rem:bmutsign1}.
The other cases can be also shown
by \eqref{1eq:pos1} as follows:
\begin{align}
\Biggl(\,
\prod_{j=1}^n
\frac
{x_j^{[-b_{jk}]+}
}
{x_j^{[b_{jk}]+}
}
\Biggr)
\frac{ 1+\hat{y}_k}{ 1+\hat{y}_k^{-1}}
\frac{ 1\oplus y_k^{-1}}{ 1\oplus y_k}
&=
\Biggl(\,
\prod_{j=1}^n x_j^{-b_{jk}}
\Biggr)
\haty_k y_k^{-1}=1,
\\
\frac
{
y_k^{[b_{ki}]_+}
}
{
y_k^{[-b_{ki}]_+}
}
\frac
{(1\oplus y_k)^{-b_{ki}}
}
{
 (1\oplus y_k^{-1})^{-b_{ki}}
}
&=
y_k^{b_{ki}}y_k^{-b_{ki}}=1.
\end{align}
The case \eqref{1eq:yhatmut2} is similar.
\end{proof}

We continue to give related notions for seeds.
\begin{defn}[$S_n$-action/Unlabeled seed]
\
\begin{itemize}
\item
\index{$S_n$-action (on seeds)}
For a seed $\Sigma=(\bfx, \bfy, B)$ in  $\calF$ and a permutation $\sigma$ of $\{1,\dots,n\}$,
we define the action of $\sigma$ on $\Sigma$ by
\begin{align}
\label{1eq:saction1}
\sigma\Sigma = (\sigma \bfx, \sigma \bfy, \sigma B),
\end{align}
where
$\sigma \bfx=\bfx'$, 
$\sigma \bfy=\bfy'$, 
$\sigma B=B'$ are defined by
\begin{align}
x'_i=x_{\sigma^{-1}(i)},
\quad
y'_i=y_{\sigma^{-1}(i)},
\quad
b'_{ij}=b_{\sigma^{-1}(i)\sigma^{-1}(j)}.
\quad
\end{align}
Clearly, $\sigma\Sigma$ is a seed in $\calF$,
and this yields a left action of the symmetric group $S_n$ of degree $n$
on the set of seeds in  $\calF$; namely,
we have $\tau (\sigma \Sigma)=\tau\sigma (\Sigma)$
for $\sigma, \tau\in S_n$.
This also induces the action of $\sigma$ on $\haty$-variables 
$\sigma \hat\bfy=\hat\bfy'$ as
\begin{align}
\haty'_i :=
y'_i
\prod_{j=1}^n x'_{j}{}^{b'_{ji}}
=
y_{\sigma^{-1}(i)}
\prod_{j=1}^n x_{\sigma^{-1}(j)}^{b_{\sigma^{-1}(j)\sigma^{-1}(i)}}
=\haty_{\sigma^{-1}(i)}.
\end{align}
\item
We introduce an  equivalence condition for (labeled) seeds in $\calF$,
 \begin{align}
 \Sigma'\sim\Sigma
 \end{align}
 if there is some permutation $\sigma\in S_n$
 such that $\Sigma'=\sigma\Sigma$.
 Then, each equivalence class $ [(\bfx,\bfy,B)]$ is called an \emph{unlabeled seed}
 in $\calF$.
 \index{seed!unlabeled}
 \end{itemize}
\end{defn}

The mutations and the action of $\sigma$ is compatible in the following sense.

\begin{prop}
\label{1prop:compat1}
The following equality holds:
\begin{align}
\mu_{\sigma(k)}(\sigma \Sigma)=\sigma (\mu_k(\Sigma)).
\end{align}
\end{prop}
\begin{proof}
We set the left and right hand sides as $(\bfx', \bfy',B')$ and $(\bfx'', \bfy'',B'')$,
respectively.
We calculate $(\bfx', \bfy',B')$ below, which
turns out to coincide with $(\bfx'', \bfy'',B'')$.
\begin{align}
\label{1eq:xmut4}
x'_i
&=
\begin{cases}
\displaystyle
x_k^{-1}
\Biggl(\, \prod_{j=1}^n x_{\sigma^{-1}(j)}^{[-b_{\sigma^{-1}(j)k}]_+}
\Biggr)
\frac{ 1+\hat{y}_k}{ 1\oplus y_k}
& i=\sigma(k),
\\
x_{\sigma^{-1}(i)}
&i\neq \sigma(k),
\end{cases}
\\
\label{1eq:ymut4}
y'_i
&=
\begin{cases}
\displaystyle
y_k^{-1}
& i=\sigma(k),
\\
y_{\sigma^{-1}(i)} y_k^{[b_{k{\sigma^{-1}(i)}}]_+} (1\oplus y_k)^{-b_{k{\sigma^{-1}(i)}}}
&i\neq \sigma(k),
\end{cases}
\\
\label{1eq:bmut4}
b'_{ij}
&=
\begin{cases}
-b_{{\sigma^{-1}(i)}{\sigma^{-1}(j)}}
&
\text{$i=\sigma(k)$ or $j=\sigma(k)$,}
\\
b_{{\sigma^{-1}(i)}{\sigma^{-1}(j)}}
+
b_{{\sigma^{-1}(i)}k} [b_{k{\sigma^{-1}(j)}}]_+
& i,j\neq \sigma(k).
\\
+[-b_{{\sigma^{-1}(i)}k}]_+b_{k{\sigma^{-1}(j)}}
&
\end{cases}
\end{align}
\end{proof}

Two mutations are not commutative, in general.
However, under some simple condition, they are commutative.
\begin{prop}
\label{1prop:bcommut1}
For a seed $\Sigma=(\bfx,\bfy,B)$ and 
a pair $k,\ell$ ($k\neq \ell$),
suppose that 
\begin{align}
\label{1eq:bkl1}
b_{k \ell}=b_{ \ell k}=0
\end{align}
holds.
Then, we have
\begin{align}
\mu_k\mu_{\ell}(\Sigma)
=
\mu_{\ell}\mu_{k}(\Sigma),
\end{align}
or equivalently,
\begin{align}
\label{1eq:A1A12}
 \mu_{\ell}\mu_k\mu_{\ell}\mu_k(\Sigma)=
\Sigma.
\end{align}
\end{prop}
\begin{proof}
We set $\Sigma'=\mu_{k}(\Sigma)$ and $\Sigma''=\mu_{\ell}(\Sigma')$,
and we show that $\Sigma''$ is symmetric with respect to $k$ and $\ell$.
By \eqref{1eq:bkl1}, we have
\begin{align}
%\label{1eq:xmut1}
x'_i
&=
\begin{cases}
\displaystyle
x_k^{-1}
\Biggl(\, \prod_{j=1}^n x_j^{[-b_{jk}]_+}
\Biggr)
\frac{ 1+\hat{y}_k}{ 1\oplus y_k}
& i=k,
\\
x_{i}
&i\neq k,
\end{cases}
\\
%\label{1eq:ymut1}
y'_i
&=
\begin{cases}
\displaystyle
y_k^{-1}
& i=k,
\\
y_{\ell}
& i = \ell,
\\
y_i y_k^{[b_{ki}]_+} (1\oplus y_k)^{-b_{ki}}
&i\neq k,\ell,
\end{cases}
\\
%\label{1eq:bmut2}
b'_{ij}&=
\begin{cases}
-b_{ij}
&
\text{$i=k$ or $j=k$,}
\\
b_{ij}
&
\text{$i=\ell$ or $j=\ell$,}
\\
b_{ij}
+
b_{ik} [b_{kj}]_+
+
[-b_{ik}]_+b_{kj}
&
\mbox{otherwise}.
\end{cases}
\end{align}
These also imply that
\begin{align}
\haty'_{\ell}=\haty_{\ell}.
\end{align}
Then, we have
\begin{align}
%\label{1eq:xmut1}
x''_i
&=
\begin{cases}
\displaystyle
x_k^{-1}
\Biggl(\, \prod_{j=1}^n x_j^{[-b_{jk}]_+}
\Biggr)
\frac{ 1+\hat{y}_k}{ 1\oplus y_k}
& i=k,
\\
\displaystyle
x_{\ell}^{-1}
\Biggl(\, \prod_{j=1}^n x_j^{[-b_{j \ell}]_+}
\Biggr)
\frac{ 1+\hat{y}_{\ell}}{ 1\oplus y_{\ell}}
& i=\ell,
\\
x_{i}
&i\neq k, \ell,
\end{cases}
\\
%\label{1eq:ymut1}
y''_i
&=
\begin{cases}
\displaystyle
y_k^{-1}
& i=k,
\\
y_{\ell}^{-1}
& i = \ell,
\\
y_i y_k^{[b_{ki}]_+}y_{\ell}^{[b_{\ell i}]_+} (1\oplus y_k)^{-b_{ki}}(1\oplus y_{\ell})^{-b_{\ell i}}
&i\neq k, \ell,
\end{cases}
\\
%\label{1eq:bmut2}
b''_{ij}&=
\begin{cases}
-b_{ij}
&
\text{$i=k$ or $j=k$,}
\\
-b_{ij}
&
\text{$i=\ell$ or $j=\ell$,}
\\
b_{ij}
+
b_{ik} [b_{kj}]_+
+
[-b_{ik}]_+b_{kj}
&
\mbox{otherwise}.
\\
\phantom{b_{ij}}
+ b_{i\ell} [b_{\ell j}]_+
+
[-b_{i\ell}]_+b_{\ell j}
\end{cases}
\end{align}
\end{proof}

\subsection{Cluster patterns and cluster algebras}

We introduce  cluster patterns and cluster algebras,
which are the main objects to be studied
in cluster algebra theory.

\pagebreak[2]
\begin{defn}[$n$-Regular tree]
\index{regular tree}
\
\begin{itemize}
\item
Let $\bbT_n$ denote the \emph{$n$-regular tree};
namely, it is a tree graph such that
each vertex has exactly $n$ edges attached to it.
Moreover,  the edges are labeled by $1$,\dots, $n$
so that the edges attached to each vertex 
are labeled without duplication.
By abusing the notation,
the set of vertices of $\bbT_n$ is also denoted by 
 $\bbT_n$.
 \item
We say that a pair of vertices $t$ and $t'$
in $\bbT_n$ are \emph{$k$-adjacent}, or $t'$ is \emph{$k$-adjacent to $t$},
\index{$k$-adjacent}
if they are connected with an edge labeled by $k$.
 \end{itemize}
\end{defn}
The graph $\bbT_n$ is finite for $n=1$ and infinite otherwise.

 \begin{ex}
 Here are the $n$-regular  trees $\bbT_n$
 for $n=1, 2, 3$.
 \smallskip
 $$
 
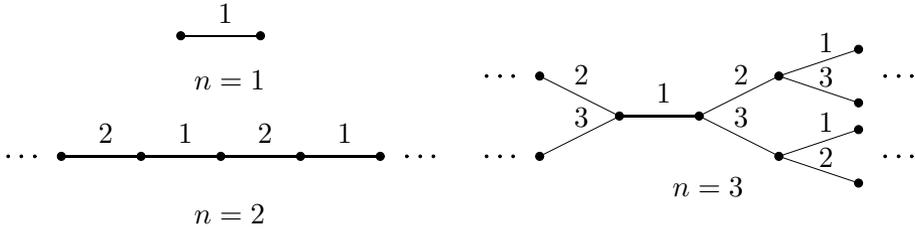
\begin{picture}(200,90)(40,5)
\put(100,40){\circle*{3}}
\put(70,40){\circle*{3}}
\put(40,40){\circle*{3}}
\put(130,40){\circle*{3}}
\put(160,40){\circle*{3}}
\put(40,40){\line(1,0){120}}
\put(20,40){\circle*{1}}
\put(25,40){\circle*{1}}
\put(30,40){\circle*{1}}
\put(170,40){\circle*{1}}
\put(175,40){\circle*{1}}
\put(180,40){\circle*{1}}
\put(90,15){$n=2$}
%\put(68,50){$t_0$}
\put(54,45){$2$}
\put(84,45){$1$}
\put(114,45){$2$}
\put(144,45){$1$}
\put(85,85){\circle*{3}}
\put(115,85){\circle*{3}}
\put(85,85){\line(1,0){30}}
\put(90,65){$n=1$}
%\put(83,95){$t_0$}
\put(99,90){$1$}
\end{picture}
\begin{picture}(100,40)(60,-10)
\put(40,25){\circle*{3}}
\put(40,55){\circle*{3}}
\put(70,40){\circle*{3}}
\put(100,40){\circle*{3}}
\put(130,55){\circle*{3}}
\put(130,25){\circle*{3}}
\put(160,65){\circle*{3}}
\put(160,45){\circle*{3}}
\put(160,35){\circle*{3}}
\put(160,15){\circle*{3}}
\put(40,25){\line(2,1){30}}
\put(40,55){\line(2,-1){30}}
\put(70,40){\line(1,0){30}}
\put(100,40){\line(2,1){30}}
\put(100,40){\line(2,-1){30}}
\put(130,55){\line(3,1){30}}
\put(130,55){\line(3,-1){30}}
\put(130,25){\line(3,1){30}}
\put(130,25){\line(3,-1){30}}
\put(20,25){\circle*{1}}
\put(25,25){\circle*{1}}
\put(30,25){\circle*{1}}
\put(20,55){\circle*{1}}
\put(25,55){\circle*{1}}
\put(30,55){\circle*{1}}
\put(170,25){\circle*{1}}
\put(175,25){\circle*{1}}
\put(180,25){\circle*{1}}
\put(170,55){\circle*{1}}
\put(175,55){\circle*{1}}
\put(180,55){\circle*{1}}
\put(90,10){$n=3$}
%\put(68,50){$t_0$}
\put(53,52){$2$}
\put(53,36){$3$}
\put(84,45){$1$}
\put(113,52){$2$}
\put(113,36){$3$}
\put(145,64){$1$}
\put(145,51){$3$}
\put(145,34){$1$}
\put(145,21){$2$}
\end{picture}
$$
 \end{ex}

\begin{defn}[Cluster pattern/$Y$-pattern/$B$-pattern]
\ \par
\begin{itemize}
\item
A collection of seeds $\bfSigma=\{ \Sigma_t
=(\bfx_t,\bfy_t,B_t)
\}_{ t\in \bbT_n} $ with coefficients in $\bbP$
 indexed by $\bbT_n$
is called a \emph{cluster pattern with coefficients in $\bbP$}
 if,
for any pair $t,t'\in \bbT_n$ that are $k$-adjacent,
the equality $\Sigma_{t'}=\mu_k( \Sigma_t)$ holds.
A cluster pattern is also called  a \emph{seed pattern}. 
 \index{cluster pattern}\index{pattern!cluster}
 \index{pattern!seed|see{cluster ---}}
\item
We replace the above $\bfSigma$
with 
a collection of seeds $\bfSigma=\{ \Sigma_t
=(\bfx_t,B_t)
\}_{ t\in \bbT_n}$ without coefficients.
Then,
it is called a \emph{cluster pattern without coefficients.} \index{cluster pattern!without coefficients}

\item
We replace the above $\bfSigma$
with 
a collection of $Y$-seeds $\bfY=\{ Y_t
=(\bfy_t,B_t)
\}_{ t\in \bbT_n}$  in $\bbP$.
Then, it is called a \emph{$Y$-pattern in $\bbP$}.  \index{pattern!$Y$-}\index{$Y$-pattern}

\item
From the above $\bfSigma$, we extract
a collection of the exchange matrices $\bfB=\{ B_t
\}_{ t\in \bbT_n}$.
We call it the \emph{$B$-pattern of $\bfSigma$}.  \index{pattern!$B$-}\index{$B$-pattern}
\end{itemize}
\end{defn}

Often, it is convenient to
 choose arbitrarily a distinguished
 vertex $t_0\in \bbT_n$ called the \emph{initial vertex}. \index{initial vertex}
 Since any cluster pattern
 $\bfSigma=\{ \Sigma_t
=(\bfx_t,\bfy_t,B_t)
\mid t\in \bbT_n 
\}$ is uniquely determined from
the \emph{initial seed\/} \index{seed!initial}
$\Sigma_{t_0}$ at $t_0$
by repeating mutations,
we  may write  $\bfSigma=\bfSigma(\Sigma_{t_0})$.

For  a seed $\Sigma_t=(\bfx_t,\bfy_t,B_t)$ 
in a cluster pattern $\bfSigma$,
we use the notation
\begin{align}
\bfx_t=(x_{1;t},\dots,x_{n;t}),
\quad
\bfy_t=(y_{1;t},\dots,y_{n;t}),
\quad
B_t=(b_{ij;t})_{i,j=1}^n.
\end{align}
Often, we omit the index $t_0$ for the initial seed  as
\begin{align}
\label{1eq:init1}
\bfx_{t_0}=\bfx=(x_{1},\dots,x_{n}),
\quad
\bfy_{t_0}=\bfy=(y_{1},\dots,y_{n}),
\quad
B_{t_0}=B=(b_{ij})_{i,j=1}^n.
\end{align}
We use  similar notations for $\haty$-variables
as well:
\begin{align}
 \hat{\bfy}_t=(\haty_{1;t},\dots,\haty_{n;t}),
 \quad
\hat\bfy_{t_0}=\hat\bfy=(\haty_{1},\dots,\haty_{n}).
\end{align}

By fixing the initial vertex,
 one may regard each coefficient $y_{i;t}$ as a rational function
with a subtraction-free expression
 in the initial coefficients
$\bfy$,
 since the mutations of coefficients in \eqref{1eq:ymut1}
 are subtraction-free.
Similarly,
 one may regard each cluster variable $x_{i;t}$  as a rational function
with a subtraction-free expression
 in the initial cluster variables
 $\bfx$;
 moreover, its coefficients are
 rational functions
with subtraction-free expressions
 in $\bfy$.

Now we  give the definition of a cluster algebra.

\begin{defn}[Cluster algebra]
 \index{cluster algebra}
For any cluster pattern $\bfSigma$,
the \emph{cluster algebra\/} $\calA=\calA(\bfSigma)$
associated with $\bfSigma$
is the $\bbZ\bbP$-subalgebra of the ambient field $\calF$
generated by all cluster variables $x_{i;t}$ ($i=1,\dots,n; t\in \bbT_n$)
of $\bfSigma$.
\end{defn}

If $\bfSigma$ has only finitely many distinct cluster variables,   
$\calA(\bfSigma)$ is clearly finitely generated.
On the other hand, if $\bfSigma$ has  infinitely many distinct cluster variables,
$\calA(\bfSigma)$ may be finitely generated or not, depending on $\bfSigma$.

\begin{ex}[Type $A_1$]
Let $n=1$.
We consider the following arrangement of
a cluster pattern $\bfSigma$ on $\bbT_1$.
$$
 \begin{xy}
(0,0)*{\bullet},
(10,0)*{\bullet},
(5,-3)*{1},
(0,-3)*{t_0},
(10,-3)*{t_1},
(1,3)*{\Sigma_{t_{0}}},
(11,3)*{\Sigma_{t_{1}}},
\ar@{-} (0,0);(10,0)
\end{xy}
$$
We have the unique choice of the exchange matrices
\begin{align}
B_{t_0}=B_{t_1}=(0).
\end{align}
Let $t_0$ be the initial vertex, and we set $x_{1;t_0}=x_1$ and $y_{1;t_0}=y_1$.
We have $\haty_{1;t_0}=\haty_1=y_1$.
Accordingly, we have
\begin{align}
x_{1;t_1}&=x_{1}^{-1}\frac{1+ y_1}{1\oplus y_1},
\quad
y_{1;t_1}=y_{1}^{-1}.
\end{align}
One can directly confirm the involution property $\mu_1(\Sigma_{t_1})=\Sigma_{t_0}$ as
follows:
\begin{align}
\label{1eq:a11}
x_{1;t_1}^{-1} \frac{1+ y_1^{-1}}{1\oplus y_1^{-1}}=
x_{1}\frac{1\oplus y_1}{1+ y_1} \frac{1+ y_1^{-1}}{1\oplus y_1^{-1}}
=x_1,
\quad
y_{1;t_1}^{-1}=y_1.
\end{align}
The associated cluster algebra $\calA(\bfSigma)$  is the $\bbZ\bbP$-algebra generated
by
\begin{align}
x_1, \ x_{1}^{-1}\frac{1+ y_1}{1\oplus y_1}.
\end{align}
This is called a cluster algebra of \emph{type $A_1$}, which depends on the choice of
a semifield $\bbP$ and also the initial coefficient $y_1$.
\end{ex}

\begin{ex}[Type $A_1\times A_1$]
\label{1ex:A1A11}
Let $n=2$.
We consider the following arrangement of
a cluster pattern on $\bbT_2$.
$$
 \begin{xy}
(-25,0)*{\cdots},
(-10,0)*{\bullet},
(0,0)*{\bullet},
(10,0)*{\bullet},
(20,0)*{\bullet},
(30,0)*{\bullet},
(45,0)*{\cdots},
(-15,-3)*{1},
(-5,-3)*{2},
(5,-3)*{1},
(15,-3)*{2},
(25,-3)*{1},
(35,-3)*{2},
(-10,-3)*{t_{-1}},
(0,-3)*{t_0},
(10,-3)*{t_1},
(20,-3)*{t_2},
(30,-3)*{t_3},
(-9,3)*{\Sigma_{t_{-1}}},
(1,3)*{\Sigma_{t_{0}}},
(11,3)*{\Sigma_{t_{1}}},
(21,3)*{\Sigma_{t_{2}}},
(31,3)*{\Sigma_{t_{3}}},
\ar@{-} (-20,0);(-10,0)
\ar@{-} (-10,0);(0,0)
\ar@{-} (0,0);(10,0)
\ar@{-} (10,0);(20,0)
\ar@{-} (20,0);(30,0)
\ar@{-} (30,0);(40,0)
\end{xy}
$$
Below we use the simplified notations such as $
\Sigma_{t_s}=\Sigma_s$, $B_{t_s}=B_s$, $\bfx_{t_s}=\bfx_s$, $ \bfy_{t_s}=\bfy_s$, etc.
Let $t_0$ be the initial vertex,
and we set $\bfx_{0}=\bfx$, $\bfy_{0}=\bfy$.
Again, we  consider the simplest case
\begin{align}
B_0=B=O,
\end{align}
so that the initial $\hat{y}$ variables are given by
\begin{align}
\hat{y}_1=y_1,
\quad
\hat{y}_2=y_2.
\end{align}
For any $s\in \bbZ$, we have
\begin{align}
B_s = O.
\end{align}
Accordingly, we have
\begin{alignat}{3}
&
\begin{cases}
\displaystyle
 x_{1;1}=x_1^{-1}\frac{1+ y_1}{1\oplus y_1},
\\
  x_{2;1}=x_2,
  \end{cases}
  &\quad
  &
  \begin{cases}
 y_{1;1}=y_1^{-1},\\ 
 y_{2;1}=y_2 ,
 \end{cases}
 \\ 
&
\begin{cases}
\displaystyle
  x_{1;2}=x_1^{-1}\frac{1+ y_1}{1\oplus y_1},
  \\
    \displaystyle
  x_{1;2}=x_2^{-1}\frac{1+ y_2}{1\oplus y_2},
  \rule{0pt}{18pt}\hskip-1pt
\end{cases}  
    &\quad 
  &\begin{cases} 
 y_{1;2}=y_1^{-1} ,\\ 
 y_{2;2}=y_2^{-1},
 \end{cases}
 \\ 
&
\begin{cases}
   x_{1;3}=x_1,\\
   \displaystyle
  x_{2;3}=x_2^{-1}\frac{1+ y_2}{1\oplus y_2},
\end{cases}
&\quad
&
\begin{cases}
 y_{1;3}=y_1 ,\\ 
 y_{2;3}=y_2^{-1},
 \end{cases}
 \\ 
&
\begin{cases}
    x_{1;4}=x_1,\\
  x_{1;4}=x_2,
  \end{cases}
  &\quad
  &
  \begin{cases}
 y_{1;4}=y_1 ,\\ 
 y_{2;4}=y_2,
 \end{cases}
\end{alignat}
where we did the same calculation as  \eqref{1eq:a11}
in the last two steps.
We observe the periodicity with period 4. Namely,
\begin{align}
\bfx_{s+4}=\bfx_s,
\quad
\bfy_{s+4}=\bfy_s,
\quad
B_{s+4}=B_s
\quad
(s\in \bbZ).
\end{align}
In other words, $\Sigma_{s+4}=\Sigma_s$.
This is regarded as a special case of \eqref{1eq:A1A12}.
The cluster algebra $\calA(\bfSigma)$ is the $\bbZ\bbP$-algebra generated
by
\begin{align}
\label{1eq:A1A11}
x_1, \ x_2, \ x_{1}^{-1}\frac{1+ y_1}{1\oplus y_1},
\  x_2^{-1}\frac{1+ y_2}{1\oplus y_2}.
\end{align}
This is called a cluster algebra of \emph{type  $A_1\times A_1$}.
\end{ex}

Observe that a cluster algebra of type  $A_1\times A_1$ obtained above is
 isomorphic to the tensor product of two cluster algebras of type $A_1$ as a $\bbZ\bbP$-algebra.
One can  extend  this result  to a  more general  situation.

\begin{prop}
\label{1prop:Bdecom2}
Let $\bfSigma=\bfSigma(\Sigma_{t_0})$ be a cluster pattern of rank $n$ with
the initial seed $\Sigma_{t_0}=(\bfx,\bfy,B)$.
Suppose that the initial exchange matrix $B_{t_0}=B$ is decomposed as
\begin{align}
\label{1eq:Bdecom3}
B
=
\begin{pmatrix}
B'  & O\\
O & B''
\end{pmatrix},
\end{align}
where  the size of $B'$ and $B''$ are $n'$ and $n''$ $(n'+n''=n)$, respectively.
Accordingly, consider cluster patterns $\bfSigma'$ and $\bfSigma''$
of rank $n'$ and $n''$
whose initial seeds are given by
\begin{align}
\Sigma'_{t_0}=((x_i)_{i=1}^{n'}, (y_i)_{i=1}^{n'}, B'),
\quad
\Sigma''_{t_0}=((x_i)_{i=n'+1}^{n}, (y_i)_{i=n'+1}^{n}, B''),
\end{align}
respectively.
Then, as $\bbZ\bbP$-algebras, we have
\begin{align}
\calA(\bfSigma)\simeq \calA(\bfSigma')\otimes_{\bbZ\bbP}\calA(\bfSigma'').
\end{align}
\end{prop}
\begin{proof}
We first note that,
by Proposition \ref{1prop:Bdecom1}, for any $t\in \bbT_n$, the exchange matrix $B_{t}$
of $\bfSigma$ is decomposed 
in the same form as \begin{align}
\label{1eq:Bdecom4}
B_{t}
=
\begin{pmatrix}
B'_{t}  & O\\
O & B''_{t}
\end{pmatrix}.
\end{align}
For a seed $\Sigma_{t}=(\bfx_t, \bfy_t, B_t)$ of $\bfSigma$,   let us consider the mutation
at $k\leq n'$. Then, the mutation of $B_t$ is given by the first case
in \eqref{1eq:Bdecom2}. Also, for $ i \leq n'$,
$x_{i;t}$ and $y_{i;t}$ mutate effectively by the sub-exchange matrix $B'_t$,
while, for $i> n' $, they are stable.
The other case $k> n'$ is similar.
It follows that
the set of  the cluster variables of $\bfSigma$ is 
the disjoint union of  those of $\bfSigma'$ or $\bfSigma''$.
Also, there is no nontrivial algebraic relation
between cluster variables of $\bfSigma'$ and those of $\bfSigma''$.
Therefore, $\calA(\bfSigma)$ is naturally isomorphic 
to $\calA(\bfSigma')\otimes_{\bbZ\bbP} \calA(\bfSigma'')$
under the identification $z z'\mapsto z\otimes_{\bbZ\bbP} z'$
($z\in\calA(\bfSigma)$, $z'\in \calA(\bfSigma')$).
\end{proof}

We say that
a square matrix $M$ is \emph{decomposable\/} (resp. \emph{indecomposable\/})
\index{matrix!decomposable} \index{matrix!indecomposable}
if, after some simultaneous permutation of the column and row indices of $M$,
$M$ is a direct sum of two square matrices
(resp. otherwise).

Thanks to Proposition \ref{1prop:Bdecom2},
in many situations it is  enough to concentrate
on cluster patterns and cluster algebras with
\emph{ indecomposable\/}  exchange matrices.

\subsection{Rank 2 periodicities}
\label{1sec:rank2}

We continue to use the parametrization and the notation for
cluster patterns of rank 2 in Example \ref{1ex:A1A11}.
 Let us take the following initial exchange matrix
\begin{align}
\label{1eq:exchange1}
B_0=B=
\begin{pmatrix}
0 & -b\\
a & 0
\end{pmatrix},
\quad
(a,b >0).
\end{align}
By the mutations \eqref{1eq:bmut2},
the exchange matrices are given by
\begin{align}
B_{s}=
\begin{cases}
B & \mbox{$s$: even},\\
-B & \mbox{$s$: odd}.\\
\end{cases}
\end{align}
Accordingly, we have, 
for even $s$,
\begin{gather}
\hat{y}_{1;s}=y_{1;s}x_{2;s}^a,
\quad
\hat{y}_{2;s}=y_{2;s}x_{1;s}^{-b},
\\
\begin{cases}
\displaystyle
x_{1;{s+1}}=x_{1;s}^{-1}\frac{ 1+\hat{y}_{1;s}}{1\oplus y_{1;s}},
\\
x_{2;{s+1}}=x_{2;s},
\end{cases}
\quad
\begin{cases}
y_{1;{s+1}}=y_{1;s}^{-1},
\\
y_{2;{s+1}}=y_{2;s} (1\oplus y_{1;s})^b,
\end{cases}
\end{gather}
and,
for odd $s$,
\begin{gather}
\hat{y}_{1;s}=y_{1;s}x_{2;s}^{-a},
\quad
\hat{y}_{2;s}=y_{2;s}x_{1;s}^{b},
\\
\label{1eq:sodd1}
\begin{cases}
x_{1;{s+1}}=x_{1;s},
\\
\displaystyle
x_{2;{s+1}}=x_{2;s}^{-1} \frac{ 1+\hat{y}_{2;s}}{1\oplus y_{2;s}},
\end{cases}
\quad
\begin{cases}
y_{1;{s+1}}=y_{1;s} (1\oplus y_{2;s})^a,
\\
y_{2;{s+1}}=y_{2;s}^{-1}.
\end{cases}
\end{gather}

Below we concentrate on the case $ab\leq 3$.
Without losing generality we may assume that $b=1$
and $a=1,2,3$,
because the opposite case is obtained by exchanging
the indices 1 and 2.
We strongly recommend the readers to carry out
the following calculations
to appreciate the systematic occurrences of  reductions (``miracles'').
It may require a few hours. So, be prepared!

\begin{ex}[Type $A_2$]
\label{1ex:A21}
Consider the case $a=1$, where
\begin{align}
B_0=B=
\begin{pmatrix}
0 & -1\\
1 & 0
\end{pmatrix},
\end{align}
so that the initial $\hat{y}$-variables are given by
\begin{align}
\label{1eq:yhatini1}
\hat{y}_{1}=y_{1}x_{2},
\quad
\hat{y}_{2}=y_{2}x_{1}^{-1}.
\end{align}
In fact, this is the case already calculated in Section \ref{1sec:first1},
where  cluster variables are without coefficients.
Here, we effectively use Proposition \ref{1prop:yhat1} so that
we do not have to look into the contents of $\hat{y}$-variables during
the calculation
 except for the initial one \eqref{1eq:yhatini1}.
We have the following result:
\begin{alignat}{3}
\label{1eq:A1mut1}
&
\begin{cases}
\displaystyle
 x_{1;1}=x_1^{-1}\frac{1+ \hat{y}_1}{1\oplus y_1},
\\
  x_{2;1}=x_2,
  \end{cases}
  &\quad
  &
  \begin{cases}
 y_{1;1}=y_1^{-1},\\ 
   y_{2;1}=y_2 (1\oplus y_1),
 \end{cases}
 \\
 &
\begin{cases}
 \displaystyle
 x_{1;2}=x_1^{-1}\frac{1+ \hat{y}_1}{1\oplus y_1},
\\
 \displaystyle
  x_{2;2}=x_2^{-1}\frac{1+ \hat{y}_2+\hat{y}_1\hat{y}_2}{1\oplus y_2\oplus y_1y_2},
  \rule{0pt}{18pt}\hskip-1pt
  \end{cases}
  &\quad
  &
  \begin{cases}
 y_{1;2}=y_1^{-1}(1\oplus y_2\oplus y_1y_2),\\ 
   y_{2;2}=y_2^{-1} (1\oplus y_1)^{-1},
 \end{cases}
 \\
  &
\begin{cases}
 \displaystyle
 x_{1;3}=x_1 x_2^{-1}\frac{1+ \hat{y}_2}{1\oplus y_2},
\\
 \displaystyle
  x_{2;3}=x_2^{-1}\frac{1+ \hat{y}_2+\hat{y}_1\hat{y}_2}{1\oplus y_2\oplus y_1y_2},
  \rule{0pt}{18pt}\hskip-1pt
  \end{cases}
  &\quad
  &
  \begin{cases}
 y_{1;3}=y_1(1\oplus y_2\oplus y_1y_2)^{-1},\\ 
   y_{2;3}=y_1^{-1}y_2^{-1} (1\oplus y_2),
 \end{cases}
 \\
  &
\begin{cases}
 \displaystyle
 x_{1;4}=x_1 x_2^{-1} \frac{1+ \hat{y}_2}{1\oplus y_2},
\\
 \displaystyle
  x_{2;4}=x_1,
  \end{cases}
  &\quad
  &
  \begin{cases}
 y_{1;4}=y_2^{-1},\\ 
   y_{2;4}=y_1y_2 (1\oplus y_2)^{-1},
 \end{cases}
  \\
  \label{1eq:A1mut2}
  &
\begin{cases}
 \displaystyle
 x_{1;5}=x_2,
\\
 \displaystyle
  x_{2;5}=x_1,
  \end{cases}
  &\quad
  &
  \begin{cases}
 y_{1;5}=y_2,\\ 
   y_{2;5}=y_1.
 \end{cases}
 \end{alignat}
 Let us extend Observation \ref{1obs:x1} in more details.
 
(a). \emph{Periodicity/Finiteness.} 
\index{pentagon periodicity}
Even with the presence of coefficients
for cluster variables,
we still have the same pentagon periodicity
 \begin{align}
\Sigma_{s+5} =\tau_{12}(\Sigma_s),
 \end{align}
 where $\tau_{12}$ is the transpose of $1$ and $2$,
 and its action was defined in \eqref{1eq:saction1}.
 In particular, the cluster variables of $\bfSigma$ are exhausted by
 \begin{align}
 \label{1eq:A1cluster1}
 x_1, \quad x_2, \quad x_1^{-1}\frac{1+ \hat{y}_1}{1\oplus y_1},\quad 
 x_2^{-1}\frac{1+ \hat{y}_2+\hat{y}_1\hat{y}_2}{1\oplus y_2\oplus y_1y_2},
 \quad
  x_1 x_2^{-1} \frac{1+ \hat{y}_2}{1\oplus y_2}.
 \end{align}
 By \eqref{1eq:yhatini1},
 as Laurent polynomials in $x_1$ and $x_2$ with coefficients in $\bbZ\bbP$,
 they are written
 as 
\begin{align}
\label{1eq:A1cluster2}
 x_1, \quad x_2, \quad x_1^{-1}\frac{1+ {y}_1x_2}{1\oplus y_1},\quad 
 x_1^{-1}x_2^{-1}\frac{x_1+ {y}_2+{y}_1{y}_2x_2}{1\oplus y_2\oplus y_1y_2},\quad
  x_2^{-1} \frac{ x_1+ {y}_2}{1\oplus y_2}.
 \end{align}
 \par
 (b). \emph{Laurent phenomenon.} 
 In \eqref{1eq:A1cluster2},
 each cluster variable $x_{i;t}$ is expressed as a Laurent polynomial in
 the initial cluster variables $\bfx$ with  coefficients in $\bbZ\bbP$.
 \par
 (c). \emph{Laurent positivity.}
  Every coefficient of the above Laurent polynomial is
 nonnegative in $\bbZ\bbP$.
\par
(d). \emph{$F$-polynomials.}
Let us rephrase
the properties (b) and (c)  in a more specific way.
 Each cluster variable $x_{i;t}$  is expressed in 
 a unified way as follows:
 \begin{align}
 x_{i;t}
 =
 \Biggl(
 \prod_{j=1}^n x_j^{g_{ji;t}}
 \Biggr)
 \frac
 {F_{i;t}(\hat\bfy)}{F_{i;t}\vert_{\bbP}(\bfy)},
 \end{align}
 where $F_{i;t}(\bfy)$ is a polynomial in formal variables $\bfy=(y_1,y_2)$
 with nonnegative integer coefficients, 
 and $F_{i;t}\vert_{\bbP}(\bfy)$ is the specialization in $\bbP$ at the initial coefficients $\bfy_{t_0}=\bfy$ in Definition \ref{1defn:sp1}.
 (Here, we conveniently abuse the symbol $\bfy$ for two different usages.)
 Later, $F_{i;t}(\bfy)$ is called an \emph{$F$-polynomial}.
 \par
 (e). \emph{Unit constant term.} 
Every $F$-polynomial $F_{i;t}(\bfy)$  has constant term 1.
 \par
 (f). \emph{Duality.}
  As a new observation, we see that the cluster variables and
 coefficients share some common/parallel structure.
In particular, the  above $F$-polynomials  also appear  for  coefficients.
 \end{ex}

\begin{ex}[Type $B_2$]
Consider the case $a=2$, where
\begin{align}
B_0=B=
\begin{pmatrix}
0 & -1\\
2 & 0
\end{pmatrix},
\quad
D=
\begin{pmatrix}
2 & 0\\
0 & 1
\end{pmatrix},
\end{align}
so that the initial $\hat{y}$-variables are given by
\begin{align}
\label{1eq:yhatini2}
\hat{y}_{1}=y_{1}x_{2}^2,
\quad
\hat{y}_{2}=y_{2}x_{1}^{-1}.
\end{align}
We have the following result:
\begin{alignat}{3}
&
\begin{cases}
\displaystyle
 x_{1;1}=x_1^{-1}\frac{1+ \hat{y}_1}{1\oplus y_1},
\\
  x_{2;1}=x_2,
  \end{cases}
  \quad
  \begin{cases}
 y_{1;1}=y_1^{-1},\\ 
   y_{2;1}=y_2 (1\oplus y_1),
 \end{cases}
 \\
 &
\begin{cases}
 \displaystyle
 x_{1;2}=x_1^{-1}\frac{1+ \hat{y}_1}{1\oplus y_1},
\\
 \displaystyle
  x_{2;2}=x_2^{-1}\frac{1+ \hat{y}_2+\hat{y}_1\hat{y}_2}{1\oplus y_2\oplus y_1y_2},
  \rule{0pt}{18pt}\hskip-1pt
  \end{cases}
  \quad
  \begin{cases}
 y_{1;2}=y_1^{-1}(1\oplus y_2\oplus y_1y_2)^2,\\ 
   y_{2;2}=y_2^{-1} (1\oplus y_1)^{-1},
 \end{cases}
 \\
  \begin{split}
  &
\begin{cases}
 \displaystyle
 x_{1;3}=x_1 x_2^{-2}\frac{1+ 2\hat{y}_2+\hat{y}_2^2 + \hat{y}_1\hat{y}_2^2}
 {1\oplus 2 y_2\oplus y_2^2 \oplus y_1y_2^2},
\\
 \displaystyle
  x_{2;3}=x_2^{-1}\frac{1+ \hat{y}_2+\hat{y}_1\hat{y}_2}{1\oplus y_2\oplus y_1y_2},
  \rule{0pt}{18pt}\hskip-1pt
  \end{cases}
\\
&
  \begin{cases}
 y_{1;3}=y_1(1\oplus y_2\oplus y_1y_2)^{-2},\\ 
   y_{2;3}=y_1^{-1}y_2^{-1} (1\oplus 2 y_2\oplus y_2^2 \oplus y_1 y_2^2 ),
 \end{cases}
 \end{split}
 \\
     \begin{split}
  &
  \label{1eq:freeex1}
\begin{cases}
 \displaystyle
 x_{1;4}=x_1 x_2^{-2}\frac{1+ 2\hat{y}_2+\hat{y}_2^2 + \hat{y}_1\hat{y}_2^2}
 {1\oplus 2 y_2\oplus y_2^2 \oplus y_1y_2^2},
\\
 \displaystyle
  x_{2;4}=x_1x_2^{-1}\frac{1+ \hat{y}_2}{1\oplus y_2},
  \rule{0pt}{18pt}\hskip-1pt
  \end{cases}
\\
&
  \begin{cases}
 y_{1;4}=y_1^{-1}y_2^{-2}(1\oplus y_2)^2,\\ 
   y_{2;4}=y_1y_2  (1\oplus 2 y_2\oplus y_2^2 \oplus y_1 y_2^2 )^{-1},
 \end{cases}
 \end{split}
  \\
  &
\begin{cases}
 \displaystyle
 x_{1;5}=x_1,
\\
 \displaystyle
  x_{2;5}=x_1x_2^{-1}\frac{1+ \hat{y}_2}{1\oplus y_2},
  \end{cases}
\quad
  \begin{cases}
 y_{1;5}=y_1y_2^2(1\oplus y_2)^{-2},\\ 
   y_{2;5}=y_2^{-1},
 \end{cases}
 \\
   &
\begin{cases}
 \displaystyle
 x_{1;6}=x_1,
\\
 \displaystyle
  x_{2;6}=x_2,
  \end{cases}
\quad
  \begin{cases}
 y_{1;6}=y_1,\\ 
   y_{2;6}=y_2.
 \end{cases}
 \end{alignat}
 During the calculation we  used the following polynomial identities:
 \begin{align}
 y_1+(1+y_2+y_1y_2)^2=(1+y_1)(1+2y_2+y_2^2+y_1y_2^2),\\
 y_1y_2+(1+2y_2+y_2^2+y_1y_2^2)=(1+y_2+y_1y_2)(1+y_2).
 \end{align}
We observe the periodicity
\begin{align}
\Sigma_{s+6}=\Sigma_s.
\end{align}
 In particular, the cluster variables of $\bfSigma$ are exhausted by
 \begin{align}
 \begin{split}
& x_1, \quad x_2, \quad   x_1^{-1}\frac{1+ \hat{y}_1}{1\oplus y_1},\quad
x_2^{-1}\frac{1+ \hat{y}_2+\hat{y}_1\hat{y}_2}{1\oplus y_2\oplus y_1y_2},
\\
& 
 x_1 x_2^{-2}\frac{1+ 2\hat{y}_2+\hat{y}_2^2 + \hat{y}_1\hat{y}_2^2}
 {1\oplus 2 y_2\oplus y_2^2 \oplus y_1y_2^2},
 \quad
  x_1 x_2^{-1} \frac{1+ \hat{y}_2}{1\oplus y_2}.
 \end{split}
 \end{align}
 As Laurent polynomials in $x_1$ and $x_2$, they are written as
 \begin{align}
 \begin{split}
&  x_1, \quad x_2, \quad   x_1^{-1}\frac{1+ y_1x_2^2}{1\oplus y_1},\quad
  x_1^{-1} x_2^{-1}\frac{x_1 + y_2+y_1y_2x_2^2}{1\oplus y_2\oplus y_1y_2},
  \\
&  
 x_1^{-1} x_2^{-2}\frac{x_1^2+ 2y_2x_1+y_2^2 + y_1y_2^2 x_2^2 }
 {1\oplus 2 y_2\oplus y_2^2 \oplus y_1y_2^2},
 \quad
  x_2^{-1} \frac{x_1 +y_2}{1\oplus y_2}.
 \end{split}
 \end{align}
 All other observations for type $A_2$ are commonly applied.
\end{ex}

\begin{ex}[Type $G_2$]
Consider the case $a=3$, where
\begin{align}
B_0=B=
\begin{pmatrix}
0 & -1\\
3 & 0
\end{pmatrix},
\quad
D=
\begin{pmatrix}
3 & 0\\
0 & 1
\end{pmatrix},
\end{align}
so that the initial $\hat{y}$-variables are given by
\begin{align}
\label{1eq:yhatini3}
\hat{y}_{1}=y_{1}x_{2}^3,
\quad
\hat{y}_{2}=y_{2}x_{1}^{-1}.
\end{align}
We have the following result:
\begin{alignat}{3}
&
\begin{cases}
\displaystyle
 x_{1;1}=x_1^{-1}\frac{1+ \hat{y}_1}{1\oplus y_1},
\\
  x_{2;1}=x_2,
  \end{cases}
  \quad
  \begin{cases}
 y_{1;1}=y_1^{-1},\\ 
   y_{2;1}=y_2 (1\oplus y_1),
 \end{cases}
 \\
 &
\begin{cases}
 \displaystyle
 x_{1;2}=x_1^{-1}\frac{1+ \hat{y}_1}{1\oplus y_1},
\\
 \displaystyle
  x_{2;2}=x_2^{-1}\frac{1+ \hat{y}_2+\hat{y}_1\hat{y}_2}{1\oplus y_2\oplus y_1y_2},
    \rule{0pt}{18pt}\hskip-1pt
  \end{cases}
  \quad
  \begin{cases}
 y_{1;2}=y_1^{-1}(1\oplus y_2\oplus y_1y_2)^3,\\ 
   y_{2;2}=y_2^{-1} (1\oplus y_1)^{-1},
 \end{cases}
 \\
     \begin{split}
  &
\begin{cases}
 \displaystyle
 x_{1;3}=x_1 x_2^{-3}\frac{1+ 3\hat{y}_2+3\hat{y}_2^2 + \hat{y}_2^3
 + 3 \hat{y}_1\hat{y}_2^2 + 2 \hat{y}_1\hat{y}_2^3 + \hat{y}_1^2\hat{y}_2^3
 }
 {1\oplus 3 y_2\oplus 3 y_2^2 \oplus y_2^3
 \oplus 3y_1y_2^2\oplus 2y_1 y_2^3\oplus y_1^2y_2^3},
\\
 \displaystyle
  x_{2;3}=x_2^{-1}\frac{1+ \hat{y}_2+\hat{y}_1\hat{y}_2}{1\oplus y_2\oplus y_1y_2},
    \rule{0pt}{18pt}\hskip-1pt
  \end{cases}
\\
&
  \begin{cases}
 y_{1;3}=y_1(1\oplus y_2\oplus y_1y_2)^{-3},\\ 
   y_{2;3}=y_1^{-1}y_2^{-1} (1\oplus 3 y_2\oplus 3 y_2^2 \oplus y_2^3
 \oplus 3y_1y_2^2\oplus 2y_1 y_2^3\oplus y_1^2y_2^3 ),
 \end{cases}
 \end{split}
 \\
     \begin{split}
  &
\begin{cases}
 \displaystyle
 x_{1;4}=x_1 x_2^{-3}\frac{1+ 3\hat{y}_2+3\hat{y}_2^2 + \hat{y}_2^3
 + 3 \hat{y}_1\hat{y}_2^2 + 2 \hat{y}_1\hat{y}_2^3 + \hat{y}_1^2\hat{y}_2^3
 }
 {1\oplus 3 y_2\oplus 3 y_2^2 \oplus y_2^3
 \oplus 3y_1y_2^2\oplus 2y_1 y_2^3\oplus y_1^2y_2^3},
\\
 \displaystyle
  x_{2;4}=x_1x_2^{-2}\frac{1+ 2 \hat{y}_2+ \hat{y}_2^2
  + \hat{y}_1 \hat{y}_2}{1\oplus 2y_2\oplus y_2^2 \oplus y_1y_2^2},
    \rule{0pt}{18pt}\hskip-1pt
  \end{cases}
\\
&
  \begin{cases}
 y_{1;4}=y_1^{-2}y_2^{-3}(1\oplus 2y_2\oplus y_2^2 \oplus y_1y_2^2)^3,\\ 
   y_{2;4}=y_1y_2  ( {1\oplus 3 y_2\oplus 3 y_2^2 \oplus y_2^3
 \oplus 3y_1y_2^2\oplus 2y_1 y_2^3\oplus y_1^2y_2^3} )^{-1},
 \end{cases}
 \end{split}
  \\
      \begin{split}
    &
\begin{cases}
 \displaystyle
 x_{1;5}=x_1^{2} x_2^{-3}\frac{1+ 3\hat{y}_2+3\hat{y}_2^2 + \hat{y}_2^3
 +  \hat{y}_1\hat{y}_2^3
 }
 {1\oplus 3 y_2\oplus 3 y_2^2 \oplus y_2^3
 \oplus y_1y_2^3},
\\
 \displaystyle
  x_{2;5}=x_1x_2^{-2}\frac{1+ 2 \hat{y}_2+ \hat{y}_2^2
  + \hat{y}_1 \hat{y}_2}{1\oplus 2y_2\oplus y_2^2 \oplus y_1y_2^2},
    \rule{0pt}{18pt}\hskip-1pt
  \end{cases}
\\
&
  \begin{cases}
 y_{1;5}=y_1^{2}y_2^{3}(1\oplus 2y_2\oplus y_2^2 \oplus y_1y_2^2)^{-3},\\ 
   y_{2;5}=y_1^{-1}y_2^{-2}  ( 1\oplus 3 y_2\oplus 3 y_2^2 \oplus y_2^3
 \oplus y_1y_2^3),
 \end{cases}
 \end{split}
  \\
      \begin{split}
      &
\begin{cases}
 \displaystyle
 x_{1;6}=x_1^{2} x_2^{-3}\frac{1+ 3\hat{y}_2+3\hat{y}_2^2 + \hat{y}_2^3
 +  \hat{y}_1\hat{y}_2^3
 }
 {1\oplus 3 y_2\oplus 3 y_2^2 \oplus y_2^3
 \oplus y_1y_2^3},
\\
 \displaystyle
  x_{2;6}=x_1x_2^{-1}\frac{1+  \hat{y}_2}{1\oplus y_2},
    \rule{0pt}{18pt}\hskip-1pt
  \end{cases}
\\
&
  \begin{cases}
 y_{1;6}=y_1^{-1}y_2^{-3}(1\oplus y_2)^{3},\\ 
   y_{2;6}=y_1y_2^{2}  ( 1\oplus 3 y_2\oplus 3 y_2^2 \oplus y_2^3
 \oplus y_1y_2^3)^{-1},
 \end{cases}
 \end{split}
  \\
  &
\begin{cases}
 \displaystyle
 x_{1;7}=x_1,
\\
 \displaystyle
  x_{2;7}=x_1x_2^{-1}\frac{1+ \hat{y}_2}{1\oplus y_2},
  \end{cases}
\quad
  \begin{cases}
 y_{1;7}=y_1y_2^3(1\oplus y_2)^{-3},\\ 
   y_{2;7}=y_2^{-1},
 \end{cases}
 \\
   &
\begin{cases}
 \displaystyle
 x_{1;8}=x_1,
\\
 \displaystyle
  x_{2;8}=x_2,
  \end{cases}
\quad
  \begin{cases}
 y_{1;8}=y_1,\\ 
   y_{2;8}=y_2.
 \end{cases}
 \end{alignat}
 During the calculation we  used the following polynomial identities:
 \begin{align}
 \begin{split}
 &y_1+(1+y_2+y_1y_2)^3\\
 &= (1+y_1)(1+3y_2+3y_2^2+
 y_2^3+ 3y_1y_2^2 + 2y_1y_2^3 + y_1^2y_2^3),
 \end{split}
 \\
 \begin{split}
& y_1y_2+(1+3y_2+3y_2^2+
 y_2^3+ 3y_1y_2^2 + 2y_1y_2^3 + y_1^2y_2^3)\\
& =(1+y_2+y_1y_2)
 (1+2y_2+y_2^2 + y_1y_2^2),
\end{split} 
\\
 \begin{split}
 \label{1eq:fac1}
& y_1^2y_2^3+ (1+2y_2+y_2^2 + y_1y_2^2)^3\\
& =(1+3y_2+3y_2^2+
 y_2^3+ 3y_1y_2^2 + 2y_1y_2^3 + y_1^2y_2^3)\\
 &\qquad\times
 (1+3y_2+3y_2^2 +y_2^3 + y_1y_2^3),
\end{split} 
\\
 \begin{split}
& y_1y_2^2+  (1+3y_2+3y_2^2 +y_2^3 + y_1y_2^3)\\
& =(1+2y_2+y_2^2 + y_1y_2^2)(1+y_2),
\end{split} 
 \end{align}
 where the most complicated identity \eqref{1eq:fac1} may be quickly checked by computer.
We observe the periodicity
\begin{align}
\Sigma_{s+8}=\Sigma_s.
\end{align}
 In particular, the cluster variables of $\bfSigma$ are exhausted by
 \begin{align}
 \begin{split}
& x_1, \quad x_2, \quad   x_1^{-1}\frac{1+ \hat{y}_1}{1\oplus y_1},\quad 
 x_2^{-1}\frac{1+ \hat{y}_2+\hat{y}_1\hat{y}_2}{1\oplus y_2\oplus y_1y_2},\quad
\\
& x_1 x_2^{-3}\frac{1+ 3\hat{y}_2+3\hat{y}_2^2 + \hat{y}_2^3
 + 3 \hat{y}_1\hat{y}_2^2 + 2 \hat{y}_1\hat{y}_2^3 + \hat{y}_1^2\hat{y}_2^3
 }
 {1\oplus 3 y_2\oplus 3 y_2^2 \oplus y_2^3
 \oplus 3y_1y_2^2\oplus 2y_1 y_2^3\oplus y_1^2y_2^3},
 \\
&
 x_1 x_2^{-2}\frac{1+ 2\hat{y}_2+\hat{y}_2^2 + \hat{y}_1\hat{y}_2^2}
 {1\oplus 2 y_2\oplus y_2^2 \oplus y_1y_2^2},
 \\
 &
x_1^{2} x_2^{-3}\frac{1+ 3\hat{y}_2+3\hat{y}_2^2 + \hat{y}_2^3
 +  \hat{y}_1\hat{y}_2^3
 }
 {1\oplus 3 y_2\oplus 3 y_2^2 \oplus y_2^3
 \oplus y_1y_2^3},
 \quad
  x_1 x_2^{-1} \frac{1+ \hat{y}_2}{1\oplus y_2}.
 \end{split}
 \end{align}
 As Laurent polynomials in $x_1$ and $x_2$, they are written as
 \begin{align}
 \begin{split}
& x_1, \quad x_2, \quad   x_1^{-1}\frac{1+ y_1x_2^3}{1\oplus y_1},\quad 
x_1^{-1} x_2^{-1}\frac{x_1+ y_2+y_1y_2x_2^3}{1\oplus y_2\oplus y_1y_2},\
\\
& x_1^{-2} x_2^{-3}\frac{x_1^3+ 3y_2x_1^{2}+3 y_2^2x_1+ y_2^3
 + 3 y_1y_2^2x_1x_2^3 + 2 y_1 y_2^3x_2^3 + y_1^2y_2^3x_2^6
 }
 {1\oplus 3 y_2\oplus 3 y_2^2 \oplus y_2^3
 \oplus 3y_1y_2^2\oplus 2y_1 y_2^3\oplus y_1^2y_2^3},
 \\
&
 x_1^{-1} x_2^{-2}\frac{x_1^2+ 2y_2x_1+ y_2^2 + y_1y_2^2x_2^3 }
 {1\oplus 2 y_2\oplus y_2^2 \oplus y_1y_2^2},
 \\
 &
x_1^{-1} x_2^{-3}\frac{x_1^3+ 3y_2x_1^{2}+3 y_2^2x_1+ y_2^3
 +  y_1y_2^3x_2^3
 }
 {1\oplus 3 y_2\oplus 3 y_2^2 \oplus y_2^3
 \oplus y_1y_2^3},
 \quad
 x_2^{-1} \frac{x_1+ y_2}{1\oplus y_2}.
 \end{split}
 \end{align}
 All other observations for type $A_2$ are commonly applied.
\end{ex}

Congratulations! You have successfully gone through the famous ``ordeal''
in cluster algebra theory.

\begin{rem}
(a). For the matrices in \eqref{1eq:exchange1} with $ab\geq 4$,
the periodicity of seeds does not occur. 
Accordingly,  we have infinitely many distinct cluster variables
for those cluster patterns.
\par
(b).
The  periods 5, 6, 8 in the above examples are interpreted as
$h+2$, where $h=3,4,6$ are the \emph{Coxeter numbers\/} \index{Coxeter number} of
the root systems of type $A_2$, $B_2$, $G_2$, respectively.
A more account on
the connection to the root systems will be given in Section \ref{1subsec:finite1}.
\end{rem}

\subsection{Free coefficients}
Let us introduce a  notion,
which is not defined in CA1-4
explicitly.

\begin{defn}
\index{coefficient!free}
\label{1defn:free1}
We say that a cluster pattern $\bfSigma$
is \emph{with free coefficients at $t_0\in \bbT_n$}
if the following conditions are satisfied:
\begin{itemize}
\item
The coefficient semifield of $\bfSigma$ is the universal semifield
$\bbQ_{\mathrm{sf}}(\bfy)$ with generators $\bfy=(y_1,\dots,y_n)$,
where $n$ is the rank of $\bfSigma$.
\item
The coefficient tuple $\bfy_{t_0}$ at $t_0$ coincides with $\bfy$.
\end{itemize}
Note that, for each $t$,
$y_{1;t}$, \dots, $y_{n;t}$ are algebraically independent.
\end{defn}

\begin{rem}
It might be natural to call such coefficients   \emph{universal coefficients}.
However, in CA4, a related but different notion is defined and called so.
Therefore, we call it differently to avoid confusion.
\end{rem}

\begin{rem}
\label{1rem:free1}
Though it is not necessary, we often identify the above base vertex $t_0$
with the initial vertex for $\bfSigma$.
In that case the notation $\bfy$ above is compatible with the convention 
for the initial coefficients \eqref{1eq:init1}.
Otherwise, we may dismiss the notation \eqref{1eq:init1} to avoid 
the conflict.
\end{rem}
For any $t\in \bbT_n$,  the coefficients $y_{1;t_0}$, \dots, $y_{n;t_0}$
in Definition \ref{1defn:free1}  are expressed as
rational functions in $\bfy_t$ with subtraction-free expressions
by doing mutations from $t$ to $t_0$.
They induce a canonical isomorphism
between the semifields
 $\bbQ_{\mathrm{sf}}(\bfy_{t_0})$ and $\bbQ_{\mathrm{sf}}(\bfy_{t})$.
Therefore,
under this identification
the base vertex can be shifted arbitrarily,
so that the  choice of the base point $t_0$ is superficial.

Recall the universality of $\bbQ_{\mathrm{sf}}(\bfy)$ in Proposition \ref{1prop:uni1}.
We are going to extend this homomorphism $\pi:\bbQ_{\mathrm{sf}}(\bfy)\rightarrow \bbP$
to a map between cluster variables  with coefficients in   $\bbQ_{\mathrm{sf}}(\bfy)$ and $\bbP$
sharing a common $B$-pattern.
However, 
there is a pitfall to avoid.
Let $\varphi:\bbP \rightarrow \bbP'$ be a semifield homomorphism.
Then, it is uniquely extended to  a ring homomorphism 
$\varphi_1:\bbZ\bbP \rightarrow \bbZ\bbP'$.
However, it can be extended to a field homomorphism
$\varphi_2: \bbQ\bbP \rightarrow \bbQ\bbP'$ only if $\varphi$
is injective. Indeed, if $\varphi$ is not injective, 
$\varphi_1$ is not injective.
Then, a fraction of $\bbZ\bbP$ whose denominator is in $\mathrm{Ker}\, \varphi_1$
does not have a well-defined image of $\varphi_2$.

Below, for $\bbP=\bbQ_{\mathrm{sf}}(\bfy)$,
we write $\bbQ\bbP$ as $\bbQ\bbQ_{\mathrm{sf}}(\bfy)$ according to our convention,
though it looks a little cumbersome.
\begin{prop}
\label{1prop:free1}
Let $\bfSigma$ be a cluster pattern with free coefficients at $t_0$.
Let $\bfSigma'$ be any cluster pattern with coefficients in any semifield $\bbP$
sharing the common $B$-pattern with $\bfSigma$.
\par
(a).
Let $\pi$ be the semifield homomorphism defined by
\begin{align}
\begin{matrix}
\pi: &\bbQ_{\mathrm{sf}}(\bfy)&\rightarrow &\bbP\\
& y_i=y_{i;t_0} & \mapsto &y'_{i;t_0}.
\end{matrix}
\end{align}
Then, we have
\begin{align}
\pi(y_{i;t})=y'_{i;t}.
\end{align}
\par
(b). 
 Let $\cal{X}(\bfSigma)$ be the set of the cluster variables of $\bfSigma$.
Let
 \begin{align}
\begin{matrix}
\varphi: &\cal{X}(\bfSigma)&\rightarrow 
&(\bbQ\bbP)(\bfx'_{t_0})
\end{matrix}
\end{align}
be  the map such that,
in each element $x_{i;t}\in (\bbQ\bbQ_{\mathrm{sf}}(\bfy))(\bfx_{t_0})$,
$x_{i;t_0}$ is replaced with $x'_{i;t_0}$, 
and $f(\bfy)\in \bbQ_{\mathrm{sf}}(\bfy)$ is replaced with $f(\bfy')=\pi(f(\bfy))$.
Then, the map is well defined; moreover, we have
\begin{align}
\varphi(x_{i;t})=x'_{i;t}.
\end{align}
\end{prop}
\begin{proof}
(a).
Any coefficient
$y_{i;t}$ is expressed as a rational function $Y_{i;t}(\bfy)\in
\bbQ_{\mathrm{sf}}(\bfy)$.
Applying $\pi$ yields the expression for $y'_{i;t}$
in $\bfy'_{t_0}$,
because both $y_{i;t}$ and $y'_{i;t}$ obey formally the same mutation formula
\eqref{1eq:ymut1}.
\par
(b).
Any cluster variable
$x_{i;t}$ is expressed as a rational function $X_{i;t}(\bfx_{t_0})\in
(\bbQ\bbQ_{\mathrm{sf}}(\bfy))(\bfx_{t_0})$.
Since the mutation \eqref{1eq:xmut1} does not involve any subtraction,
$X_{i;t}(\bfx_{t_0})$ has a subtraction-free expression in 
$(\bbQ\bbQ_{\mathrm{sf}}(\bfy))(\bfx_{t_0})$.
In particular,  for the denominator of any coefficient in $X_{i;t}(\bfx_{t_0})$,
its image by $\pi$ does not vanish.
Thus, the map $\varphi$ is well defined.
Applying $\varphi$ yields the expression for  $x'_{i;t}$
in $\bfx'_{t_0}$,
because both $x_{i;t}$ and $x'_{i;t}$ obey formally the same mutation formula
\eqref{1eq:xmut1}.
\end{proof}

The result is plainly rephrased that  cluster variables with \emph{any specific choice 
of coefficients\/}
 can be obtained from cluster variables
with \emph{free coefficients\/} by \emph{specializing the coefficients}.
To see how it works more concretely,
let us look at the formulas in Section \ref{1sec:rank2} (e.g.,
\eqref{1eq:A1mut1}--\eqref{1eq:A1mut2}).
They can
be viewed as expressions
of $x_{i;t}$
 in $(\bbQ\bbQ_{\mathrm{sf}}(\bfy))(\bfx_{t_0})$
and $y_{i;t}$ in $\bbQ_{\mathrm{sf}}(\bfy)$,
because we did not use any specific property of a given 
semifield $\bbP$ and initial coefficients $\bfy$ therein.
Meanwhile, they can be viewed  also as expressions 
for any specific choice of coefficients.
(Viewing so is nothing but the specialization by $\varphi$ and $\pi$.)

The following is an immediate consequence
of Proposition \ref{1prop:free1}.

\begin{prop}
\label{1prop:period1}
Let $\bfSigma$ be a cluster pattern with free coefficients at $t_0$.
Let $\bfSigma'$ be a cluster pattern with  coefficients in any semifield $\bbP$
sharing the common $B$-pattern with $\bfSigma$.
Then, for any $t,t'\in \bbT_n$ and a permutation $\sigma\in S_n$,
the following fact holds:
\begin{align}
\label{1eq:period1}
\Sigma_t=
\sigma \Sigma_{t'}
\quad
\Longrightarrow
\quad
\Sigma'_t=
\sigma \Sigma'_{t'}.
\end{align}
In other words, any periodicity of $\bfSigma$ implies
the same periodicity of $\bfSigma'$.
\end{prop}
\begin{proof}
The equality $\Sigma_t=
\sigma \Sigma_{t'}$ reduces
to the equality $\Sigma'_t=
\sigma \Sigma'_{t'}$
under the specialization in 
Proposition \ref{1prop:free1}.
\end{proof}

Whether the opposite implication of \eqref{1eq:period1} holds
or not is an important issue posed in CA4.
Namely,   
there is a possibility that  additional periodicities
of $\bfSigma'$
occur
 under some specialization of coefficients,
 though it does not happen in the rank 2 examples
 in  Section \ref{1sec:rank2}.
We will discuss more about the problem later
in Section \ref{1subsec:further1}.

\newpage

\section{Fundamental Results}

In this section we present some of the most fundamental
results in cluster algebra theory.

\subsection{Laurent phenomenon}
\label{1subsec:Laurent1}

We prove the Laurent phenomenon observed in rank 2 examples.
This is  the main result in CA1
and also the most fundamental fact on cluster algebras.

\begin{thm}[{Laurent phenomenon \cite{Fomin02}}]
\index{Laurent phenomenon}
\label{1thm:Laurent1}
Let $\bfSigma$ be any  cluster pattern with coefficients in any semifield $\bbP$.
Let $t_0,t\in \bbT_n$ be any vertices.
Then, any cluster variable $x_{i;t}$ is expressed as a Laurent polynomial in 
$\bfx_{t_0}$ with coefficients in $\bbZ\bbP$.
\end{thm}

Here we present the ``classic proof'' in CA1.
In view of Proposition \ref{1prop:free1},
 it is enough to prove Theorem \ref{1thm:Laurent1} for 
 any cluster pattern $\bfSigma$ with \emph{free coefficients at $t_0$}.
 This is because 
 the Laurent polynomial expression of $x_{i;t}$ in $\bfx_{t_0}$
 for $\bfSigma$ reduces to
 the one for any choice of  coefficients
under the specialization
in  Proposition \ref{1prop:free1}.

\begin{defn}[Coprime]
\index{coprime!for Laurent polynomials}
We say that two Laurent polynomials in $\bbZ\bbP[\bfx_{t_0}^{\pm1}]$ are
\emph{coprime\/} if there is no common factor except for Laurent monomials
in $\bfx_{t_0}$ with coefficients in $\bbZ\bbP^{\times}=\{\pm1\}\bbP$.
\end{defn}

From now on, we assume that $\bfSigma$ has free coefficients at $t_0$,
that is, $\bbP=\bbQ_{\rmsf}(\bfy)$ and $\bfy_{t_0}=\bfy$.
The proof in CA1 relies on the  following lemma.
(The assumption of coefficients is necessary only for the claim (b).)

\begin{lem}
\label{1lem:coprime1}
Let $t_1, t_2, t_3\in \bbT_n$ be vertices that are sequentially adjacent
to $t_0$ in the following way,
where $k \neq \ell$:
\vskip-5pt
\begin{align}
\label{1eq:subtree0}
\begin{picture}(90,0)(40,35)
\put(100,40){\circle*{3}}
\put(70,40){\circle*{3}}
\put(40,40){\circle*{3}}
\put(130,40){\circle*{3}}
\put(40,40){\line(1,0){90}}
\put(38,30){$t_0$}
\put(68,30){$t_1$}
\put(98,30){$t_2$}
\put(128,30){$t_3$}
\put(54,45){$k$}
\put(84,45){$\ell$}
\put(114,45){$k$}
\end{picture}
\end{align}
Then, the following facts hold:
\par
(a). The cluster variable $x_{k;t_3}$ can be expressed as a Laurent polynomial in
$\bfx_{t_0}$ with coefficients in $\bbZ\bbP$.
\par
(b). As elements in $\bbZ\bbP[\bfx_{t_0}^{\pm1}]$,
$x_{k;t_1}$ is coprime with  $x_{k;t_3}$  and $x_{\ell;t_3}$.
\end{lem}

Let us temporarily assume Lemma \ref{1lem:coprime1}
and prove Theorem \ref{1thm:Laurent1}.

Let $d(t,t')$ denote the distance between $t$ and $t'$ in $\bbT_n$,
that is, the number of edges of $\bbT_n$ between $t$ and $t'$.

\begin{proof}[Proof of Theorem \ref{1thm:Laurent1}]
Let $d=d(t_0,t)$. Let $t_1$, $t_2$, $t_3$ be the ones in \eqref{1eq:subtree0}.
Recall that
\begin{align}
\label{1eq:xmutrecall1}
x_{i;t_1}
&=
\begin{cases}
\displaystyle
x_{k;t_0}^{-1}\Biggl(\, \prod_{j=1}^n x_{j;t_0}^{[-b_{jk;t_0}]_+}
\Biggr)
\frac{ 1+\hat{y}_{k;t_0}}{ 1\oplus y_{k;t_0}}
& i=k,
\\
x_{i;t_0}
&i\neq k.
\end{cases}
\end{align}
They are certainly Laurent polynomials in $\bfx_{t_0}$.
Thus, the claim  holds for $d=1$.
The case $d=2$ is similar, where only $x_{\ell;t_2}$ is concerned.
For $d=3$, the problem arises 
for $x_{k;t_3}$, because 
the term $x_{k;t_2}=x_{k;t_1}$ in the denominator is no longer  a monomial in $\bfx_{t_0}$.
However, this case is covered by Lemma \ref{1lem:coprime1} (a).

We prove the claim by the induction on $d=d(t_0,t)\geq 3$,
\emph{where we fix $t$ and vary $t_0$}.
Assume that the claim hold up to $d$.
We consider the  situation in the following graph,
where $d(t_1,t)=d(t_3,t)=d$ and $d(t_0,t)=d+1$:
\vskip0pt
\begin{align}
\label{eq:subtree1}
\raisebox{-22pt}{
\begin{picture}(180,40)(40,0)
\put(100,40){\circle*{3}}
\put(100,10){\circle*{3}}
\put(70,40){\circle*{3}}
\put(40,40){\circle*{3}}
\put(130,40){\circle*{3}}
\put(160,40){\circle*{3}}
\put(190,40){\circle*{3}}
\put(220,40){\circle*{3}}
\put(40,40){\line(1,0){120}}
\put(190,40){\line(1,0){30}}
\put(100,40){\line(0,-1){30}}
\put(170,40){\circle*{1}}
\put(175,40){\circle*{1}}
\put(180,40){\circle*{1}}
\put(38,30){$t_0$}
\put(68,30){$t_1$}
\put(102,30){$t_2$}
\put(103,7){$t_3$}
\put(218,30){$t$}
\put(54,45){$k$}
\put(84,45){$\ell$}
\put(88,22){$k$}
\end{picture}
}
\end{align}
By the induction assumption,
$x_{i;t}$ is expressed as a Laurent polynomial in $\bfx_{t_1}$.
We write
\begin{align}
\label{1eq:xL1}
x_{i;t}=x_{k;t_1}^{-a} f(\bfx_{t_1}),
\quad
 f(\bfx_{t_1})\in \bbZ\bbP[\bfx_{t_1}^{\pm1}],
\end{align}
where $a\geq 0 $ is a sufficiently large integer
such that $f(\bfx_{t_1})$ does not contain any negative power of $x_{k;t_1}$.
Substituting the expression  \eqref{1eq:xmutrecall1} for  $x_{k;t_1}$ in 
$ f(\bfx_{t_1})$, we have
\begin{align}
\label{1eq:xL3}
x_{i;t}=x_{k;t_1}^{-a} \tilde{f}(\bfx_{t_0}),
\quad
\tilde{f}(\bfx_{t_0})\in \bbZ\bbP[\bfx_{t_0}^{\pm1}].
\end{align}
Meanwhile,
$x_{i;t}$ is also expressed  as a Laurent polynomial in $\bfx_{t_3}$.
Similarly, we write
\begin{align}
x_{i;t}=x_{k;t_3}^{-b}x_{\ell;t_3}^{-c} g(\bfx_{t_3}),
\quad
 g(\bfx_{t_3})\in \bbZ\bbP[\bfx_{t_3}^{\pm1}],
\end{align}
where $b,c \geq 0 $ are  sufficiently large integers
such that $g(\bfx_{t_3})$ does not contain any negative powers of $x_{k;t_3}$
and $x_{\ell;t_3}$.
Recall that $x_{k;t_3}$ and $x_{\ell;t_3}$ can be expressed as  Laurent polynomials 
in $\bfx_{t_0}$.
Substituting these expressions  for   $x_{k;t_3}$ and $x_{\ell;t_3}$ in 
$g(\bfx_{t_3})$,
we have
\begin{align}
\label{1eq:xL2}
x_{i;t}=x_{k;t_3}^{-b}x_{\ell;t_3}^{-c} \tilde{g}(\bfx_{t_0}),
\quad
 \tilde{g}(\bfx_{t_0})\in \bbZ\bbP[\bfx_{t_0}^{\pm1}].
\end{align}
Comparing \eqref{1eq:xL1} and \eqref{1eq:xL2}, we obtain the equality
\begin{align}
x_{k;t_3}^{b}x_{\ell;t_3}^{c}\tilde{f}(\bfx_{t_0})=x_{k;t_1}^{a}   \tilde{g}(\bfx_{t_0})
\quad \text{in $\bbZ\bbP[\bfx_{t_0}^{\pm1}]$}.
\end{align}
Then, by Lemma \ref{1lem:coprime1} (b), 
$\tilde{f}(\bfx_{t_0})$ is divisible by $x_{k;t_1}^{a} $ in 
$\bbZ\bbP[\bfx_{t_0}^{\pm1}]$.
Thus, by \eqref{1eq:xL3}, we conclude that
$x_{i;t}\in \bbZ\bbP[\bfx_{t_0}^{\pm1}]$.
\end{proof}

\begin{rem}
The above method  is commonly applicable to various systems showing the Laurent phenomenon beyond cluster algebras \cite{Fomin01}.
\end{rem}

Let us go back and prove Lemma \ref{1lem:coprime1}.

\begin{proof}[Proof of Lemma \ref{1lem:coprime1}]
Below we temporarily view any cluster variable as an element in $\bbQ\bbP(\bfx_{t_0})$.
For $a,b\in \bbQ\bbP(\bfx_{t_0})$,
we write $a\sim b$ if there is a  Laurent monomial $m$ in $\bfx_{t_0}$ with
 coefficients in $\bbZ\bbP^{\times}$ such that $a=mb$.

(a).
We need to compare two cluster variables
\begin{align}
x_{k;t_1}&=x_{k;t_0}^{-1}\Biggl(\, \prod_{j=1}^n x_{j;t_0}^{[-b_{jk;t_0}]_+}
\Biggr)
\frac{ 1+\hat{y}_{k;t_0}}{ 1\oplus y_{k;t_0}},
\\
x_{k;t_3}&=x_{k;t_2}^{-1}\Biggl(\, \prod_{j=1}^n x_{j;t_2}^{[-\varepsilon b_{jk;t_2}]_+}
\Biggr)
\frac{ 1+\hat{y}_{k;t_2}^\varepsilon }{ 1\oplus y_{k;t_2}^\varepsilon },
\end{align}
where for $x_{k;t_3}$ we employ the $\varepsilon$-expression
in Proposition \ref{1prop:epsilon1}.
Recall that $x_{k;t_2}=x_{k;t_1}$ and $x_{i;t_2}=x_{i;t_0}$ for $i\neq k,\ell$.
Thus, 
 we have
\begin{align}
x_{k;t_1}&\sim 1+\hat{y}_{k;t_0}.
\end{align}
If $b_{\ell k;t_0}\neq 0$,
we set
\begin{align}
\label{1eq:epsilon2}
\varepsilon = \mathrm{sign}(b_{\ell k;t_0})
=
- \mathrm{sign}(b_{\ell k;t_1})
=
 \mathrm{sign}(b_{\ell k;t_2})
 \in \{1, -1\},
\end{align}
where $\mathrm{sign}(a)=1$ if $a>0$,
and $-1$ if $a<0$.
If $b_{\ell k;t_0}=0$,  we choose $\varepsilon=\pm1$ arbitrarily.
Then, $[-\varepsilon b_{\ell k;t_2}]_+=0$, so that we have
\begin{align}
\label{1eq:yhatex4}
x_{k;t_3}&
\sim x_{\ell;t_2}^{[-\varepsilon b_{\ell k;t_2}]_+}
\frac{ 1+\hat{y}_{k;t_2}^\varepsilon}{x_{k;t_1}}
\sim 
\frac{ 1+\hat{y}_{k;t_2}^\varepsilon}{1+\hat{y}_{k;t_0}}.
\end{align}
Meanwhile, by Proposition \ref{1prop:epsilon1}, with the same $\varepsilon$
in \eqref{1eq:epsilon2},
we have
\begin{align}
\label{1eq:yhatex6}
\hat{y}_{\ell;t_1}&=
\hat{y}_{\ell;t_0}
\haty_{k;t_0}^{[\varepsilon b_{k\ell;t_0}]_+}
(1+\haty_{k;t_0}^{\varepsilon})^{-b_{k\ell;t_0}}
=
\hat{y}_{\ell;t_0}(1+\haty_{k;t_0}^{\varepsilon})^{- b_{k\ell;t_0}},
\\
\label{1eq:yhatex1}
\hat{y}_{k;t_2}&=\hat{y}_{k;t_1}\hat{y}_{\ell;t_1}^{[\varepsilon b_{\ell k;t_1}]_+}
(1+\hat{y}_{\ell;t_1}^{\varepsilon})^{-b_{\ell k;t_1}}
=
\hat{y}_{k;t_0}^{-1}(1+\hat{y}_{\ell;t_1}^{\varepsilon})^{-b_{\ell k;t_1}}.
\end{align}
We note that $-\varepsilon b_{k\ell;t_0}$ and $-\varepsilon b_{\ell k ;t_1}$
are both nonnegative.
Thus, ${ 1+\hat{y}_{k;t_2}^\varepsilon}\in \bbZ\bbP[\bfx_{t_0}^{\pm1}]$.
Therefore, it is enough to prove that $1+\hat{y}_{k;t_2}^\varepsilon$ is divisible
by $ 1+\hat{y}_{k;t_0}$ in $\bbZ\bbP[\bfx_{t_0}^{\pm1}]$.
If $b_{\ell k;t_0}= 0$, we have $\hat{y}_{k;t_2}^\varepsilon
=\hat{y}_{k;t_0}^{-\varepsilon}$
by \eqref{1eq:yhatex1}.
Therefore, the claim holds.
Suppose that $b_{\ell k;t_0}\neq 0$.
Let $I$ be the
 ideal  of $\bbZ\bbP[\bfx_{t_0}^{\pm1}]$
generated by $ 1+\hat{y}_{k;t_0}$.
By \eqref{1eq:yhatex6},
we have
$\hat{y}_{\ell;t_1}^{\varepsilon}\equiv 0\mod I$.
Thus, by \eqref{1eq:yhatex1}, we have
$\hat{y}_{k;t_2}^\varepsilon
\equiv 
\hat{y}_{k;t_0}^{-\varepsilon}
\mod I
$.
Therefore,  we have
\begin{align}
\label{1eq:yhatex2}
1+\hat{y}_{k;t_2}^\varepsilon
\equiv
1+\hat{y}_{k;t_0}^{-\varepsilon}
=
\hat{y}_{k;t_0}^{-\varepsilon}
(1+\hat{y}_{k;t_0}^{\varepsilon}
)
\equiv 0
\mod I.
\end{align}
\par
(b).
 First we  note that
\begin{align}
x_{\ell;t_3}=x_{\ell;t_2}\sim 
 x_{k;t_1}^{[-\varepsilon  b_{k \ell ;t_1}]_+}
(  1+\hat{y}_{\ell;t_1}^{\varepsilon})
=
 1+\hat{y}_{\ell;t_1}^{\varepsilon},
\end{align}
where $\varepsilon$ is the same sign in \eqref{1eq:epsilon2}.
Recall that by the assumption of free coefficients,
${y}_{\ell;t_0}=y_{\ell}$ and ${y}_{k;t_0}=y_k$ are algebraically independent.
By \eqref{1eq:yhatex6} and \eqref{1eq:yhatex1},
the following facts hold:
\begin{itemize}
\item
$ 1+\hat{y}_{k;t_0}$ is a constant with respect to $y_{\ell}^{\varepsilon}$.
\item
$1+\hat{y}_{\ell;t_1}^{\varepsilon}$ 
 is a binomial with respect to 
 ${y}_{\ell}^{\varepsilon}$
whose constant term is 1.
\item
$1+\hat{y}_{k;t_2}^{\varepsilon}$ is a polynomial  with respect to 
${y}_{\ell}^{\varepsilon}$
whose constant term is $1+\hat{y}_{k;t_0}^{-\varepsilon}$
\end{itemize}
It follows that
$x_{k;t_1}\sim 1+\hat{y}_{k;t_0}$ and $x_{\ell;t_3}\sim 1+\hat{y}_{\ell;t_1}^{\varepsilon}$ are coprime.
Also, 
it follows that
 $(1+\hat{y}_{k;t_2}^{\varepsilon})/(1+\hat{y}_{k;t_0}) $, which is in $\bbZ\bbP[\bfx_{t_0}^{\pm1}]$ by (a),
 is a polynomial with respect to 
${y}_{\ell}^{\varepsilon}$
whose constant term is $\hat{y}_{k;t_0}^{-1}$ if $\varepsilon=1$
and $1$ if $\varepsilon=-1$.
Therefore,
 $x_{k;t_3}\sim (1+\hat{y}_{k;t_2}^{\varepsilon})/(1+\hat{y}_{k;t_0})$ and $x_{k;t_1}\sim  1+\hat{y}_{k;t_0}$  are coprime.
\end{proof}

\subsection{Finite type classification}
\label{1subsec:finite1}

We present the finite type classification
of cluster algebras and cluster patterns without proofs.
This is the main result in CA2.

To state the result,
we briefly explain the background in Lie theory. 

The following definition is the counterpart of Definition \ref{1defn:skew1}.

\begin{defn}[Symmetrizable matrix]
\index{matrix!symmetrizable}
\label{1defn:symmetrizalbe1}
An integer square matrix $A=(a_{ij})_{i,j=1}^n$
is said to be \emph{symmetrizable\/}
if there is a  diagonal matrix $D=(d_i\delta_{ij})_{i,j=1}^n$
whose diagonal entries $d_i$ are positive integers
such that $DA$ is symmetric,
i.e,
\begin{align}
\label{1eq:ss2}
d_i a_{ij}= d_ja_{ji}.
\end{align}
The matrix $D$ is called a \emph{(left) symmetrizer\/} of $A$.
In particular, any symmetric matrix is symmetrizable.
\end{defn}

\begin{defn}[Cartan matrix] 
\index{Cartan matrix}
An $n\times n$ (integer) square matrix $A=(a_{ij})_{i,j=1}^n$
is called a \emph{(generalized) Cartan matrix\/} if
it satisfies the following conditions:
\begin{itemize}
\item For any $i$, we have $a_{ii}=2$.
\item For any $i,j$ ($i\neq j$),
we have $a_{ij}\leq 0$; moreover, $a_{ij}< 0$ if and only if $a_{ji}<0$.
\end{itemize}
\end{defn}

For each symmetrizable Cartan matrix $X$, one can define the associated
 Kac-Moody algebra, root system, and Weyl group \cite{Kac90}.

There is a natural (many-to-one) correspondence from a skew-symmetriz\-able matrix to a symmetrizable Cartan matrix.

\begin{defn}[Cartan counterpart]
\index{Cartan counterpart}
\label{1defn:Cartan1}
With any skew-symmetrizable matrix $B=(b_{ij})_{i,j=1}^n$,
we associate a symmetrizable Cartan matrix $A=A(B)$ as
\begin{align}
\label{1eq:AB1}
a_{ij}=
\begin{cases}
2 & i=j,\\
-|b_{ij}| & i\neq j.
\end{cases}
\end{align}
The matrix $A(B)$ is called the \emph{Cartan counterpart\/} of $B$.
\end{defn}

\begin{ex}
For the initial exchange matrices $B_0$
of type $A_2$, $B_2$, $G_2$ in Section \ref{1sec:rank2},
the associated Cartan matrices $A(B_0)$ are given as follows:
\begin{align}
\begin{pmatrix}
2 & -1\\
-1 & 2\\
\end{pmatrix},
\quad
\begin{pmatrix}
2 & -1\\
-2 & 2\\
\end{pmatrix},
\quad
\begin{pmatrix}
2 & -1\\
-3 & 2\\
\end{pmatrix}.
\end{align}
They are the Cartan matrices of type $A_2$, $B_2$, $G_2$, respectively,
as  explained below.
\end{ex}

We recall the Dynkin diagrams of finite type,
which is well-known in the context of semisimple Lie algebras and
crystallographic root systems \cite{Bourbaki68, Humphreys90}.

\begin{figure}
\centering
\begin{picture}(283,185)(-15,-175)
%
% A_r
\put(0,0){\circle{4}}
\put(20,0){\circle{4}}
\put(80,0){\circle{4}}
\put(100,0){\circle{4}}
\put(45,0){\circle*{1}}
\put(50,0){\circle*{1}}
\put(55,0){\circle*{1}}
\put(2,0){\line(1,0){16}}
\put(22,0){\line(1,0){16}}
\put(62,0){\line(1,0){16}}
\put(82,0){\line(1,0){16}}
\put(-30,-2){$A_n$}
\put(-2,-13){\small $1$}
\put(18,-13){\small $2$}
\put(70,-13){\small $ n-1$}
\put(98,-13){\small $n$}
%
% B_r
\put(180,0){
\put(0,0){\circle{4}}
\put(20,0){\circle{4}}
\put(80,0){\circle{4}}
\put(100,0){\circle{4}}
\put(45,0){\circle*{1}}
\put(50,0){\circle*{1}}
\put(55,0){\circle*{1}}
\put(2,0){\line(1,0){16}}
\put(22,0){\line(1,0){16}}
\put(62,0){\line(1,0){16}}
\put(81,-2){\line(1,0){18}}
\put(81,2){\line(1,0){18}}
\put(93,0){\line(-1,1){6}}
\put(93,0){\line(-1,-1){6}}
\put(-30,-2){$B_n$}
\put(-2,-13){\small $1$}
\put(18,-13){\small $2$}
\put(70,-13){\small $ n-1$}
\put(98,-13){\small $n$}
}
%
% C_r

\put(0,-40){
\put(0,0){\circle{4}}
\put(20,0){\circle{4}}
\put(80,0){\circle{4}}
\put(100,0){\circle{4}}
\put(45,0){\circle*{1}}
\put(50,0){\circle*{1}}
\put(55,0){\circle*{1}}
\put(2,0){\line(1,0){16}}
\put(22,0){\line(1,0){16}}
\put(62,0){\line(1,0){16}}
\put(81,-2){\line(1,0){18}}
\put(81,2){\line(1,0){18}}
\put(87,0){\line(1,1){6}}
\put(87,0){\line(1,-1){6}}
\put(-30,-2){$C_n$}
\put(-2,-13){\small $1$}
\put(18,-13){\small $2$}
\put(70,-13){\small $ n-1$}
\put(98,-13){\small $n$}
}

%
% D_r
\put(180,-40){
\put(0,0){\circle{4}}
\put(20,0){\circle{4}}
\put(80,0){\circle{4}}
\put(100,10){\circle{4}}
\put(100,-10){\circle{4}}
\put(45,0){\circle*{1}}
\put(50,0){\circle*{1}}
\put(55,0){\circle*{1}}
\put(2,0){\line(1,0){16}}
\put(22,0){\line(1,0){16}}
\put(62,0){\line(1,0){16}}
\put(81.4,0.7){\line(2,1){17}}
\put(81.4,-0.7){\line(2,-1){17}}
\put(-30,-2){$D_n$}
\put(-2,-13){\small $1$}
\put(18,-13){\small $2$}
\put(83,15){\small $ n-1$}
\put(70,-13){\small $ n-2$}
\put(98,-20){\small $n$}
}
%
% E_6
\put(0,-80){
\put(0,0){\circle{4}}
\put(20,0){\circle{4}}
\put(40,0){\circle{4}}
\put(60,0){\circle{4}}
\put(80,0){\circle{4}}
\put(40,20){\circle{4}}
\put(2,0){\line(1,0){16}}
\put(22,0){\line(1,0){16}}
\put(42,0){\line(1,0){16}}
\put(62,0){\line(1,0){16}}
\put(40,2){\line(0,1){16}}
\put(-30,-2){$E_6$}
\put(-2,-13){\small $1$}
\put(18,-13){\small $2$}
\put(38,-13){\small $3$}
\put(58,-13){\small $5$}
\put(78,-13){\small $6$}
\put(47,18){\small $4$}
}
%
% E_7
\put(180,-80){
\put(0,0){\circle{4}}
\put(20,0){\circle{4}}
\put(40,0){\circle{4}}
\put(60,0){\circle{4}}
\put(80,0){\circle{4}}
\put(100,0){\circle{4}}
\put(40,20){\circle{4}}
\put(2,0){\line(1,0){16}}
\put(22,0){\line(1,0){16}}
\put(42,0){\line(1,0){16}}
\put(62,0){\line(1,0){16}}
\put(82,0){\line(1,0){16}}
\put(40,2){\line(0,1){16}}
\put(-30,-2){$E_7$}
\put(-2,-13){\small $1$}
\put(18,-13){\small $2$}
\put(38,-13){\small $3$}
\put(58,-13){\small $4$}
\put(78,-13){\small $5$}
\put(98,-13){\small $6$}
\put(47,18){\small $7$}
}
%
% E_8
\put(0,-120){
\put(0,0){\circle{4}}
\put(20,0){\circle{4}}
\put(40,0){\circle{4}}
\put(60,0){\circle{4}}
\put(80,0){\circle{4}}
\put(100,0){\circle{4}}
\put(120,0){\circle{4}}
\put(80,20){\circle{4}}
\put(2,0){\line(1,0){16}}
\put(22,0){\line(1,0){16}}
\put(42,0){\line(1,0){16}}
\put(62,0){\line(1,0){16}}
\put(82,0){\line(1,0){16}}
\put(102,0){\line(1,0){16}}
\put(80,2){\line(0,1){16}}
\put(-30,-2){$E_8$}
\put(-2,-13){\small $1$}
\put(18,-13){\small $2$}
\put(38,-13){\small $3$}
\put(58,-13){\small $4$}
\put(78,-13){\small $5$}
\put(98,-13){\small $6$}
\put(118,-13){\small $7$}
\put(87,18){\small $8$}
}
%
% F_4
\put(0,-160){
\put(0,0){\circle{4}}
\put(20,0){\circle{4}}
\put(40,0){\circle{4}}
\put(60,0){\circle{4}}
\put(2,0){\line(1,0){16}}
\put(42,0){\line(1,0){16}}
\put(21,-2){\line(1,0){18}}
\put(21,2){\line(1,0){18}}
\put(33,0){\line(-1,1){6}}
\put(33,0){\line(-1,-1){6}}
\put(-30,-2){$F_4$}
\put(-2,-13){\small $1$}
\put(18,-13){\small $2$}
\put(38,-13){\small $3$}
\put(58,-13){\small $4$}
}
%
% G_2
\put(180,-160){
\put(0,0){\circle{4}}
\put(20,0){\circle{4}}
\put(1,-2){\line(1,0){18}}
\put(1,2){\line(1,0){18}}
\put(2,0){\line(1,0){16}}
\put(13,0){\line(-1,1){6}}
\put(13,0){\line(-1,-1){6}}
\put(-30,-2){$G_2$}
\put(-2,-13){\small $1$}
\put(18,-13){\small $2$}
}
\end{picture}
\vskip-10pt
\caption{Dynkin diagrams of finite type}
\label{1fig:Dynkin1}
\end{figure}
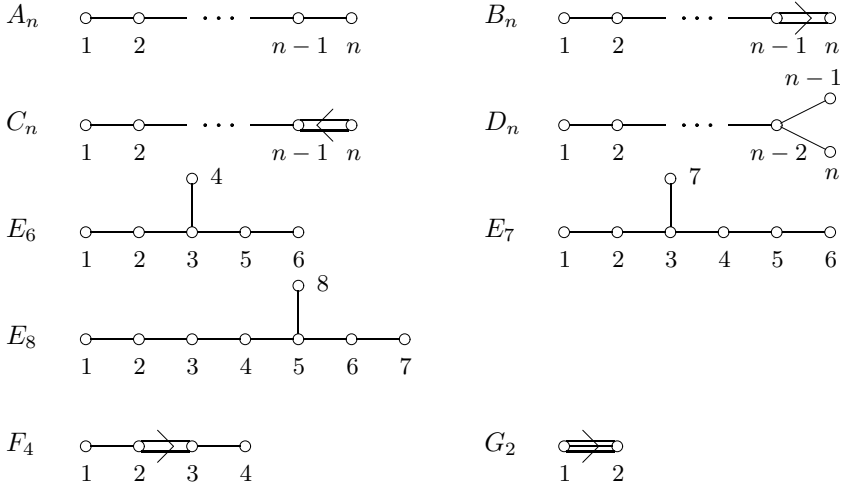

\begin{defn}
\label{1defn:Cartan2}
We call the graphs given in the list in Figure \ref{1fig:Dynkin1}
the \emph{Dynkin diagrams of finite type}.
\index{Dynkin diagram}
Here, $A_n$ ($n\geq 1$), $B_n$ ($n\geq 2$),
$C_n$ ($n\geq 2$), $D_n$ ($n\geq 4$),
and we identify $B_2$ and $C_2$ (up to the relabeling of the vertices).
\end{defn}
 
For each Dynkin diagram of type $X_n$ in Figure \ref{1fig:Dynkin1},
we define the $n\times n$ matrix $A=A(X_n)$ such that
$a_{ii}=2$, and nondiagonal entries are given as follows,
where we follow the convention in \cite{Kac90}, which is the transpose to 
the one  in \cite{Bourbaki68}
:
\begin{align}
\begin{cases}
a_{ij}=a_{ji}=0
&
\raisebox{4pt}{
\begin{picture}(20,10)(0,0)
\put(0,0){\circle{4}}
\put(20,0){\circle{4}}
\put(-2,-10){\small $i$}
\put(18,-10){\small $j$}
\end{picture}
}
\\
a_{ij}=a_{ji}=-1
&\raisebox{4pt}{
\begin{picture}(20,10)(0,0)
\put(0,0){\circle{4}}
\put(20,0){\circle{4}}
\put(2,0){\line(1,0){16}}
\put(-2,-10){\small $i$}
\put(18,-10){\small $j$}
\end{picture}
}
\\
a_{ij}=-1, \
a_{ji}=-2
&\raisebox{4pt}{
\begin{picture}(20,10)(0,0)
\put(0,0){\circle{4}}
\put(20,0){\circle{4}}
\put(1,-2){\line(1,0){18}}
\put(1,2){\line(1,0){18}}
\put(13,0){\line(-1,1){6}}
\put(13,0){\line(-1,-1){6}}
\put(-2,-10){\small $i$}
\put(18,-10){\small $j$}
\end{picture}
}
\\
a_{ij}=-1, \
a_{ji}=-3
&\raisebox{4pt}{
\begin{picture}(20,10)(0,0)
\put(0,0){\circle{4}}
\put(20,0){\circle{4}}
\put(1,-2){\line(1,0){18}}
\put(1,2){\line(1,0){18}}
\put(2,0){\line(1,0){16}}
\put(13,0){\line(-1,1){6}}
\put(13,0){\line(-1,-1){6}}
\put(-2,-10){\small $i$}
\put(18,-10){\small $j$}
\end{picture}
}
\end{cases}
\end{align}
The resulting matrix $A$ is a symmetrizable Cartan matrix.
These matrices (up to the simultaneous  permutations
of row and column indices) are called the \emph{Cartan matrices of finite type (of type $X_n$)}.
\index{Cartan matrix!of finite type}

\begin{ex}
The Cartan matrices of type $A_3$, $B_3$, $C_3$  (up to the simultaneous  permutations
of row and column indices) are given, respectively, by
\begin{align}
\left(
\begin{matrix}
2 & -1 & 0\\
-1 & 2 & -1\\
0 & -1 & 2
\end{matrix}
\right),
\quad
\left(
\begin{matrix}
2 & -1 & 0\\
-1 & 2 & -1\\
0 & -2 & 2
\end{matrix}
\right),
\quad
\left(
\begin{matrix}
2 & -1 & 0\\
-1 & 2 & -2\\
0 & -1 & 2
\end{matrix}
\right).
\end{align}
\end{ex}
For an indecomposable symmetrizable Cartan matrix $A$, it is known that
the following finiteness conditions are
equivalent \cite[Prop.~4.9]{Kac90}:
\begin{itemize}
\item[(a).]
$A$ is of finite type.
\item[(b).]
The Kac-Moody algebra associated with $A$ is finite-dimensional.
\item[(c).]
The root system associated with $A$ is a finite set.
\item[(d).]
The Weyl group associated with $A$ is a finite group.
\end{itemize}
Moreover, the isomorphism classes of such
Kac-Moody algebras, root systems, and Weyl groups,
respectively,
are classified by the Dynkin diagrams of finite type.
In addition, the above finiteness conditions
are equivalent to the following condition \cite[Prop.~4.7]{Kac90}:
\begin{itemize}
\item[(e).] All principal minors of $A$ are positive.
\end{itemize}

Now we consider the counterpart in cluster algebras.

\begin{defn}[Finite type]
A cluster pattern $\bfSigma$ is said to be \emph{of finite type\/}
if there are only finitely many distinct seeds 
of $\bfSigma$. \index{cluster pattern! of finite type}
The cluster algebra $\calA(\bfSigma)$ associated with a cluster pattern $\bfSigma$
 is said to be \emph{of finite type\/} if the cluster pattern $\bfSigma$
 is of finite type. \index{cluster algebra! of finite type}
\end{defn}

\begin{rem}
It is known \cite{Fomin03a} that
a cluster pattern $\bfSigma$ is of finite type
if and only if
  there are finitely many distinct \emph{clusters\/} 
  (equivalently, finitely many distinct \emph{cluster variables\/})
  of $\bfSigma$.
The only if part holds trivially,
while the if part requires a detailed study.
\end{rem}

\begin{defn}[Strongly isomorphic]
Let $\bfSigma$ and $\bfSigma'$ be 
cluster patterns of a common rank $n$ with a common coefficient semifield $\bbP$.
\begin{itemize}
\item
We say that cluster patterns $\bfSigma$ and $\bfSigma'$
are \emph{isomorphic\/} if there are some $t,t'\in \bbT_n$
and a permutation $\sigma\in S_n$ such that
$(\bfy'_{t'}, B'_{t'})=(\sigma(\bfy_t), \sigma(B_t))$.
\index{cluster pattern!isomorphic}
\item
We say that the cluster algebras $\calA(\bfSigma)$ and $\calA(\bfSigma')$
are \emph{strongly isomorphic\/} if
the underlying cluster patterns $\bfSigma$ and $\bfSigma'$
are isomorphic.
\index{cluster pattern!strongly isomorphic}
\end{itemize}
\end{defn}

\begin{prop}
If $\calA(\bfSigma)$ and $\calA(\bfSigma')$ are \emph{strongly isomorphic},
they are isomorphic 
as $\bbZ\bbP$-algebras.
\end{prop}
\begin{proof}
By the assumption, there are
some $t_0,t_1\in \bbT_n$
and a permutation $\sigma\in S_n$ such that
$(\bfy'_{t_1}, B'_{t_1})=(\sigma(\bfy_{t_0}), \sigma(B_{t_0}))$.
Let us consider the cluster pattern $\sigma\bfSigma=\{ \sigma\Sigma_t \}_{t\in \bbT_n}$.
Then, the correspondence $\bfx'_{t_1} \mapsto \sigma\bfx_{t_0}$ induces
a $\bbZ\bbP$-algebra isomorphism $\varphi:\calA(\bfSigma')\buildrel \sim \over \rightarrow
\calA(\sigma\bfSigma)$.
Meanwhile, thanks to Proposition \ref{1prop:compat1},
the set of cluster variables of $\sigma\bfSigma$ coincides 
with the one of $\bfSigma$.
Therefore, $\calA(\sigma\bfSigma)=\calA(\bfSigma)$.
\end{proof}

Now we are ready to state the classification of cluster patterns/algebras of
finite type.
In view of Proposition \ref{1prop:Bdecom2},
it is enough to concentrate on the cluster patterns
whose exchange matrices are indecomposable.

\begin{thm}[Finite type classification \cite{Fomin03a}]
\label{1thm:finite1}
For  cluster patterns with indecomposable exchange matrices,
the following facts hold.
\par
(a). The cluster pattern $\bfSigma$ is of finite type
if and only if there is some $t\in \bbT_n$ such that
the Cartan matrix $A(B_t)$ associated with $B_t$ is of finite type.
\par
(b). The Dynkin type of $A(B_t)$ does not depend on the choice of
such $t$.
\end{thm}

By Theorem \ref{1thm:finite1},
any isomorphism class of cluster patterns of finite type
with indecomposable exchange matrices
is uniquely labeled by the Dynkin type of the corresponding Cartan matrix $A(B_t)$ therein.
We have already used this labeling for the rank 2 examples in
Section \ref{1sec:rank2}.

For a cluster pattern of finite type,
the root system  of the corresponding type naturally
 parametrizes the cluster variables.

\begin{defn}[Denominator vector/Non-initial cluster variable] Let  $\bfSigma$ be any cluster pattern with a given initial vertex $t_0$.
\begin{itemize}
\item
For any cluster variable $x_{i;t}$,
we define an integer vector $\bfd_{i;t}=(d_{ji;t})_{j=1}^n$
such that $-d_{ji;t}$ is the lowest degree  of $x_{j;t_0}$
in the Laurent polynomial expression of $x_{i;t}$ in $\bfx_{t_0}$.
The vector $\bfd_{i;t}=(d_{ji;t})_{j=1}^n$ is called the \emph{denominator vector\/} (or \emph{$d$-vector}, for short)
 \index{vector!denominator ($d$-)}\index{$d$-vector}
 of $x_{i;t}$.
In particular, for an initial variable $x_{i;t_0}$, we have $\bfd_{i;t_0}=-\bfe_i$.
\item
We say that  a cluster variable $x_{i;t}$  is \emph{non-initial\/} \index{cluster variable!non-initial} if it does not coincide with
any initial cluster variables $x_{1;t_0}$, \dots, $x_{n;t_0}$. 
\end{itemize}
\end{defn}

\begin{ex}
In the examples in Section \ref{1sec:rank2},
we see the following list of the denominator vectors for the non-initial cluster variables.
\begin{align}
\mbox(A_2)&:
(1,0), (1,1), (0,1),
\\
\mbox(B_2)&:
(1,0), (1,1), (1,2), (0,1),
\\
\mbox(G_2)&:
(1,0), (1,1), (1,2), (1,3),(2,3), (0,1).
\end{align}
They are naturally identified with the \emph{positive roots\/} of the 
corresponding root systems.
See Figure \ref{1fig:root1}.

\begin{figure}
\centering
\leavevmode
\begin{xy}
(14,0)*{\mbox{\small$\alpha_1$}},
(6,12)*{\mbox{\small$\alpha_1+\alpha_2$}},
(-6.5,12)*{\mbox{\small$\alpha_2$}},
(0,-14)*{A_2}
\ar@{->} (0,0);(10,0) %
\ar@{->} (0,0);(5,8.66) %
\ar@{->} (0,0);(-5,8.66) %
\ar@{->} (0,0);(-10,0) %
\ar@{->} (0,0);(-5,-8.66) %
\ar@{->} (0,0);(5,-8.66) %
\end{xy}
\hskip43pt
\begin{xy}
(14,0)*{\mbox{\small$\alpha_1$}},
(0,12)*{\mbox{\small$\alpha_1+2\alpha_2$}},
(10,7.5)*{\mbox{\small$\alpha_1+\alpha_2$}},
(-7.5,7.5)*{\mbox{\small$\alpha_2$}},
(0,-14)*{B_2}
\ar@{->} (0,0);(10,0) %
\ar@{->} (0,0);(5,5) %
\ar@{->} (0,0);(0,10) %
\ar@{->} (0,0);(-5, 5) %
\ar@{->} (0,0);(-10,0) %
\ar@{->} (0,0);(-5,-5) %
\ar@{->} (0,0);(0,-10) %
\ar@{->} (0,0);(5,-5) %
\end{xy}
\qquad
\begin{xy}
(14,0)*{\mbox{\small$\alpha_1$}},
(-12,8)*{\mbox{\small$\alpha_1+3\alpha_2$}},
(13,8)*{\mbox{\small$2\alpha_1+3\alpha_2$}},
(0,12)*{\mbox{\small$\alpha_1+2\alpha_2$}},
(12,4)*{\mbox{\small$\alpha_1+ \alpha_2$}},
(-8.3,4)*{\mbox{\small$\alpha_2$}},
(0,-14)*{G_2}
\ar@{->} (0,0);(10,0) %
\ar@{->} (0,0);(5,8.66) %
\ar@{->} (0,0);(-5,8.66) %
\ar@{->} (0,0);(-10,0) %
\ar@{->} (0,0);(-5,-8.66) %
\ar@{->} (0,0);(5,-8.66) %
\ar@{->} (0,0);(5,2.89) %
\ar@{->} (0,0);(5,-2.89) %
\ar@{->} (0,0);(-5,2.89) %
\ar@{->} (0,0);(-5,-2.89) %
\ar@{->} (0,0);(0,5.77) %
\ar@{->} (0,0);(0,-5.77) %
\end{xy}
\caption{Root systems of rank 2.}
\label{1fig:root1}
\end{figure}
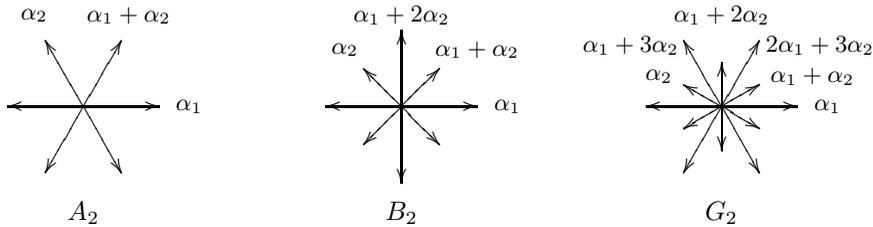
\end{ex}

This phenomenon can be fully generalized to any cluster pattern of finite type
if we properly choose the initial vertex $t_0$.

\begin{thm}[\cite{Fomin03a}]
Let $\bfSigma$ be any cluster pattern of finite type with a initial vertex
$t_0$ such that $A(B_{t_0})$ is a Cartan matrix of finite type.
Then, the denominator vector of any non-initial cluster variable
is identified with a positive root of the corresponding root system.
Moreover, this correspondence is one-to-one.
\end{thm}

\subsection{Cluster algebras of geometric type}
\label{1sec:geometric1}

Below we focus on a class of cluster patterns whose
coefficient semifields are \emph{tropical semifields}.
They are important in  various applications of cluster algebras.

Let $\mathrm{Trop}(\bfu)$ be the tropical semifield in  Example \ref{1ex:trop1}.

\begin{defn}[Geometric type]
\index{cluster pattern!of geometric type}\index{cluster algebra!of geometric type}
\label{1defn:geom1}
A cluster pattern (resp. cluster algebra) with coefficients
in a tropical semifield $\mathrm{Trop}(\bfu)$ is called
a cluster pattern (resp. cluster algebra) of \emph{geometric type},
where $\bfu=(u_1,\dots, u_m)$ and $m$ is taken independently of
the rank $n$ of a cluster pattern.
\end{defn}

\begin{rem}
The terminology originates in the fact
that  cluster algebras of this type typically arise as the coordinating rings
of certain algebraic varieties,
which are called the \emph{geometric realization\/} of cluster algebras
in \cite{Fomin03a}.
\end{rem}

Let $\bfSigma$ be a cluster pattern of geometric type.
Recall that
an element of $ \mathrm{Trop}(\bfu)$ is a Laurent monomial in $\bfu$
with  coefficient 1.
Therefore, any coefficient $y_{i;t}$ of $\bfSigma$ at $t\in \bbT_n$ is represented as
\begin{align}
\label{eq:cexp1}
y_{i;t}=\prod_{j=1}^m u_j^{c_{ji;t}},
\quad
i=1,\dots,n.
\end{align}
This determines an $m\times n$ integer matrix $C_t=(c_{ij;t})$
for each $t \in \mathbb{T}_n$,
and  the coefficient tuple $\bfy_t$ and the matrix $C_t$
are identified.
Under this identification, the mutation in \eqref{1eq:ymut1}
is translated as follows.

\begin{prop}
Let $t,t'\in \bbT_n$ be $k$-adjacent vertices. 
Then, we have the following mutation formula of matrices $C_t$:
\begin{align}
\label{1eq:cmut1}
c_{ij;t'}=
\begin{cases}
-c_{ik;t} & j = k,\\
c_{ij;t} +  c_{ik;t}[b_{kj;t}]_+  + [-c_{ik;t}]_+b_{kj;t} 
& j \neq k.
\end{cases}
\end{align}
\end{prop}
\begin{proof}
For any integer $a$, the following identity holds:
\begin{align}
\min(0,a)=-[-a]_+.
\end{align}
Thus, by the definition of the tropical sum in \eqref{1eq:ts1}, we have
\begin{align}
\label{1eq:1+y1}
{1\oplus y_{k;t}}=
 \prod_{i=1}^m u_i^{-[-c_{ik;t}]_+}.
\end{align}
By replacing the index $i$ in \eqref{1eq:ymut1} with $j$,
and comparing the power of $u_i$, we obtain
\eqref{1eq:cmut1}.
\end{proof}

Observe that the formula \eqref{1eq:cmut1}
is  parallel  to the mutation \eqref{1eq:bmut2}
of the exchange matrix $B_t$.
This motivates us to define the following $(n+m)\times n$
matrix
\begin{align}
\label{1eq:extended1}
\tilde{B_t}=
\left(
\begin{matrix}
B_t\\
C_t
\end{matrix}
\right),
\end{align}
which is called an \emph{extended exchange matrix\/} of $\bfSigma$.
\index{matrix!extended exchange}
Abusing the notation,
let us write the $(i,j)$-entry of $\tilde{B_t}$ as $b_{ij}$.
Then, we can unify the mutations  \eqref{1eq:bmut2}
and
\eqref{1eq:cmut1}
as a single matrix mutation of $\tilde{B_t}$ defined by
\begin{align}
\label{1eq:bmut5}
b_{ij;t'}=
\begin{cases}
-b_{ij;t} & \mbox{$i=k$ or $j=k$},\\
b_{ij;t}+ b_{ik;t}[b_{kj;t}]_+  + [-b_{ik;t}]_+b_{kj;t} 
& i,j\neq k,
\end{cases}
\end{align}
where $i=1,\dots,n+m$, and $j,k=1,\dots,n$.
Next,
let us observe that
the mutations of cluster variables
can be also expressed in terms of the matrix $\tilde{B}_t$.
By \eqref{1eq:1+y1}, we have
\begin{align}
\label{1eq:1+y2}
\frac{y_{k;t}}{1\oplus y_{k;t}}= \prod_{j=1}^m u_j^{[c_{jk;t}]_+},
\quad
\frac{1}{1\oplus y_{k;t}}= \prod_{j=1}^m u_j^{[-c_{jk;t}]_+}.
\end{align}
Thus,
the mutation in the form of \eqref{1eq:xmutstandard1}
is written as
\begin{align}
\label{1eq:xmut5}
x_{i;t'}&=
\begin{cases}
\displaystyle
\frac{1}{x_{k;t}}
\Biggl(\, \prod_{j=1}^m u_j^{[c_{jk;t}]_+}\prod_{j=1}^{n} x_{j;t}^{[b_{jk;t}]_+}
+ \prod_{j=1}^m u_j^{[-c_{jk;t}]_+}\prod_{j=1}^{n} x_{j;t}^{[-b_{jk;t}]_+}
\Biggr)
&
i=k,
\\
x_{i;t}& i \neq k.
\end{cases}
\end{align}
To summarize the result so far, the coefficients $\bfy_t$
of geometric type
can be safely replaced with
the lower part $C_t$ of the extended exchange matrices $\tilde{B_t}$.

To make the picture complete, 
we introduce an \emph{extended cluster}, \index{cluster!extended}
\begin{align}
\tilde\bfx_t=(x_{1;t},\dots,x_{n+m:t}):=(x_{1;t},\dots,x_{n;t}, u_1,\dots,u_m).
\end{align}
Then, together with the extended exchange matrix notation,
the mutation \eqref{1eq:xmut5} is written as,
 for $k=1,\dots,n$,
\begin{align}
\label{1eq:xmut6}
x_{i;t'}&=
\begin{cases}
\displaystyle
\frac{1}{x_{k;t}}
\Biggl(\,
\prod_{j=1}^{n+m} x_{j;t}^{[b_{jk;t}]_+}
+ \prod_{j=1}^{n+m} x_{j;t}^{[-b_{jk;t}]_+}
\Biggr)
&
i=k,
\\
x_{i;t}& i \neq k.
\end{cases}
\end{align}
This formally coincides with the mutation of cluster variables
\emph{without coefficients}.

Let us formulate this observation more completely in terms of
 a cluster pattern of rank $n+m$ \emph{without coefficients}.
We consider the extension of $\tilde{B}_{t_0}$ 
 to a full $(n+m)\times (n+m)$
skew-symmetrizable matrix.
Such an extension (together with its skew-symmetrizer) is not unique.
To given an example,
let $D$ be a  skew-symmetrizer of $B_{t_0}$,
and let $d$ be the least common multiple
of  the diagonal entries of $D$.
Then, the following matrix gives a skew-symmetrizable extension
of  $\tilde{B}_{t_0}$:
\begin{align}
\hat{B}_{t_0}=
\begin{pmatrix}
B_{t_0} & -D^{-1} d C_{t_0}^T\\
C_{t_0} & O
\end{pmatrix},
\end{align}
where  a skew-symmetrizer is given by $\hat{D}=D\oplus d I_m$.
(The factor  $d$ ensures that the matrix $\hat{B}_{t_0}$ is an integer matrix.)
For a given initial vertex $t_0$ in $\bbT_{n+m}$,  we regard $\bbT_n$
as a subtree of $\bbT_{n+m}$
such that $\bbT_n$ contains $t_0$ and all vertices 
that are sequentially adjacent to $t_0$ with edges labeled by $\{1, \dots, n\}$
in $\bbT_{n+m}$.
Then,  the above cluster pattern $\bfSigma$ of geometric type
 is equivalent to consider a cluster pattern of rank $n+m$
 \emph{without coefficients}, restricted to the subtree $\bbT_n \subset \bbT_{n+m}$,
 and with the  initial seed 
 $ (\tilde\bfx_{t_0}, \hat{B}_{t_0})$.
The elements $x_{n+1;t}=x_{n+1}$, \dots, $x_{n+m;t}=x_{n+m}$ are called the \emph{frozen variables},
\index{cluster variable!frozen variable}
because they are not mutated.
 We note that the right half of $\hat{B}_t$
 does not influence neither \eqref{1eq:bmut5}
 nor \eqref{1eq:xmut6}.
Thus, we may safely replace seeds $(\tilde\bfx_t, \hat{B}_t)$ with
$(\tilde\bfx_t, \tilde{B}_t)$ without losing information.

Let us summarize the results obtained above as a proposition:
\begin{prop}
\label{1prop:geometric1}
We have the following   three  equivalent presentations
of a cluster pattern of geometric type:
\begin{itemize}
\item
a cluster pattern  consisting seeds $(\bfx_t, \bfy_t, B_t)$
with coefficients in $\mathrm{Trop}(\bfu)$, by definition.
\item
a cluster pattern   consisting seeds $(\bfx_t, \tilde{B}_t)$
whose coefficients  in $\mathrm{Trop}(\bfu)$  are encoded in
the extend exchange matrices $\tilde{B}_t$.
\item
a cluster pattern   consisting extended seeds $(\tilde\bfx_t, \hat{B}_t)$
or $(\tilde\bfx_t, \tilde{B}_t)$
without coefficients and with frozen variables  $x_{n+1}$, \dots, $x_{n+m}$
that are identified with $u_1$, \dots, $u_m$.
(In this case we still say that it is of rank $n$.)
\end{itemize}
\end{prop}

Even if the above three pictures are equivalent as \emph{cluster patterns},
there is some discrepancy  for the corresponding \emph{cluster algebras\/}
due to the difference of the coefficient semifields.
Note that, for $\bbP=\mathrm{Trop}(\bfu)$, the group ring
$\bbZ\bbP$ is identified with the Laurent polynomial ring
$\bbZ[\bfu^{\pm1}]$.
Then, for the first two pictures,
the  corresponding cluster algebras are
given by
\begin{align}
\label{1eq:ext1}
 \calA=\bbZ[\bfu^{\pm1}][\bfx_t]_{t\in \bbT_n}
=
\bbZ[\bfx_t, \bfu^{\pm1}]_{t\in \bbT_n},
\end{align}
while, for the last one,
including the frozen variables as honorary members,
the  corresponding cluster algebra is
given by
\begin{align}
\label{1eq:geometric1}
 \calA=\bbZ[\tilde\bfx_t]_{t\in \bbT_n}
=
\bbZ[\bfx_t, \bfu]_{t\in \bbT_n}.
\end{align}
Therefore, one has to be careful in which sense cluster algebra is considered.

For a cluster pattern of geometric type,
the Laurent phenomenon in Theorem \ref{1thm:Laurent1}
claims
\begin{align}
x_{i;t} \in \bbZ[\bfu^{\pm1}][\bfx_{t_0}^{\pm1}] = \bbZ[\bfx_{t_0}^{\pm1}, \bfu^{\pm1}].
\end{align}
Actually,
 the following strong version of the
Laurent phenomenon holds.

\begin{thm}[{\cite{Fomin03a}}]
\label{1thm:Laurent2}
For any cluster pattern of geometric type, we have
\begin{align}
\label{1eq:Laurent1}
x_{i;t} \in 
\bbZ[\bfx_{t_0}^{\pm1}, \bfu].
\end{align}
\end{thm}

\begin{proof}
We consider the following claim, which is slightly stronger than the desired result.
\par\noindent
{\bf Claim.} For each $j$, $x_{i;t}$ is a polynomial in $u_j$
whose constant  term (with respect to $u_j$)  is a \emph{nonzero\/} polynomial
with a subtraction-free rational expression in  $\bfx_{t_0}^{\pm1}$
and $\bfu$ other than $u_j$.
(See Definition \ref{1ex:sf1} for a subtraction-free expression.)

We prove the claim by the induction on the distance $d=d(t_0, t)$ in $\bbT_n$.
For $d=0$, we have $x_{i;t}=x_i$, so that the claim holds.
 Suppose that the claim holds for $d=d(t_0,t)$.
Let $t'\in \bbT_n$ be $k$-adjacent to $t$.
It is enough to prove the claim for $x_{k;t'}$.
Let us look at the first case of \eqref{1eq:xmut5}.
Note that, for each $j$, either $[c_{jk;t}]_+$ or $[-c_{jk;t}]_+$ is zero.
Thus, by the induction assumption, the numerator in 
\eqref{1eq:xmut5} is a polynomial in $u_j$
whose constant term
 is a nonzero polynomial
 with a subtraction-free rational expression
in  $\bfx_{t_0}^{\pm1}$
and $\bfu$ other than $u_j$.
(Here,  the condition of having a subtraction-free expression 
guarantees 
that,  even in the special case $c_{jk;t}=0$, the constant terms from
the two products in \eqref{1eq:xmut5}  never cancel accidentally.)
Also,  the denominator $x_{k;t}$ has the same property.
Since we already know that $x_{k;t'}$ is in 
$\bbZ[\bfx_{t_0}^{\pm1}, \bfu^{\pm1}]$,
the denominator divides the numerator such that
the result is a polynomial in $u_j$
having the same property.
\end{proof}

\begin{rem}
The above claim in the proof does not imply that  $x_{i;t}$ has a nonzero constant term
in $\bfu$.
\end{rem}

\subsection{Grassmannian $\rmGr(2,5)$}

Let us present a prototypical example 
of a cluster algebra (of geometric type) appearing in Lie theory.

First, we briefly recall basic definitions/facts on Grassmannian $\rmGr(2,5)$.
See \cite[Section 9]{Fulton97} for the background and more information.

\begin{defn}[Grassmannian $\rmGr(2,5)$]
\index{Grassmannian}
\label{1defn:Grassmannian1}
\
\begin{itemize}
\item
In short, the \emph{Grassmannian $\rmGr(2,5)$}
is the complex projective variety consisting of all 2-dimensional subspaces in a
5-dimensional vector space $V$ over $\bbC$.
\item
An element of $\rmGr(2,5)$ is identified with
its basis, which is represented by
a $2\times 5$ complex matrix $M=(m_{ij})$ of rank 2 modulo the left action of $GL(2,\bbC)$.
This tells us that $\dim \rmGr(2,5)=10-4=6$.
\item For each pair $1\leq i < j\leq 5$,
let
\begin{align}
p_{ij}=p_{ij}(M):=
\left|
\begin{matrix}
m_{1i} & m_{1j}\\
m_{2i} & m_{2j} 
\end{matrix}
\right|.
\end{align}
Then, we define a map $i: \rmGr(2,5)\rightarrow \bbC\bbP^{9}$ by
\begin{align}
i: M \mapsto [p_{12}:p_{13}:\dots :p_{45}],
\end{align}
where $[z_1:\dots :z_{10}]$ is the homogeneous coordinates of 
the complex projective space $\bbC\bbP^{9}$.
The map is injective, and it is called the \emph{Pl\"ucker embedding},
while $p_{ij}$ are called the \emph{Pl\"ucker coordinates}.
 \index{Grassmannian!Pl\"ucker coordinates}
\item  For  $1\leq i<j<k<\ell\leq 5$, the Pl\"ucker coordinates satisfy the relations
\begin{align}
\label{1eq:PluR1}
R_{ijk\ell}:\ p_{ij}p_{k\ell}-p_{ik}p_{j\ell} +p_{i\ell}p_{j k}  =0.
\end{align}
 They are called the \emph{Pl\"ucker relations}.
 \index{Grassmannian!Pl\"ucker relation}
 Note that there are five Pl\"ucker relations, $R_{1234}$, $R_{1235}$, $R_{1245}$,
 $R_{1345}$, and $R_{2345}$.   
\item The homogeneous coordinate ring of $\rmGr(2,5)$
is given by
\begin{align}
\label{1eq:cr1}
\bbC[\rmGr(2,5)]=\bbC[\bfP]/I_R,
\quad \bfP=(p_{ij})_{1\leq i<j\leq 5},
\end{align}
where $I_R$ is 
 the homogeneous ideal generated by the Pl\"ucker relations;
 moreover, $I_R$ is prime.
 Thus, $\rmGr(2,5)$ is a projective variety in  $\bbC\bbP^{9}$.
 \item
 Alternatively, we may consider the \emph{affine cone} $\widehat{\rmGr}(2,5)$ over
 $\rmGr(2,5)$, that is,
 the affine variety in the complex affine space $\bbA^{10}$ defined by 
 the same  Pl\"ucker relations \eqref{1eq:PluR1}.
 The coordinate ring of  $\widehat{\rmGr}(2,5)$ has the same form as
 \eqref{1eq:cr1},
 \begin{align}
\label{1eq:cr2}
\bbC[\widehat{\rmGr}(2,5)]=\bbC[\bfP]/I_R.
\end{align}
However, the ring in  \eqref{1eq:cr1} is regarded as a graded ring,
while it is a (not graded) ring here.
\end{itemize}
\end{defn}

We are going to show that  $\bbC[{\rmGr}(2,5)]$
(or $\bbC[\widehat{\rmGr}(2,5)]$ by forgetting the grading)
has the cluster algebra structure.
To do that, we consider a cluster pattern of
geometric type, where there are
two cluster variables $x_{i;t}$ ($i=1,2$) and five frozen variables
$x_a$, \dots, $x_e$.
The cluster pattern  $\bfSigma$  we consider is encoded in the following initial 
extended exchange matrix $\tilde{B}_0$, or equivalently, the corresponding
initial \emph{extended quiver\/} $\tilde{Q}_0$, where the \index{quiver!extended}
vertices $a$, \dots, $e$ are  \emph{frozen vertices\/}
and  we will mutate  only at the vertices 1 and 2.

\begin{align}
\label{1eq:BQ1}
\tilde{B}_0
=
\begin{pmatrix}
0 & -1\\
1 & 0\\
\hline
-1 & 0\\
1 & 0\\
-1 & 1 \\
0 & -1 \\
0 & 1\\
\end{pmatrix},
\qquad
\tilde{Q}_0
=
\raisebox{-8pt}
{\mbox{
\begin{xy}
(0,-3)*{a},
(10,-3)*{1},
(20,-3)*{2},
(30,-3)*{e},
(5,11)*{b},
(15,11)*{c},
(25,11)*{d},
(0,0)*{\bullet},
(10,0)*\cir<2pt>{},
(20,0)*\cir<2pt>{},
(30,0)*{\bullet},
(5,8)*{\bullet},
(15,8)*{\bullet},
(25,8)*{\bullet},
\ar@{->} (9,0);(1,0)
\ar@{->} (19,0);(11,0)
\ar@{->} (29,0);(21,0)
\ar@{->} (5.5,7);(9.5,1)
\ar@{<-} (14.5,7);(10.5,1)
\ar@{->} (15.5,7);(19.5,1)
\ar@{<-} (24.5,7);(20.5,1)
\end{xy}
}}
\end{align}
In particular, this is a special case of cluster pattern of type $A_2$.

In addition, we introduce an alternative presentation of the above data
by a \emph{triangulation of a pentagon},
where the frozen vertices of the quiver correspond to the sides of the pentagon,
while the unfrozen vertices correspond to  two diagonals
as follows:
\begin{align}
\label{1eq:pentagon1}
T_0=
\raisebox{-20pt}
{\mbox{
\begin{xy}
(-6,18)*{a},
(6,18)*{b},
(10.5,5.5)*{c},
(-10.5,5)*{e},
(0,-3)*{d},
(0,14)*{1},
(0,7)*{2},
(-12,12)*{1},
(12,12)*{3},
(0,22)*{2},
(8,-1)*{4},
(-8,-1)*{5},
\ar@{-} (-6.5,0);(6.5,0)
\ar@{-} (6.5,0);(10.3,11.7)
\ar@{-} (0,20);(10.3,11.7)
\ar@{-} (-6.5,0);(-10.3,11.7)
\ar@{-} (0,20);(-10.3,11.7)
\ar@{-} (10.3,11.7);(-10.3,11.7)
\ar@{-} (6.5,0);(-10.3,11.7)
\end{xy}
}}
\end{align}
To explain in more detail, we identify a subquiver in $\tilde{Q}_0$ and a triangle in $T_0$ by the following  correspondence:
\begin{align}
\raisebox{-8pt}
{\mbox{
\begin{xy}
(10,-3)*{i},
(20,-3)*{j},
(15,11)*{k},
(10,0)*\cir<2pt>{},
(20,0)*\cir<2pt>{},
(15,8)*\cir<2pt>{},
(42,6)*{i},
(54,6)*{k},
(48,-5)*{j},
\ar@{->} (19,0);(11,0)
\ar@{<-} (14.5,7);(10.5,1)
\ar@{->} (15.5,7);(19.5,1)
\ar@{-} (40,-2);(56,-2)
\ar@{-} (48,11);(56,-2)
\ar@{-} (48,11);(40,-2)
\end{xy}
}}
\end{align}
where we regard that
the arrows among frozen vertices (i.e, $a\rightarrow b$ and
$d\rightarrow e$) are omitted
in  $\tilde{Q}_0$.
(This is reasonable, because they do not influence  the mutations of
 seeds $(\tilde\bfx_t, \tilde{B}_t)$ at all.)

The advantage of the presentation by  a  triangulation is that
one can manage a mutation of a seed pictorially as a \emph{flip\/} of a diagonal. 
For example, the mutation of the quiver $\tilde{Q}_0$ at the vertex 1
gives the following quiver $\tilde{Q}_1$,
where again the arrows among frozen vertices are omitted:
\begin{align}
%\label{1eq:Q2}
\tilde{Q}_1
=
\raisebox{-8pt}
{\mbox{
\begin{xy}
(0,-3)*{a},
(10,-3)*{1},
(20,-3)*{2},
(30,-3)*{e},
(5,11)*{b},
(15,11)*{c},
(25,11)*{d},
(0,0)*{\bullet},
(10,0)*\cir<2pt>{},
(20,0)*\cir<2pt>{},
(30,0)*{\bullet},
(5,8)*{\bullet},
(15,8)*{\bullet},
(25,8)*{\bullet},
(19,-1); (1,-1) **\crv{(10,-10)};
\ar@{<-} (9,0);(1,0)
\ar@{<-} (19,0);(11,0)
\ar@{->} (29,0);(21,0)
\ar@{<-} (5.5,7);(9.5,1)
\ar@{->} (14.5,7);(10.5,1)
\ar@{<-} (24.5,7);(20.5,1)
\ar@{->} (2,-1.9);(1,-1)
\end{xy}
}}
\end{align}
The corresponding triangulation $T_1$ is given by
the flip of the diagonal 1 inside the
quadrilateral with sides $a$, $b$, $c$, 2
depicted as below:
\begin{align}
T_1=
\raisebox{-20pt}
{\mbox{
\begin{xy}
%edge
(-6,18)*{a},
(6,18)*{b},
(10.5,5.5)*{c},
(-10.5,5)*{e},
(0,-3)*{d},
%diagonal
(5,10)*{1},
(0,7)*{2},
(-12,12)*{1},
(12,12)*{3},
(0,22)*{2},
(8,-1)*{4},
(-8,-1)*{5},
%edges
\ar@{-} (-6.5,0);(6.5,0)
\ar@{-} (6.5,0);(10.3,11.7)
\ar@{-} (0,20);(10.3,11.7)
\ar@{-} (-6.5,0);(-10.3,11.7)
\ar@{-} (0,20);(-10.3,11.7)
%diagonal
\ar@{-} (0,20);(6.5,0)
\ar@{-} (6.5,0);(-10.3,11.7)
\end{xy}
}}
\end{align}
Meanwhile, the mutation of the cluster variable $x_1$ is given by
\begin{align}
\label{1eq:pto1}
x'_1=\frac{x_b x_2+x_a x_c}{x_1}.
\end{align}
By rewriting it in the form
\begin{align}
\label{1eq:pto2}
x_1 x'_1={x_b x_2+x_a x_c},
\end{align}
we can recognize it as \emph{Ptolemy's theorem\/} for 
cyclic quadrilaterals.

\begin{table}
\begin{tabular}{ccccc}
$t$ &  $\tilde{B}_t$ & $\tilde{Q}_t$ & $T_t$& $x_{1;t}, x_{2;t}$
\\
\hline
0
&
\raisebox{0pt}[50pt]
{
$
{\small
\begin{pmatrix}
0 & -1\\
1 & 0\\
\hline
-1 & 0\\
1 & 0\\
-1 & 1 \\
0 & -1 \\
0 & 1\\
\end{pmatrix}
}
$
}
&
{\small
\begin{xy}
(0,-3)*{a},
(10,-3)*{1},
(20,-3)*{2},
(30,-3)*{e},
(5,11)*{b},
(15,11)*{c},
(25,11)*{d},
(0,0)*{\bullet},
(10,0)*\cir<2pt>{},
(20,0)*\cir<2pt>{},
(30,0)*{\bullet},
(5,8)*{\bullet},
(15,8)*{\bullet},
(25,8)*{\bullet},
\ar@{->} (9,0);(1,0)
\ar@{->} (19,0);(11,0)
\ar@{->} (29,0);(21,0)
\ar@{->} (5.5,7);(9.5,1)
\ar@{<-} (14.5,7);(10.5,1)
\ar@{->} (15.5,7);(19.5,1)
\ar@{<-} (24.5,7);(20.5,1)
\end{xy}
}
&
{
\small
\raisebox{-20pt}
{\begin{xy}
(-6,18)*{a},
(6,18)*{b},
(10.5,5.5)*{c},
(-10.5,5)*{e},
(0,-3)*{d},
(0,14)*{1},
(0,7)*{2},
(-12,12)*{1},
(12,12)*{3},
(0,22)*{2},
(8,-1)*{4},
(-8,-1)*{5},
\ar@{-} (-6.5,0);(6.5,0)
\ar@{-} (6.5,0);(10.3,11.7)
\ar@{-} (0,20);(10.3,11.7)
\ar@{-} (-6.5,0);(-10.3,11.7)
\ar@{-} (0,20);(-10.3,11.7)
\ar@{-} (10.3,11.7);(-10.3,11.7)
\ar@{-} (6.5,0);(-10.3,11.7)
\end{xy}}
}
&
$x_{13}$, $x_{14}$
\\
1
&
\raisebox{0pt}[50pt]
{
$
{\small
\begin{pmatrix}
0 & 1\\
-1 & 0\\
\hline
1 & -1\\
-1 & 0\\
1 & 0 \\
0 & -1 \\
0 & 1\\
\end{pmatrix}
}
$
}
&
{\small
\begin{xy}
(0,-3)*{a},
(10,-3)*{1},
(20,-3)*{2},
(30,-3)*{e},
(5,11)*{b},
(15,11)*{c},
(25,11)*{d},
(0,0)*{\bullet},
(10,0)*\cir<2pt>{},
(20,0)*\cir<2pt>{},
(30,0)*{\bullet},
(5,8)*{\bullet},
(15,8)*{\bullet},
(25,8)*{\bullet},
(19,-1); (1,-1) **\crv{(10,-10)};
\ar@{<-} (9,0);(1,0)
\ar@{<-} (19,0);(11,0)
\ar@{->} (29,0);(21,0)
\ar@{<-} (5.5,7);(9.5,1)
\ar@{->} (14.5,7);(10.5,1)
\ar@{<-} (24.5,7);(20.5,1)
\ar@{->} (2,-1.9);(1,-1)
\end{xy}
}
&
{
\small
\raisebox{-20pt}
{\begin{xy}
%edge
(-6,18)*{a},
(6,18)*{b},
(10.5,5.5)*{c},
(-10.5,5)*{e},
(0,-3)*{d},
%diagonal
(5,10)*{1},
(0,7)*{2},
(-12,12)*{1},
(12,12)*{3},
(0,22)*{2},
(8,-1)*{4},
(-8,-1)*{5},
%edges
\ar@{-} (-6.5,0);(6.5,0)
\ar@{-} (6.5,0);(10.3,11.7)
\ar@{-} (0,20);(10.3,11.7)
\ar@{-} (-6.5,0);(-10.3,11.7)
\ar@{-} (0,20);(-10.3,11.7)
%diagonal
\ar@{-} (0,20);(6.5,0)
\ar@{-} (6.5,0);(-10.3,11.7)
\end{xy}}
}
&
$x_{24}$, $x_{14}$
\\
2
&
\raisebox{0pt}[50pt]
{
$
{\small
\begin{pmatrix}
0 & -1\\
1 & 0\\
\hline
0 & 1\\
-1 & 0\\
1 & 0 \\
-1 & 1 \\
0 & -1\\
\end{pmatrix}
}
$
}
&
{\small
\begin{xy}
(0,-3)*{a},
(10,-3)*{1},
(20,-3)*{2},
(30,-3)*{e},
(5,11)*{b},
(15,11)*{c},
(25,11)*{d},
(0,0)*{\bullet},
(10,0)*\cir<2pt>{},
(20,0)*\cir<2pt>{},
(30,0)*{\bullet},
(5,8)*{\bullet},
(15,8)*{\bullet},
(25,8)*{\bullet},
(19,-1); (1,-1) **\crv{(10,-10)};
\ar@{->} (19,0);(11,0)
\ar@{<-} (29,0);(21,0)
\ar@{<-} (5.5,7);(9.5,1)
\ar@{->} (14.5,7);(10.5,1)
\ar@{->} (24.5,7);(20.5,1)
\ar@{->} (11,0.5);(24,7.5)
\ar@{->} (18,-1.9);(19,-1)
\end{xy}
}
&
{
\small
\raisebox{-20pt}
{\begin{xy}
%edge
(-6,18)*{a},
(6,18)*{b},
(10.5,5.5)*{c},
(-10.5,5)*{e},
(0,-3)*{d},
%diagonal
(5,10)*{1},
(-2,8)*{2},
(-12,12)*{1},
(12,12)*{3},
(0,22)*{2},
(8,-1)*{4},
(-8,-1)*{5},
%edges
\ar@{-} (-6.5,0);(6.5,0)
\ar@{-} (6.5,0);(10.3,11.7)
\ar@{-} (0,20);(10.3,11.7)
\ar@{-} (-6.5,0);(-10.3,11.7)
\ar@{-} (0,20);(-10.3,11.7)
%diagonal
\ar@{-} (0,20);(6.5,0)
\ar@{-} (0,20);(-6.5,0)
\end{xy}}
}
&
$x_{24}$, $x_{25}$
\\
3
&
\raisebox{0pt}[50pt]
{
$
{\small
\begin{pmatrix}
0 & 1\\
-1 & 0\\
\hline
0 & 1\\
1 & -1\\
-1 & 0 \\
1 & 0 \\
0 & -1\\
\end{pmatrix}
}
$
}
&
{\small
\begin{xy}
(0,-3)*{a},
(10,-3)*{1},
(20,-3)*{2},
(30,-3)*{e},
(5,11)*{b},
(15,11)*{c},
(25,11)*{d},
(0,0)*{\bullet},
(10,0)*\cir<2pt>{},
(20,0)*\cir<2pt>{},
(30,0)*{\bullet},
(5,8)*{\bullet},
(15,8)*{\bullet},
(25,8)*{\bullet},
(19,-1); (1,-1) **\crv{(10,-10)};
\ar@{<-} (19,0);(11,0)
\ar@{<-} (29,0);(21,0)
\ar@{->} (5.5,7);(9.5,1)
\ar@{<-} (14.5,7);(10.5,1)
\ar@{<-} (11,0.5);(24,7.5)
\ar@{->} (19,0.5);(6,7.5)
\ar@{->} (18,-1.9);(19,-1)
\end{xy}
}
&
{
\small
\raisebox{-20pt}
{\begin{xy}
%edge
(-6,18)*{a},
(6,18)*{b},
(10.5,5.5)*{c},
(-10.5,5)*{e},
(0,-3)*{d},
%diagonal
(3,4)*{1},
(-2,8)*{2},
(-12,12)*{1},
(12,12)*{3},
(0,22)*{2},
(8,-1)*{4},
(-8,-1)*{5},
%edges
\ar@{-} (-6.5,0);(6.5,0)
\ar@{-} (6.5,0);(10.3,11.7)
\ar@{-} (0,20);(10.3,11.7)
\ar@{-} (-6.5,0);(-10.3,11.7)
\ar@{-} (0,20);(-10.3,11.7)
%diagonal
\ar@{-} (10.3,11.7);(-6.5,0)
\ar@{-} (0,20);(-6.5,0)
\end{xy}}
}
&
$x_{35}$, $x_{25}$
\\
4
&
\raisebox{0pt}[50pt]
{
$
{\small
\begin{pmatrix}
0 & -1\\
1 & 0\\
\hline
0 & -1\\
0 & 1\\
-1 & 0 \\
1 & 0 \\
-1 & 1\\
\end{pmatrix}
}
$
}
&
{\small
\begin{xy}
(0,-3)*{a},
(10,-3)*{1},
(20,-3)*{2},
(30,-3)*{e},
(5,11)*{b},
(15,11)*{c},
(25,11)*{d},
(0,0)*{\bullet},
(10,0)*\cir<2pt>{},
(20,0)*\cir<2pt>{},
(30,0)*{\bullet},
(5,8)*{\bullet},
(15,8)*{\bullet},
(25,8)*{\bullet},
(19,-1); (1,-1) **\crv{(10,-10)};
(29,-1); (11,-1) **\crv{(20,-10)};
\ar@{->} (19,0);(11,0)
\ar@{->} (29,0);(21,0)
\ar@{<-} (14.5,7);(10.5,1)
\ar@{<-} (11,0.5);(24,7.5)
\ar@{<-} (19,0.5);(6,7.5)
\ar@{->} (2,-1.9);(1,-1)
\ar@{->} (28,-1.9);(29,-1)
\end{xy}
}
&
{
\small
\raisebox{-20pt}
{\begin{xy}
%edge
(-6,18)*{a},
(6,18)*{b},
(10.5,5.5)*{c},
(-10.5,5)*{e},
(0,-3)*{d},
%diagonal
(3,4)*{1},
(0,9.5)*{2},
(-12,12)*{1},
(12,12)*{3},
(0,22)*{2},
(8,-1)*{4},
(-8,-1)*{5},
%edges
\ar@{-} (-6.5,0);(6.5,0)
\ar@{-} (6.5,0);(10.3,11.7)
\ar@{-} (0,20);(10.3,11.7)
\ar@{-} (-6.5,0);(-10.3,11.7)
\ar@{-} (0,20);(-10.3,11.7)
%diagonal
\ar@{-} (10.3,11.7);(-6.5,0)
\ar@{-} (10.3,11.7);(-10.3,11.7)
\end{xy}}
}
&
$x_{35}$, $x_{13}$
\\
5
&
\raisebox{0pt}[50pt]
{
$
{\small
\begin{pmatrix}
0 & 1\\
-1 & 0\\
\hline
0 & -1\\
0 & 1\\
1 & -1\\
-1 & 0 \\
1 & 0\\
\end{pmatrix}
}
$
}
&
{\small
\begin{xy}
(0,-3)*{a},
(10,-3)*{1},
(20,-3)*{2},
(30,-3)*{e},
(5,11)*{b},
(15,11)*{c},
(25,11)*{d},
(0,0)*{\bullet},
(10,0)*\cir<2pt>{},
(20,0)*\cir<2pt>{},
(30,0)*{\bullet},
(5,8)*{\bullet},
(15,8)*{\bullet},
(25,8)*{\bullet},
(19,-1); (1,-1) **\crv{(10,-10)};
(29,-1); (11,-1) **\crv{(20,-10)};
\ar@{<-} (19,0);(11,0)
\ar@{->} (14.5,7);(10.5,1)
\ar@{<-} (15.5,7);(19.5,1)
\ar@{->} (11,0.5);(24,7.5)
\ar@{<-} (19,0.5);(6,7.5)
\ar@{->} (2,-1.9);(1,-1)
\ar@{->} (12,-1.9);(11,-1)
\end{xy}
}
&
{
\small
\raisebox{-20pt}
{\begin{xy}
%edge
(-6,18)*{a},
(6,18)*{b},
(10.5,5.5)*{c},
(-10.5,5)*{e},
(0,-3)*{d},
%diagonal
(-3,4)*{1},
(0,9.5)*{2},
(-12,12)*{1},
(12,12)*{3},
(0,22)*{2},
(8,-1)*{4},
(-8,-1)*{5},
%edges
\ar@{-} (-6.5,0);(6.5,0)
\ar@{-} (6.5,0);(10.3,11.7)
\ar@{-} (0,20);(10.3,11.7)
\ar@{-} (-6.5,0);(-10.3,11.7)
\ar@{-} (0,20);(-10.3,11.7)
%diagonal
\ar@{-} (-10.3,11.7);(6.5,0)
\ar@{-} (10.3,11.7);(-10.3,11.7)
\end{xy}}
}
&
$x_{14}$, $x_{13}$
\end{tabular}
\medskip
\caption{Mutations of extended exchange matrices.}
\label{1tab:Grass1}
\end{table}

The mutations of these data are presented up to $t=5$ in Table \ref{1tab:Grass1}, where
we  observe the familiar pentagon periodicity of type $A_2$ as the periodicity of the alternative flips in a pentagon.
In fact, this is the reason why we call this periodicity so.

\begin{rem}
This is indeed a prototypical example of the \emph{surface realization\/} of cluster patterns/algebras
developed by Fock-Goncharov \cite{Fock05} and Fomin-Shapiro-Thurston \cite{Fomin08}.
Only a  limited class of cluster patterns/algebras are realized in this way.
In \cite{Fomin08} a complete classification of the class of cluster algebras
  that admit the surface realization is given.
\end{rem}

Let $\calA$ be the cluster algebra
associated with the cluster pattern $\bfSigma$ 
in the sense of \eqref{1eq:geometric1}.
In order to match $\bbC[{\rmGr}(2,5)]$ in \eqref{1eq:cr1},  we replace the ground ring $\bbZ$
 with $\bbC$.
In view of the result in Example \ref{1ex:A21}, 
it is given by
\begin{align}
\label{1eq:AC1}
\calA=\bbC[x_{1;0}, x_{2;0}, x_{1;1}, x_{2;2},x_{1;3}, x_a, \dots, x_e].
\end{align}
We regard it as a graded ring such that every generator is of degree 1.

We have the following conclusion.

\begin{thm}[\cite{Fomin02,Fomin03a}]
The graded ring
$\bbC[{\rmGr}(2,5)]$ is isomorphic to  the cluster algebra $\calA$
in \eqref{1eq:AC1}.
In particular, $\bbC[{\rmGr}(2,5)]\rule{0pt}{12pt}$ has a cluster algebra structure.
\end{thm}

\begin{proof}

Let  $[ij]$ be the chord (diagonal or side) 
whose end points are labeled by $i$ and $j$
in the underlying pentagon 
of $T_t$.
We give the following labeling of cluster/frozen variables of $\calA$:
\begin{itemize}
\item
Each frozen variable $x_{i}$ ($i=a,\dots,e$) is labeled as
 $x_{jk}$ if
the side labeled by $i$ is $[jk]$.
Explicitly,
\begin{align}
x_a=x_{12},\ 
x_b=x_{23},\ 
x_c=x_{34},\ 
x_d=x_{45},\ 
x_e=x_{15}.
\end{align}
\item
Each cluster variable $x_{i;t}$ ($i=1,2$) is labeled as
$x_{jk}$ if
the diagonal labeled by $i$ in $T_t$ is $[jk]$ in
Table \ref{1tab:Grass1}.
In particular,
\begin{align}
 x_{1;0}=x_{13},\
 x_{2;0}=x_{14},\ 
 x_{1;1}=x_{24},\ 
 x_{2;2}=x_{25},\ 
 x_{1;3}=x_{35}.
 \end{align}
\end{itemize}
The cluster algebra $\calA$ in \eqref{1eq:AC1} is now written as
\begin{align}
\calA=\bbC[\bfX],
\quad
\bfX=(x_{ij})_{1\leq i<j\leq 5}.
\end{align}
On the other hand,
the mutation from $t=0$ to $1$
in \eqref{1eq:pto2} is identified with the Pl\"ucker relation
\begin{align}
\label{1eq:pto3}
R_{1234}: \quad x_{13}  x_{24}=x_{14} x_{23} +x_{12}x_{34}.
\end{align}
Similarly, the mutations from $t=1$ to $5$
are identified with the Pl\"ucker relations $R_{1245}$, $R_{2345}$, $R_{1235}$, $R_{1345}$.
Therefore, we have a graded ring homomorphism
\begin{align}
\label{eq:Ap1}
\begin{matrix}
\varphi:
&\bbC[{\rmGr}(2,5)]
=\bbC[\bfP]/I_R
&\rightarrow
&\calA=
\bbC[\bfX]
\\
&p_{ij} &\mapsto & x_{ij}.
\end{matrix}
\end{align}
The map $\varphi$ is clearly surjective.
To show that $\varphi$ is injective, we rely on the following facts:
\begin{itemize}
\item
The  ideal $I_R$ is prime  as stated in Definition \ref{1defn:Grassmannian1}.
\item
In particular, the  ideal $I_R$ contains no monomial in $\bfP$.
(If there is a monomial in $I_R$, it should be at least quadratic, which contradict the primeness.)
\end{itemize}
Suppose that a polynomial $F(\bfP)$ in $\bfP=(p_{ij})$ is in $\mathrm{Ker}\, \varphi$.
Namely, $F(\bfX)=0$ in $\calA$.
By multiplying some monomial $m(\bfX)$ to $F(\bfX)$ and applying the mutation relations
(the Pl\"ucker relations for $\bfX$), one can eliminate all non-initial cluster variables
$x_{24}$, $x_{25}$, $x_{35}$ and obtain  a polynomial $G(\bfx_0)$ in the initial cluster variables $\bfx_0$. 
Namely,
we have $G(\bfx_0)=m(\bfX)F(\bfX)=0$ in $\calA$.
Then, due to the algebraic independence of  $\bfx_0$, 
 $G(\bfx_0)$ is the zero polynomial.
 Applying the same operation to $F(\bfP)$ in $\bbC[\bfP]$,  we have $m(\bfP)F(\bfP)\equiv 0$ mod $I_R$.
Since $m(\bfP)\not\equiv 0$, we have $F(\bfP)\equiv 0$ by the primeness of $I_R$.
\end{proof}

\begin{rem}
It was shown by Scott \cite{Scott03} that 
 the homogeneous coordinate ring  of a general Grassmannian $\mathrm{Gr}(k,m)$ also has the cluster algebra structure.
\end{rem}

\newpage

\section{Separation formulas}

In this section we present the \emph{separation formulas}
given in CA4.
The are particularly important to study the
structure of seeds in a cluster pattern systematically.

\subsection{Principal coefficients}

We introduce the notion of cluster patterns with principal coefficients.
They are a special class of cluster patterns of geometric type
and  play an important role in studying  cluster patterns.

\begin{defn}
\index{coefficient!principal}\index{cluster pattern!with principal coefficients}
\label{1defn:principal1}
We say that a cluster pattern $\bfSigma$
is \emph{with principal coefficients at $t_0\in \bbT_n$}
if the following conditions are satisfied:
\begin{itemize}
\item
The coefficient semifield of $\bfSigma$ is a tropical semifield
$\mathrm{Trop}(\bfy)$ with generators $\bfy=(y_1,\dots,y_n)$,
where $n$ is the rank of $\bfSigma$.
\item
The coefficient tuple $\bfy_{t_0}$ at $t_0$ coincides with $\bfy$.
\end{itemize}
The same remark in Remark \ref{1rem:free1} is applicable
for the notation.
\end{defn}

In other words, a cluster pattern with principal coefficients at $t_0$
is a cluster pattern of geometric type such that
the extended exchange matrix in \eqref{1eq:extended1} at $t_0$
is given by
 \begin{align}
\label{1eq:Cinit1}
\tilde{B}_{t_0}
=
\left(
\begin{matrix}
B_{t_0}\\
I
\end{matrix}
\right).
\end{align}
All properties of cluster patterns of geometric type
hold to cluster patterns with principal coefficients.

\begin{rem}
In contrast to a  cluster pattern of free coefficients in
Definition \ref{1defn:free1},
the notion of principal coefficients crucially depends on the choice
of the base vertex $t_0$.
To explain it in more detail,
let $\bfy_t$ and $\bfy'_t$ ($t\in \bbT_n$) are principal coefficients at
base vertices $t_0$ and $t'_0$, respectively.
Then, the correspondence
 $\varphi: y_{i;t}\mapsto y'_{i;t}$ \emph{cannot} be
extended to a semifield homomorphism from 
$\mathrm{Trop}(\bfy_{t_0})$ to $\mathrm{Trop}(\bfy_{t'_0})$,
in general.
\end{rem}

Let $\bfSigma$ be a  cluster pattern  with principal coefficients
at $t_0$.
One may regard the coefficients of $\bfSigma$
 as the
\emph{tropicalization
of  free coefficients\/} at $t_0$.
\index{tropicalization!of free coefficients}
Namely,
let $\bfSigma'$ be a cluster pattern  with free coefficients 
at $t_0$
such that it shares the common $B$-pattern
 with  $\bfSigma$.
Then,
as a special case of Proposition \ref{1prop:free1},
 $\bfy'_t$ and $\bfy_t$  are related by the
tropicalization homomorphism in \eqref{1eq:trophom1} by
\begin{align}
\label{1eq:trophom2}
\begin{matrix}
\pi_{\mathrm{trop}}: & \mathbb{Q}_{\mathrm{sf}}(\bfy)
&\rightarrow& \mathrm{Trop}(\bfy)\\
& {y}'_{i;t}
& \mapsto
&y_{i;t}.
\end{matrix}
\end{align}

\subsection{$C$- and $G$-matrices, and $F$-polynomials}

For a given cluster pattern $\bfSigma$ with principal coefficients
at $t_0$, one can 
 define important quantities called
  $C$-matrices, $G$-matrices, and $F$-polynomials.
  They are the \emph{building blocks} of seeds of a cluster pattern with coefficients in any semifield $\bbP$,
  and together with the separation formula  in Theorem \ref{1thm:sep1},
  they clarify the structure of seeds as well as the relation between cluster variables and coefficients.
 
\medskip
\par
 \noindent
 {\bf (a). $C$-matrices.} Let us start with $C$-matrices.
 They are  the matrices defined
 in \eqref{eq:cexp1} 
 specialized
 for a cluster pattern $\bfSigma$ with principal coefficients
at $t_0$.

\begin{defn}[$C$-matrix/$c$-vector]
For a given cluster pattern $\bfSigma$ with principal coefficients
at $t_0$,
the \emph{$C$-matrix\/} $C_t=(c_{ij;t})$ \index{matrix!$C$-}\index{$C$-matrix}
of $\bfy_t$
is the $n \times n$ integer matrix defined by
\begin{align}
y_{i;t}=\prod_{j=1}^n y_j^{c_{ji;t}}.
\end{align}
Equivalently,
they are  defined as the lower half of the
$2n \times n$ extended exchange matrix $\tilde{B}_t$
of  \eqref{1eq:extended1} with the initial condition 
\eqref{1eq:Cinit1}.
The $i$th column vector $\bfc_{i;t}=(c_{ji})_{j=1}^n$ of $C_t$ is called
the \emph{$c$-vector\/} of $y_{i;t}$.  \index{vector!$c$-}\index{$c$-vector}
We  call 
the collection of the $C$-matrices $\bfC^{t_0}=\{ C_t\}_{t\in \bbT_n}$ 
the \emph{$C$-pattern\/} of $\bfSigma$.  \index{pattern!$C$-}\index{$C$-pattern}
\end{defn}

The $C$-matrices are uniquely determined by the underlying
 $B$-pattern $\bfB$ of $\bfSigma$ and $t_0$, thus,
  eventually by $t_0$ and $B_{t_0}$ only.

\begin{prop}
\label{1prop:Crel1}
The $C$-pattern  $\bfC^{t_0}$ of $\bfSigma$ is uniquely determined by the following
initial condition and the mutation formula:
\begin{align}
\label{1eq:Cmat1}
C_{t_0}&=I,
\\
 \label{1eq:cmut2}
c_{ij;t'}&=
\begin{cases}
-c_{ik;t} & j = k,\\
c_{ij;t} 
+ c_{ik;t}[b_{kj;t}]_+
+ [-c_{ik;t}]_+b_{kj;t} 
& j \neq k,
\end{cases}
 \end{align}
where $t$ and $t'$ are $k$-adjacent.
\end{prop}
\begin{proof}
The initial condition \eqref{1eq:Cmat1} holds by \eqref{1eq:Cinit1}.
The mutation \eqref{1eq:cmut2} was established in \eqref{1eq:cmut1}.
\end{proof}

For the second case of \eqref{1eq:cmut2},
we have an alternative expression,
which is parallel to the $\varepsilon$-expressions in Proposition \ref{1prop:epsilon1}.

\begin{prop}
\index{$\varepsilon$-expression!for $C$-matrices}
The following expression does not depend on the choice of
$\varepsilon \in \{1, -1\}$:
\begin{align}
\label{1eq:Ceps1}
c_{ij;t} 
+ c_{ik;t}[\varepsilon b_{kj;t}]_+
+ [-\varepsilon c_{ik;t}]_+b_{kj;t} .
\end{align}
\end{prop}
\begin{proof}
This follows from \eqref{1eq:pos2}.
Alternatively, 
one can derive each expression from \eqref{1eq:ymut2}
in the same way as \eqref{1eq:cmut1}.
\end{proof}

 \medskip
 \noindent
 {\bf (b). $G$-matrices.}
 Next, we define $G$-matrices.
By Theorem \ref{1thm:Laurent2}, we have
 \begin{align}
 \label{1eq:Laurent2}
x_{i;t} \in \bbZ[\bfx^{\pm1}, \bfy],
\end{align}
where we set $\bfx_{t_0}=\bfx$.
We introduce a certain degree vector in $\bbZ^n$ for each cluster variables $x_{i;t}$.
 
 \begin{defn}[Principal $\bbZ^n$-grading]
 \label{1defn:principalg1}
  For each monomial
  in $\bbZ[\bfx^{\pm1}, \bfy]$,
  we define its  $\bbZ^n$-grading by
  \begin{align}
    \label{1eq:deg1}
  \deg(x_{i})=\bfe_i,
  \quad
    \deg(y_{i})=-\bfb_{i;t_0},
  \end{align}
  where $\bfe_i$ is the $i$th unit vector and $\bfb_{i;t_0}$ is the $i$th column vector of $B_{t_0}$.
  We call it the \emph{principal\/ \hbox{$\bbZ^n$-grading}}.
  \index{principal $\bbZ^n$-grading}
  \end{defn}
  
  A seemingly artificial degree of $y_i$ is designed to ensures the following property:
  \begin{align}
  \label{1eq:deg2}
      \deg(\haty_{i})=
      \deg\Biggl(y_i
      \prod_{j=1}^n x_{j}^{b_{ji;t_0}}
      \Biggr)
      =-\bfb_{i;t_0}+\bfb_{i;t_0}=
      \bfzero.
  \end{align}
 
\begin{lem}
\label{1lem:hom1}
Every $x_{i;t} \in \bbZ[\bfx^{\pm1}, \bfy]$ is  homogenous with respect to
the principal $\bbZ^n$-grading.
\end{lem}
\begin{proof}
We first note that, thanks to Proposition \ref{1prop:yhat1},
any $\haty$-variable $ \haty_{i;t}$ is written as a rational function
 of the initial $\haty$-variables $\hat\bfy$.
Then, we can effectively treat it as $    \deg(\haty_{i;t})=\bfzero$ 
in the following calculation, even though $\haty_{i;t}$ does not belong to $ \bbZ[\bfx^{\pm1}, \bfy]$.
We prove the claim by the induction on $d=d(t_0,t)$.
The claim is trivial for the initial cluster variables $x_{i}$.
Suppose that claim is true up to $d=d(t_0,t)$.  
Let $t'\in \bbT_n$ be $k$-adjacent to $t$.
We look at the mutation formula in the form \eqref{1eq:xmut1},
\begin{align}
\label{1eq:xmut7}
x_{k;t'}
&=
x_{k;t}^{-1}\Biggl(\, \prod_{j=1}^n x_{j;t}^{[-b_{jk;t}]_+}
\Biggr)
\frac{ 1+\hat{y}_{k;t}}{ 1\oplus y_{k;t}}.
\end{align}
Then, the binomial $ 1+\hat{y}_{k;t}$ is homogeneous with degree $\bfzero$.
Also, $1\oplus y_{k;t}$ is actually a monomial in $\bfy$ because it belongs to $\mathrm{Trop}(\bfy)$.
Thus, the right hand side of \eqref{1eq:xmut7} is  homogeneous thanks to the induction assumption.
\end{proof}

Based on Lemma \ref{1lem:hom1}, we define $G$-matrices as follows.

\begin{defn}[$G$-matrix/$g$-vector]
For a given cluster pattern $\bfSigma$ with principal coefficients
at $t_0$,
we define a $\bbZ^n$-vector
\begin{align}
\bfg_{i;t}=\deg (x_{i;t})\in \bbZ^n,
\end{align}
where $\deg$ is the principal $\bbZ^n$-grading in Definition \ref{1defn:principalg1}.
We call it the \emph{\hbox{$g$-vector}} of $x_{i;t}$.  \index{vector!$g$-}\index{$g$-vector}
For each $t\in \bbT_n$,
we  define the \emph{$G$-matrix\/} $G_t=(g_{ij;t})$  of $\bfx_t$ such that  \index{matrix!$G$-}\index{$G$-matrix}
its $i$th column vector is $\bfg_{i;t}=(g_{ji;t})_{j=1}^n$.
We  call 
the collection of the $G$-matrices $\bfG^{t_0}=\{ G_t\}_{t\in \bbT_n}$ 
the \emph{$G$-pattern\/} of $\bfSigma$.  \index{pattern!$G$-}\index{$G$-pattern}
\end{defn}

Again,
the $G$-pattern $\bfG^{t_0}$  is uniquely determined by 
 the underlying $B$-pattern $\bfB$ and $t_0$ (together with the $C$-pattern $\bfC^{t_0}$ determined from them).

\begin{prop}
\label{1prop:Grel1}
The  $G$-pattern $\bfG^{t_0}$ of $\bfSigma$ is uniquely determined by the following
initial condition and the mutation formula:
\begin{align}
\label{1eq:Gmat1}
G_{t_0}&=I,
\\
 \label{1eq:gmut1}
 g_{ij;t'}&=
 \begin{cases}
 \displaystyle
 -g_{ik;t}
 + \sum_{\ell=1}^n g_{i\ell;t} [-b_{\ell k;t}]_+
 -  \sum_{\ell=1}^nb_{i\ell;t_0}  [-c_{\ell k;t}]_+ 
 &
 j= k,
 \\
 g_{ij;t}  &
 j \neq k,
 \end{cases}
 \end{align}
 where $t$ and $t'$ are $k$-adjacent.
\end{prop}
\begin{proof}
The initial condition \eqref{1eq:Gmat1} follows from the fact $\deg x_i = \bfe_i$.
The mutation \eqref{1eq:gmut1} follows from \eqref{1eq:xmut7} and
the formula
\begin{align}
\deg\biggl( \frac{1}{1\oplus y_{k;t}}\biggr)= 
\deg \Biggl(\, \prod_{j=1}^n y_j^{[-c_{jk;t}]_+}\Biggr)
=
- \sum_{j=1}^n  [-c_{jk;t}]_+ \bfb_{j;t_0}
,
\end{align}
where we used \eqref{1eq:1+y2} in the first equality.
\end{proof}

We  have the following duality relation between
 $C$- and $G$-matrices.
\begin{prop} The following equality holds:
\begin{align}
\label{1eq:GBBC1}
G_t B_t = B_{t_0} C_{t}.
\end{align}
\end{prop}
\begin{proof}
This follows from the formula
 \begin{align}
 \begin{split}
  \label{1eq:deg3}
 \deg(\haty_{i;t})&=
      \deg\Biggl(y_{i;t}
      \prod_{j=1}^n x_{j;t}^{b_{ji;t}}
      \Biggr)
      =
            \deg\Biggl(\,\prod_{j=1}^n y_j^{c_{ji;t}}
      \prod_{j=1}^n x_{j;t}^{b_{ji;t}}
      \Biggr)\\
      &=-B_{t_0}\bfc_{i;t}+G_t\bfb_{i;t}.
      \end{split}
  \end{align}
and the fact $\deg(\haty_{i;t})=\bfzero$.
\end{proof}

For the first case of \eqref{1eq:gmut1},
 we have an alternative expression,
which is parallel to the $\varepsilon$-expressions in Proposition \ref{1prop:epsilon1}.

\begin{prop}
\index{$\varepsilon$-expression!for $G$-matrices}
The following expression does not depend on the choice of
$\varepsilon \in \{1, -1\}$:
\begin{align}
\label{1eq:Geps1}
 -g_{ik;t}
 + \sum_{\ell=1}^n g_{i\ell;t} [-\varepsilon b_{\ell k;t}]_+
 -  \sum_{\ell=1}^nb_{i\ell;t_0}  [-\varepsilon c_{\ell k;t}]_+ .
\end{align}
\end{prop}
\begin{proof}
We take the difference of two expressions.
Then, it vanishes thanks to \eqref{1eq:pos1}
and \eqref{1eq:GBBC1}.
\end{proof}

 \medskip
 \noindent
 {\bf (c). $F$-polynomials.}
 
 The following definition relies on the fact 
\eqref {1eq:Laurent2}.
 \begin{defn}[$F$-polynomial]
 For a given cluster pattern $\bfSigma$ with principal coefficients
at $t_0$,
we define  polynomials $F_{i;t}(\bfy)\in \bbZ[\bfy]$
in formal variables $\bfy=(y_1,\dots,y_n)$
by specializing Laurent polynomials $x_{i;t}(\bfx,\bfy) \in \bbZ[\bfx^{\pm1}, \bfy]$
with $x_{1}=\dots = x_n=1$.
We call  $F_{i;t}(\bfy)$  the \emph{$F$-polynomial\/} of $x_{i;t}$.  \index{$F$-polynomial}
Let $\bfF_t=(F_{1;t}(\bfy),\dots, F_{n;t}(\bfy))$.
We  call 
the collection of the $F$-polynomials $\bfF^{t_0}=\{ \bfF_t\}_{t\in \bbT_n}$ 
the \emph{$F$-pattern\/} of $\bfSigma$.  \index{pattern!$F$-}\index{$F$-pattern}
\end{defn}

Again, the $F$-pattern $\bfF^{t_0}$  is
 uniquely determined by
  the underlying $B$-pattern $\bfB$ and $t_0$ (together with the $C$-pattern $\bfC^{t_0}$ determined from them).

\begin{prop}
\label{1prop:Frel1}
The  $F$-pattern $\bfF^{t_0}$ of $\bfSigma$ is uniquely determined by the following
initial condition and the mutation formula:
\begin{align}
\label{1eq:Finit1}
F_{i;t_0}(\bfy)&=1,
\\
 \label{1eq:Fmut1}
 F_{i;t'}(\bfy)&=
 \begin{cases}
\frac
 {
  \displaystyle
 M_{k;t}(\bfy)
 }
{ \displaystyle
 F_{k;t}(\bfy)
 }
   &
 i= k,
 \\
 F_{i;t}(\bfy)  &
 i \neq k,
 \end{cases}
 \end{align}
 where $t$ and $t'$ are $k$-adjacent,
and
\begin{align}
\label{1eq:M1}
 M_{k;t}(\bfy)
 &=
    \prod_{j=1}^{n}
  y_j^{[c_{jk;t}]_+}
    \prod_{j=1}^{n}
  F_{j;t}(\bfy)^{[b_{jk;t}]_+}
+
    \prod_{j=1}^{n}
  y_j^{[-c_{jk;t}]_+}
    \prod_{j=1}^{n}
  F_{j;t}(\bfy)^{[-b_{jk;t}]_+}
  \\
  \label{1eq:M2}
   &=
\Biggl(\,
    \prod_{j=1}^{n}
  y_j^{[-c_{jk;t}]_+}
    \prod_{j=1}^{n}
  F_{j;t}(\bfy)^{[-b_{jk;t}]_+}
  \Biggr)
  \Biggl(
  1
 + 
 \prod_{j=1}^{n}
  y_j^{c_{jk;t}}
    \prod_{j=1}^{n}
  F_{j;t}(\bfy)^{b_{jk;t}}
  \Biggr).
\end{align}
\end{prop}
\begin{proof}
This is obtained by specializing
the mutation  of $\bfx_t$ in the form \eqref{1eq:xmut5}
with $x_{1}=\dots = x_n=1$.
\end{proof}

We may regard $F$-polynomials also as elements in $\mathbb{Q}_{\mathrm{sf}}(\bfy)$
because the mutation in \eqref{1eq:Fmut1} is a subtraction-free operation.
Then, we may apply 
the tropicalization homomorphism 
$\pi_{\mathrm{trop}}:\mathbb{Q}_{\mathrm{sf}}(\bfy)
\rightarrow \mathrm{Trop}(\bfy)$
in \eqref{1eq:trophom2}.

\begin{prop} The following fact holds:
\begin{align}
\label{1eq:F1}
\pi_{\mathrm{trop}}(F_{i;t}(\bfy))=1.
\end{align}
\end{prop}
\begin{proof}
We prove it by the induction on the distance $d=d(t_0,t)$.
For $t=t_0$, the claim holds by \eqref{1eq:Finit1}.
Assume that the claim holds for 
$d=d(t_0,t)$.
Let $t'$ be $k$-adjacent to  $t$.
Note that, for each $j$, either $[c_{jk;t}]_+$ or $[-c_{jk;t}]_+$ is zero.
Then, $\pi_{\mathrm{trop}}( M_{k;t}(\bfy))=1$ for
$ M_{k;t}(\bfy)$ in \eqref{1eq:M1}.
So, the claim also holds for $t'$
by \eqref{1eq:Fmut1}
\end{proof}

\begin{rem}
The fact \eqref{1eq:F1} does not imply that $F_{i;t}(\bfy)$ has a nonzero constant term
as a polynomial in $\bfy$. For example,
$\pi_{\mathrm{trop}}(y_1+y_2)=1$.

\end{rem}

In summary,
even if the patterns $\bfC^{t_0}$, $\bfG^{t_0}$, $\bfF^{t_0}$
are originally extracted from a cluster pattern $\bfSigma$ with principal coefficients
at $t_0$, it turned out that they are uniquely and directly defined 
from the underlying $B$-pattern $\bfB$ and the base vertex $t_0$.
Therefore, one can associate these patterns $\bfC^{t_0}$, $\bfG^{t_0}$, $\bfF^{t_0}$
to any cluster pattern $\bfSigma$ with
coefficients in 
 \emph{any  semifield $\bbP$}.

\subsection{Separation formulas}

Let $\bfSigma$ be any cluster pattern with coefficients in any semifield $\bbP$.
Let $t_0$ be a given initial vertex $t_0$,
and let $\bfC^{t_0}$, $\bfG^{t_0}$, $\bfF^{t_0}$ be
the associated $C$-, $G$-, $F$-patterns, respectively.
As already mentioned,
each $F$-polynomial $F_{i;t}(\bfy)$ belongs to $\bbQ_{\mathrm{sf}}(\bfy)$.
Following Definition \ref{1defn:sp1},
let $F_{i;t}{\vert}_{\bbP}(\bfy)$ be the specialization of
$F_{i;t}(\bfy)$ in $\bbP$ at 
the initial coefficient $\bfy_{t_0}=\bfy$ of $\bfSigma$,
where the abuse of the symbol $\bfy$ does not cause a serious problem.

The cluster variables and coefficients
of any cluster pattern are expressed
by the initial cluster variables and coefficients
together with $C$- and $G$-matrices, and $F$-polynomials.
This is one of the most fundamental and useful
properties of cluster patterns.

\begin{thm}[{Separation Formulas \cite{Fomin07}}]
\index{separation formula}
\label{1thm:sep1}
Let $\bfSigma$ be any cluster pattern
with coefficients in any semifield $\bbP$ and with a given initial vertex $t_0$.
Let 
\begin{align}
\bfx_{t_0}=\bfx,
\quad
\bfy_{t_0}=\bfy,
\quad
\hat\bfy_{t_0}=\hat\bfy
\end{align}
be the initial cluster variables,
coefficients, and $\haty$-variables.
Then,
the following formulas hold.
\begin{align}
\label{1eq:sep1}
x_{i;t}&=
\Biggl(\,
\prod_{j=1}^n
x_j^{g_{ji;t}}
\Biggr)
\frac{F_{i;t}(\hat{\bfy})}{F_{i;t}\vert_{\bbP}(\bfy)},
\\
\label{1eq:sep2}
y_{i;t}&=
\Biggl(\,
\prod_{j=1}^n
y_j^{c_{ji;t}}
\Biggr)
\prod_{j=1}^n
F_{j;t}\vert_{\bbP}(\bfy)^{b_{ji;t}},
\\
\label{1eq:sep3}
\haty_{i;t}&=
\Biggl(\,
\prod_{j=1}^n
\haty_j^{c_{ji;t}}
\Biggr)
\prod_{j=1}^n
F_{j;t}(\hat\bfy)^{b_{ji;t}}.
\end{align}
\end{thm}
\begin{proof}
The formulas hold for $t=t_0$
due to the initial condition of $C$- and $G$-matrices and $F$-polynomials
in Propositions \ref{1prop:Crel1}, \ref{1prop:Grel1}, and \ref{1prop:Frel1}.
Therefore, it is enough to show that
their right hand sides 
mutate in the same way as $x_{i;t}$, $y_{i;t}$, $\haty_{i;t}$
based on the mutations of $C_t$, $G_t$, and $F_{i;t}$.

We first treat  the formula \eqref{1eq:sep2}.
Let $y_{i;t}$ temporarily denote the one in  the right hand side of \eqref{1eq:sep2}.
Let $t'\in \bbT_n$ be  $k$-adjacent to $t$.
For $i= k$, 
\begin{align}
y_{k;t'}
&=
\Biggl(\,
\prod_{j=1}^n
y_j^{-c_{jk;t}}
\Biggr)
\prod_{j=1}^n
F_{j;t}\vert_{\bbP}(\bfy)^{
-b_{jk;t}
}
=
y_{k;t}^{-1},
\end{align}
where in the first equality we used the fact $b_{kk;t}=0$.
(Below we do not repeat this remark.)
For $i\neq k$,
\begin{align}
\begin{split}
y_{i;t'}
&=
\Biggl(\,
\prod_{j=1}^n
y_j^{c_{ji;t}
+
c_{jk;t}[b_{ki;t}]_+
+
[-c_{jk;t}]_+b_{ki;t}
}
\Biggr)
\\
&
\qquad
\times
\Biggl(\,
\prod_{j\neq k}
F_{j;t}\vert_{\bbP}(\bfy)^{
b_{ji;t}
+
b_{jk;t} [b_{ki;t}]_+
+
[-b_{jk;t}]_+b_{ki;t}
}
\Biggr)
\\
&
\qquad
\times
\Biggl\{
 F_{k;t}\vert_{\bbP}(\bfy)^{-1}
\Biggl(\,
    \prod_{j=1}^{n}
  y_j^{[-c_{jk;t}]_+}
    \prod_{j=1}^{n}
  F_{j;t}\vert_{\bbP}(\bfy)^{[-b_{jk;t}]_+}
\Biggr)
  \\
  &
\qquad  \times
\Biggl(
  1
 \oplus
 \prod_{j=1}^{n}
  y_j^{c_{jk;t}}
    \prod_{j=1}^{n}
  F_{j;t}\vert_{\bbP}(\bfy)^{b_{jk;t}}
\Biggr)
  \Biggr\}^{-b_{ki;t}}
  \\
  &=
  y_{i;t}y_{k;t}^{[b_{ki;t}]_+}(1\oplus y_{k;t})^{-b_{ki;t}}.
  \end{split}
\end{align}
The formula  \eqref{1eq:sep3} is proved 
 in the exactly same manner as above.

Now we treat \eqref{1eq:sep1}.
Let $x_{i;t}$ temporarily denote the one in  the right hand side of \eqref{1eq:sep1}.
For $i\neq k$, $x_{i;t'}=x_{i;t}$ holds.
For $x_{k;t'}$, we have
\begin{align*}
\qquad x_{k;t'}
&=
\Biggl(\,
\prod_{j=1}^n
x_j^{ -g_{jk;t}
 + \sum_{\ell=1}^n g_{j\ell;t} [-b_{\ell k;t}]_+
 -  \sum_{\ell=1}^nb_{j\ell;t_0}  [-c_{\ell k;t}]_+ }
\Biggr)
\\
&\qquad
\times
\frac{
\displaystyle
F_{k;t}(\hat{\bfy})^{-1}
    \prod_{j=1}^{n}
  \hat{y}_j^{[-c_{jk;t}]_+}
    \prod_{j=1}^{n}
  F_{j;t}(\hat{\bfy})^{[-b_{jk;t}]_+}
  }
  {
  \displaystyle
  F_{k;t}\vert_{\bbP}(\bfy)^{-1}
    \prod_{j=1}^{n}
 y_j^{[-c_{jk;t}]_+}
    \prod_{j=1}^{n}
  F_{j;t}\vert_{\bbP}(\bfy)^{[-b_{jk;t}]_+}
  }
  \\
  &\qquad
\times
  \frac{
  \displaystyle
  1
 + 
 \prod_{j=1}^{n}
  \hat{y}_j^{c_{jk;t}}
    \prod_{j=1}^{n}
  F_{j;t}(\hat{\bfy})^{b_{jk;t}}
  }
  {
  \displaystyle
  1
 \oplus
 \prod_{j=1}^{n}
 y_j^{c_{jk;t}}
    \prod_{j=1}^{n}
  F_{j;t}\vert_{\bbP}(\bfy)^{b_{jk;t}}
  }
  \\
  &=
  x_{k;t}^{-1}
\Biggl(\,
  \prod_{j=1}^n x_{j;t}^{[-b_{jk;t}]_+}
\Biggr)
  \frac{1+\hat{y}_{k;t}}
  {1\oplus y_{k;t}}
,
\end{align*}
where we used
\eqref{1eq:sep2} and \eqref{1eq:sep3}
in the last equality.
\end{proof}

\begin{rem}
Originally in \cite{Fomin07},
the formula \eqref{1eq:sep1} was referred to as the separation formula,
where the additions $\oplus/+$ in $\bbP/\calF$ are \emph{separated\/}
in the denominator and the numerator.
Here we also included the formulas \eqref{1eq:sep2} and \eqref{1eq:sep3}
in the family,
because they together exhibit a  duality relation between cluster variables
and coefficients/$\haty$-variables.
\end{rem}

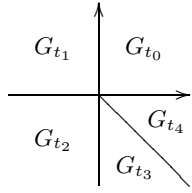
\begin{figure}
\centering
\leavevmode
\begin{xy}
0;/r1.2mm/:,
(5, 5)*{\text{\small $G_{t_{0}}$}},
(-5, 5)*{\text{\small $G_{t_{1}}$}},
(-5, -5)*{\text{\small $G_{t_{2}}$}},
(4, -8)*{\text{\small $G_{t_{3}}$}},
(7.5, -3)*{\text{\small $G_{t_{4}}$}},
(0,0)="A",
\ar "A"+(0,0); "A"+(10,0)
\ar "A"+(0,0); "A"+(0,10)
\ar@{-} "A"+(0,0); "A"+(-10,0)
\ar@{-} "A"+(0,0); "A"+(0,-10)
\ar@{-} "A"+(0,0); "A"+(10,-10)
\end{xy}
\caption{The $G$-fan of  a cluster pattern of type $A_2$.}
\label{1fig:Gfan1}
\end{figure}

\begin{ex}[Type $A_2$]
\label{1ex:typeA22}
Let us consider the cluster pattern of type $A_2$
in Example \ref{1ex:A21}.
By comparing the result therein with the separation formulas 
\eqref{1eq:sep1} and \eqref{1eq:sep2},
one can easily read off $C$- and $G$-matrices and $F$-polynomials
 as follows:
\begin{alignat}{5}
C_{0}&=
\begin{pmatrix}
1 & 0\\
0 & 1
\end{pmatrix},
\
&
G_{0}&=
\begin{pmatrix}
1 & 0\\
0 & 1
\end{pmatrix},
\
&
&
\begin{cases}
F_{1;0}(\bfy)=1,\\
F_{2;0}(\bfy)=1,\\
\end{cases}
\\
C_{1}&=
\begin{pmatrix}
-1 & 0\\
0 & 1
\end{pmatrix},
\
&
G_{1}&=
\begin{pmatrix}
-1 & 0\\
0 & 1
\end{pmatrix},
\
&&
\begin{cases}
F_{1;1}(\bfy)=1+y_1,\\
F_{2;1}(\bfy)=1,\\
\end{cases}
\\
C_{2}&=
\begin{pmatrix}
-1 & 0\\
0 & -1
\end{pmatrix},
\quad
&
G_{2}&=
\begin{pmatrix}
-1 & 0\\
0 & -1
\end{pmatrix},
\quad
&&
\begin{cases}
F_{1;2}(\bfy)=1+y_1,\\
F_{2;2}(\bfy)=1+y_2+y_1y_2,
\end{cases}
\\
C_{3}&=
\begin{pmatrix}
1 & -1\\
0 & -1
\end{pmatrix},
\
&
G_{3}&=
\begin{pmatrix}
1 & 0\\
-1 & -1
\end{pmatrix},
\
&&
\begin{cases}
F_{1;3}(\bfy)=1+y_2,\\
F_{2;3}(\bfy)=1+y_2+y_1y_2,\\
\end{cases}
\\
C_{4}&=
\begin{pmatrix}
0 & 1\\
-1 & 1
\end{pmatrix},
\
&
G_{4}&=
\begin{pmatrix}
1 & 1\\
-1 & 0
\end{pmatrix},
\
&&
\begin{cases}
F_{1;4}(\bfy)=1+y_2,\\
F_{2;4}(\bfy)=1,\\
\end{cases}
\\
C_{5}&=
\begin{pmatrix}
0 & 1\\
1 & 0
\end{pmatrix},
\
&
G_{5}&=
\begin{pmatrix}
0 & 1\\
1 & 0
\end{pmatrix},
\
&&
\begin{cases}
F_{1;5}(\bfy)=1,\\
F_{2;5}(\bfy)=1.\\
\end{cases}
\end{alignat}
We have the following observations.
\par
(a). \emph{Duality.}
The following relation between $C$- and $G$-matrices holds:
\begin{align}
\label{1eq:dual1}
G_t^T C_t=I,
\end{align}
where $G_t^T$ is the transpose of $G_t$.
\par
(b). \emph{$G$-fan.}
\index{fan!$G$-}\index{$G$-fan}
For each $G$-matrix $G_t$, let
\begin{align}
\sigma(G_t):=
\bbR_{\geq 0} \bfg_{1;t} 
+ \bbR_{\geq 0} \bfg_{2;t}
\end{align}
be the cone 
in $\bbR^2$ spanned by its $g$-vectors,
which we call a \emph{$G$-cone}.
$G$-cones intersect only in their boundaries, thereby forming a \emph{fan\/}
(the  \emph{$G$-fan\/}, or  the \emph{$g$-vector fan}  of $\bfSigma$).
See Figure \ref{1fig:Gfan1}.
\end{ex} 

In view of the definitions of
$C$- and $G$-matrices and $F$-polynomials,
we define the action of a permutation $\sigma\in S_n$
on them as follows:
\begin{gather}
\sigma C_t=C',\quad c'_{ij}:=c_{i\sigma^{-1}(j);t},
\quad
\sigma G_t=G',\quad g'_{ij}:=g_{i\sigma^{-1}(j);t},
\\
\sigma \bfF_t=\bfF',\quad
F'_{i}(\bfy):=F_{\sigma^{-1}(i);t}(\bfy).
\end{gather}

The following is a consequence of
the separation formulas.

\begin{prop}
\label{1eq:period2}
Let $\bfSigma$ be a cluster pattern with free coefficients at $t_0$.
Let $\bfSigma'$ be a cluster pattern with  principal coefficients at $t_0$
sharing the common $B$-pattern with $\bfSigma$.
Then, for any $t,t'\in \bbT_n$ and a permutation $\sigma\in S_n$,
the following fact holds:
\begin{align}
\label{1eq:period3}
\Sigma_t=
\sigma \Sigma_{t'}
\quad
\Longleftrightarrow
\quad
\Sigma'_t=
\sigma \Sigma'_{t'}.
\end{align}
In other words, the periodicities of $\bfSigma$ and
 $\bfSigma'$ coincide.
\end{prop}
\begin{proof}
The implication $\Longrightarrow$ was given in Proposition \ref{1prop:period1}.
Let us show the  implication $\Longleftarrow$.
Assume the equality $\Sigma'_t=
\sigma \Sigma'_{t'}$. 
Then, it implies the same periodicity for
$C$- and $G$-matrices and $F$-polynomials,
\begin{align}
C_t=\sigma C_{t'},
\quad
G_t=\sigma G_{t'},
\quad
\bfF_t=\sigma \bfF_{t'}.
\end{align}
Also,  we have $B_t=\sigma B_{t'}$.
Then, applying the separation formulas \eqref{1eq:sep1} and \eqref{1eq:sep2}
to $\bfSigma$,
we have
\begin{align}
\bfx_t= \sigma \bfx_{t'},
\quad
\bfy_t= \sigma \bfy_{t'}.
\end{align}
\end{proof}

In summary, through the separation formulas,
a cluster pattern of principal coefficients controls
a cluster pattern of any coefficients
(including free coefficients) sharing the common $B$-pattern.
This is the reason why it is called \emph{principal\/} in CA4.

\subsection{Further results}
\label{1subsec:further1}
In the rest of the section, we present without proof some advanced results
on $C$- and $G$-matrices and $F$-polynomials,
which are fundamentally important in cluster algebra theory.

\begin{defn}
We say that a vector $\bfv\in \bbZ^n$ is \emph{positive (resp.~negative)\/}
\index{vector!positive/negative}
if it is a nonzero vector, and all nonzero components are positive (resp.~negative).
We say that a square matrix is \emph{column sign-coherent\/} if each column vector
\index{matrix!column/row sign-coherent}
is either positive or negative.
\end{defn}

The following  theorems were conjectured by Fomin-Zelevinsky \cite{Fomin02,Fomin07}
 and proved by Gross-Hacking-Keel-Kontsevich \cite{Gross14} by the scattering diagram method.

\begin{thm}[{Sign-coherence of $C$-matrices \cite{Gross14}}]
\index{sign-coherence!of $C$-matrix}
\label{1thm:sign1}
Every $c$-vector $\bfc_{i;t}$ is either positive or negative.
In other words, every $C$-matrix $C_t$ is  column sign-coherent.
\end{thm}

\begin{thm}[Unit constant property {\cite{Gross14}}]
\label{1thm:const1}
Every $F$-polynomial $F_{i;t}(\bfy)$ has constant term $1$.
\end{thm}

\begin{thm}[Laurent positivity {\cite{Gross14}}]
\index{Laurent positivity}
\label{1thm:positive1}
Every $F$-polynomial $F_{i;t}(\bfy)$ has no negative coefficients.
\end{thm}

It was shown in  \cite{Fomin07} that
Theorems \ref{1thm:sign1} and \ref{1thm:const1}
are equivalent to each other.

\begin{proof}[Proof of equivalence between Theorems \ref{1thm:sign1} and  \ref{1thm:const1}]
Let $t, t'\in \bbT_n$ be $k$-adjacent.
Assume that 
the constant term of $F_{i;t}(\bfy)$
is 1.
Then, by $\eqref{1eq:Fmut1}$, the following two conditions
are equivalent:
\par
(a).  The constant term of $F_{k;t'}(\bfy)$
is 1.
\par
(b).  The $c$-vector $\bfc_{k;t}$ is either positive or negative.
\par\noindent
Therefore, Theorem \ref{1thm:sign1} follows
from Theorem \ref{1thm:const1}.
Conversely, Theorem \ref{1thm:const1}
can be shown
from Theorem \ref{1thm:sign1} and the fact $F_{i;t_0}(\bfy)=1$
by  the induction on $d=d(t_0,t)$.
\end{proof}

There are several important consequences of
Theorems  \ref{1thm:sign1}, \ref{1thm:const1},  \ref{1thm:positive1}
together with the separation formulas.

The  duality between $C$- and $G$-matrices
observed in Example \ref{1ex:typeA22}
 is a consequence of Theorem \ref{1thm:sign1}.

\begin{thm}[Duality \cite{Nakanishi11a}] 
\index{duality}
\label{1thm:dual2}
Let $D$ be a common
skew-symmetrizer of $B$-pattern $\bfB$ of $\bfSigma$.
Then, the following equality holds:
\begin{align}
\label{1eq:dual2}
D^{-1} G_t^T D C_t=I.
\end{align}
\end{thm}
\begin{proof}
Thanks to Theorem \ref{1thm:sign1},
for each $t\in \bbT_n$ and $k\in \{1,\dots,n\}$,
the sign $\varepsilon_{k;t}\in \{1, -1\}$ of the $c$-vector
$\bfc_{k;t}$ is ambiguously determined.
We set $\varepsilon=\varepsilon_{k;t}$ in 
\eqref{1eq:Ceps1} and \eqref{1eq:Geps1}.
Note that $[-\varepsilon_{k;t} c_{\ell k}]_+=0$ for any
$\ell=1,\dots, n$ due to the sign-coherence.
Then, the mutation formulas in \eqref{1eq:cmut2}
and \eqref{1eq:gmut1} are simplified as follows:
\begin{align}
 \label{1eq:cmute2}
c_{ij;t'}&=
\begin{cases}
-c_{ik;t} & j = k,\\
c_{ij;t} 
+ c_{ik;t}[\varepsilon_{k;t} b_{kj;t}]_+
& j \neq k,
\end{cases}
\\ \label{1eq:gmute1}
 g_{ij;t'}&=
 \begin{cases}
 \displaystyle
 -g_{ik;t}
 + \sum_{\ell=1}^n g_{i\ell;t} [-\varepsilon_{k;t} b_{\ell k;t}]_+
 &
 j= k,
 \\
 g_{ij;t}  &
 j \neq k.
 \end{cases}
 \end{align}
 They are  written in the matrix form
 \begin{gather}
 \label{1eq:CP1}
 C_{t'}=C_tP,
\quad
 G_{t'}=G_t Q,
\\
\small
P=\begin{pmatrix}
1 & \\
& \ddots\\
&&1\\
[\varepsilon_{k;t}b_{k1;t}]_+& \dots & [\varepsilon_{k;t}b_{k,k-1;t}]_+ & -1
& [\varepsilon_{k;t}b_{k,k+1;t}]_+&\dots& [\varepsilon_{k;t}b_{kn;t}]_+\\
&&&&1\\
&&&&&\ddots \\
&&&&&&1
\end{pmatrix} ,
\\
\small
Q=\begin{pmatrix}
1 &&&[-\varepsilon_{k;t}b_{1k;t}]_+ \\
& \ddots&&\vdots\\
&&1& [-\varepsilon_{k;t}b_{k-1,k;t}]_+\\
&  & & -1
&&& \\
&&& [-\varepsilon_{k;t}b_{k+1,k;t}]_+&1\\
&&&\vdots&&\ddots \\
&&& [-\varepsilon_{k;t}b_{nk;t}]_+&&&1
\end{pmatrix} .
 \end{gather}
 By \eqref{1eq:ss1}, we have
 \begin{align}
 DP=Q^T D.
 \end{align}
 Also, it is easy to see that
 \begin{align}
 P^2=I.
 \end{align}
 Now, we prove \eqref{1eq:dual2} by the induction on the distance $d=d(t_0,t)$.
For $t=t_0$, the claim holds by \eqref{1eq:Cmat1} and 
\eqref{1eq:Gmat1}.
Assume that the claim holds for 
$d=d(t_0,t)$.
Let $t'$ be $k$-adjacent to  $t$.
Then, we have
\begin{align}
D^{-1} G_{t'}^T D C_{t'}=
D^{-1} (Q^T G_t^T) D (C_t P)=
P(D^{-1} G_t^T D C_t)P=I.
\end{align}
\end{proof}

We have the following immediate consequences of
Theorem \ref{1thm:dual2}.
\begin{prop}[\cite{Nakanishi11a}]
\index{unimodularity}
 (a). (Unimodularity)
\begin{align}
\label{1eq:unimod1}
\det C_t = \det G_t \in \{1, -1\}.
\end{align}
\par
(b). The matrix $B_t$ is determined by $B_{t_0}$ and $C_t$ as follows:
\begin{align}
DB_t= C_t^T (DB_{t_0}) C_t.
\end{align}
\end{prop}
\begin{proof}
(a). Since $C_t$ and $G_t$ are integer matrices, 
the equality \eqref{1eq:dual2} implies the claim \eqref{1eq:unimod1}.
(This can be also shown more directly by \eqref{1eq:dual2} and the fact
$\det P = \det Q = -1$.)
(b).
By \eqref{1eq:GBBC1} and \eqref{1eq:dual2}, we have
\begin{align}
C_t^T DB_{t_0} C_t = C_t^T DG_t B_t
=DB_t.
\end{align}
\end{proof}

The following result was shown by Cao-Huang-Li \cite{Cao17} using
Theorems \ref{1thm:sign1}, \ref{1thm:const1}, \ref{1thm:positive1}, 
\ref{1thm:dual2}
and  the separation formulas all together.

\begin{thm}[{Detropicalization \cite{Cao17}}]
\label{thm:detrop1}
\index{detropicalization}
Let  $\bfSigma$ be  any cluster pattern with a given initial vertex $t_0$.
For  any  $t,t'\in \bbT_n$,
the following implications hold:
\begin{align}
G_{t}= G_{t'}
\quad
&\Longrightarrow
\quad
\bfx_{t}= \bfx_{t'},
\\ 
C_{t}= C_{t'}
\quad
&\Longrightarrow
\quad
\bfy_{t}= \bfy_{t'}.
\end{align}
\end{thm}

In other words, the cluster variables  and coefficients are uniquely
determined by their tropical parts, respectively.

The above theorem is refined and completed
as the following equivalences of periodicities.

\begin{thm}[Synchronicity \cite{Nakanishi19}]
\index{synchronicity}
\label{1thm:synchro1}
Let  $\bfSigma$ be  any cluster pattern with a given initial vertex $t_0$.
For  any  $t,t'\in \bbT_n$ and any permutation $\sigma\in S_n$,
the following four conditions are equivalent:
\begin{itemize}
\item[(a).] $G_{t}=\sigma G_{t'}$.
\item[(b).] $C_{t}=\sigma C_{t'}$.
\item[(c).] $\bfx_{t}=\sigma \bfx_{t'}$.
\item[(d).] $\Sigma_{t}=\sigma \Sigma_{t'}$.
\end{itemize}
Moreover, if 
$\bfSigma$ has free or principal coefficients at some $t'_0$,
 then the above conditions are also
equivalent to the following condition:
\begin{itemize}
\item[(e).] $\bfy_{t}=\sigma \bfy_{t'}$.
\end{itemize}
\end{thm}
As a corollary of Theorem \ref{1thm:synchro1},
we obtain the following result,
which answer the problem in CA4.

\begin{cor}[\cite{Nakanishi19}]
In  Proposition \ref{1prop:period1},
the opposite implication
in \eqref{1eq:period1}
holds for any cluster pattern $\bfSigma'$.
\end{cor}
\begin{proof}
By Theorem \ref{1thm:synchro1},
the both conditions in \eqref{1eq:period1}
is equivalent to the condition $G_{t}=\sigma G_{t'}$
 for their common $G$-matrices.

\end{proof}

\newpage
\section{Upper cluster algebras}

So far, we mainly consider the properties of seeds and cluster patterns.
In this section we study the structure of cluster algebras
via the notion of \emph{upper cluster algebras}.
All results are taken from CA3.

\subsection{Upper cluster algebras}

Upper cluster algebras was introduced in CA3.
The motivation was two-fold as follows:
\begin{itemize}
\item
They are useful to study the structure of cluster algebras.
\item
The coordinate rings of some algebraic varieties
(e.g., double Bruhat cells) are isomorphic to some upper cluster algebras,
not cluster algebras.
In other words, cluster algebras are too small in some situation.
\end{itemize}

Let $\bfSigma$ be any cluster pattern with coefficients in any
 semifield $\bbP$,
and let $\calA=\calA(\bfSigma)$ be the associated cluster algebra
in the ambient filed $\calF$.

\begin{defn}[Upper cluster algebra]
\index{cluster algebra!upper}
The \emph{upper cluster algebra\/} $\overline{\calA}=\overline{\calA}(\bfSigma)$
is a $\bbZ\bbP$-subalgebra of $\calF$
defined by
\begin{align}
\overline{\calA}=\bigcap_{t\in \bbT_n} \bbZ\bbP[\bfx_{t}^{\pm1}].
\end{align}
\end{defn}

Thanks to the Laurent phenomenon in Theorem \ref{1thm:Laurent1},
any cluster variable $x_{i;t}$ belong to $\overline{\calA}$.
Therefore,
\begin{align}
\label{1eq:AA1}
\calA \subset \overline{\calA}.
\end{align}
In some cases, $ \calA =  \overline{\calA} $ occurs, but,
in general, $ \calA \neq  \overline{\calA} $.

Let us introduce some related notions.
\begin{defn}
\label{1defn:UL1}
For each $t\in \bbT_n$, let $t_i$ ($i=1,\dots,n$) be the vertex
that is $i$-adjacent to $t$.
We define  $\bbZ\bbP$-subalgebras of $\calF$,
\begin{align}
\label{1eq:upper1}
\calU_t&=\bbZ\bbP[\bfx_t^{\pm1}] \cap \bigcap_{i=1}^n \bbZ\bbP[\bfx_{t_i}^{\pm1}],
\\
\label{1eq:lower1}
\calL_t&=\bbZ\bbP[\bfx_t, \bfx_{t_1}, \dots, \bfx_{t_n} ].
\end{align}
We call them the \emph{upper bound\/} and \emph{lower bound}
of $\calA$ at $t$, respectively.
\index{cluster algebra!upper bound}\index{upper bound}\index{cluster algebra!lower bound}\index{lower bound}
\end{defn}
Clearly, we have
\begin{align}
\label{1eq:lower2}
\calL_t \subset \calA \subset \overline{\calA} \subset \calU_t,
\end{align}
which explains their names.

We introduce the following technical condition.
\begin{defn}
For a given seed $\Sigma_t$ and $k=1,\dots,n$,
let
\begin{align}
\label{1eq:PP1}
P_{k;t}:=
\frac{ 1}{ 1\oplus y_{k;t}}
\Biggl(
 y_{k;t}
\prod_{j=1}^n
x_{j;t}^{[b_{jk;t}]_+}
+
\prod_{j=1}^n
x_{j;t}^{[-b_{jk;t}]_+}
\Biggr)
\in \bbZ\bbP[\bfx_t]
\end{align}
be the polynomial appearing in the mutation formula \eqref{1eq:xmutstandard1}.
We say that a seed $\Sigma_t$ is \emph{coprime\/} \index{coprime!for seeds}
if $P_{1;t}$, \dots, $P_{n;t}$ are pairwise coprime in $\bbZ\bbP[\bfx_t]$,
namely, any common factor of $P_{i;t}$ and $P_{j;t}$ ($i\neq j$)
belongs to $\bbZ\bbP^{\times}=\{\pm1 \}\bbP$,
where $\bbZ\bbP^{\times}$ is the unit group of the ring $\bbZ\bbP$.
\end{defn}

For our purpose, the following  sufficient condition 
for coprimeness 
is useful.
\begin{lem}
\label{1lem:coprime2}
Let $\bfSigma$ be any cluster pattern 
with free coefficients at $t_0$.
Then, $\Sigma_{t}$ is coprime for any $t\in \bbT_n$.
\end{lem}
\begin{proof}
We see in  \eqref{1eq:PP1} that,
up to a multiplicative factor $1\oplus y_{k;t}$
in $\bbZ\bbP[\bfx_t]$,
$P_{k;t}$ is a binomial of degree 1
in $y_{k;t}$.
Thus,
for $\bbP=\bbQ_{\rmsf}(\bfy_{t_0})\simeq \bbQ_{\rmsf}(\bfy_{t})$,
$P_{k;t}$ is clearly irreducible.
Therefore, $P_{1;t}$, \dots,  $P_{n;t}$ are pairwise coprime.
\end{proof}
%
%\begin{proof}
%We see in  \eqref{1eq:PP1} that,
%up to a multiplicative factor $1\oplus y_{k;t}$
%in $\bbZ\bbP[\bfx_t]$,
%$P_{k;t}$ is a binomial
%in $y_{k;t}$.
%Also, we note that $y_{1;t}$, \dots, $y_{n;t}$ are algebraically independent.
%First, we show that
%$P_{k;t}$ is irreducible in $\bbZ\bbP[\bfx_t]$.
%To see it, we define the \emph{tropical degree} for $p\in \bbP$
% by the total degree of the monomial $\pi_{\rm trop}(p)$.
%In particular, the tropical degrees of  the coefficients of two monomials in  $P_{k;t}$ are 0 and 1.
%It follows  that, if $P_{k;t}=Q P$ for some $Q,P\in \bbZ\bbP[\bfx_t]$,
% either $Q$  or $P$ should be a polynomial whose coefficients in $\bbP$
% have a common tropical degree.
% Moreover, it should actually consist of a single term in $\{\pm1\}\bbP$.
% Thus, $P_{k;t}$ is irreducible.
% It follows that
%$P_{1;t}$, \dots,  $P_{n;t}$ are pairwise coprime.
%\end{proof}

\begin{rem}
In \cite{Berenstein05}, instead of free coefficients, the following condition is used
to guarantee the coprimeness of $\Sigma_t$:
\begin{itemize}
\item
A cluster pattern $\Sigma$ is of geometric type
and its initial extended exchange matrix $\tilde{B}_{t_0}$ has
\emph{full\/} rank.
\end{itemize}
\end{rem}

\subsection{Rank 2 case}

Let us concentrate on a cluster pattern $\bfSigma$ of rank 2
with a given initial vertex $t_0$.
let $t_i$ ($i=1, 2$) be the vertex
that is $i$-adjacent to $t_0$.
Let
\begin{align}
\bfx_{t_0}=\bfx=(x_1,x_2),
\quad x'_1=x_{1;t_1},\
x'_2=x_{2;t_2},
\end{align}
so that 
\begin{align}
\label{1eq:upper2}
\calU_{t_0}&=\bbZ\bbP[x_1^{\pm1}, x_2^{\pm1}]
 \cap
\bbZ\bbP[x'_1{}^{\pm1}, x_2^{\pm1}]
 \cap
\bbZ\bbP[x_1^{\pm1}, x'_2{}^{\pm1}],
\\
\label{1eq:lower3}
\calL_{t_0}&=\bbZ\bbP[x_1, x_2, x'_1, x'_2].
\end{align}
We take the initial exchange matrix $B_{t_0}$ as in \eqref{1eq:exchange1},
where we include the case $a=b=0$ therein.
Then, the mutations at $t_0$ in directions 1 and 2 are written explicitly as
\begin{gather}
\label{1eq:relation1}
x_1 x'_1 =P_{1;t_0}=
p_1^+ x_2^a + p_1^-,
\quad
x_2 x'_2 =P_{2;t_0}=
p_2^+
+
p_2^-
x_1^{b},
\\
p_k^+=\frac{ y_k}{ 1\oplus y_k},
\quad
p_k^-=\frac{ 1}{ 1\oplus y_k},
\end{gather}
where we set $\bfy_{t_0}=\bfy$.
For $a,b>0$, it is clear that $\Sigma_{t_0}$ is coprime.
For $a=b=0$, there are cases that $\Sigma_{t_0}$ is not coprime.
\begin{ex}
Suppose that $a=b=0$ and $y_2=y_1^3$ in $\bbP$.
Then,
\begin{align}
 P_{1;t_0}=
\frac{1+y_1}{1\oplus y_1}
,
\quad
P_{2;t_0}=
\frac{1+y_1^3}{1\oplus y_1^3}
\end{align}
have a non-monomial common factor $1+y_1$.
Thus, $\Sigma_{t_0}$ is not coprime.
\end{ex}

\begin{prop}[\cite{Berenstein05}]
\label{1prop:LU1}
Let $\bfSigma$ be any cluster pattern of rank 2 
such that the initial seed $\Sigma_{t_0}$ is coprime.
Then, the following
equality holds:
\begin{align}
\label{1eq:LU1}
\calL_{t_0} = \calU_{t_0}.
\end{align}
This further implies that
\begin{align}
\label{1eq:LAU1}
\calL_{t_0}=\calA=\overline{\calA}=\calU_{t_0}.
\end{align}
\end{prop}

It is clear that \eqref{1eq:LU1} implies
\eqref{1eq:LAU1} thanks to \eqref{1eq:lower2}.

First we  present a  consequence
of Proposition \ref{1prop:LU1}.

\begin{thm}[\cite{Berenstein05}]
\label{1thm:generate1}
Let $\bfSigma$ be any cluster pattern of rank 2 with coefficients
in any semifield $\bbP$.
Then, the cluster algebra $\calA$ is generated by $x_1,x_2, x'_1, x'_2$.
In particular, it is finitely generated.
\end{thm}
\begin{proof}
First we assume that $\bfSigma$ has free coefficients
at $t_0$.
Then, by Lemma \ref{1lem:coprime2}, $\Sigma_{t_0}$
is coprime.
Thus, by Proposition \ref{1prop:LU1},
$\calA=\calL_{t_0}=\bbZ\bbP[x_1, x_2, x'_1, x'_2]$.
Therefore, the claim holds.
This means that
any cluster variable $x_{i;t}$ is expressed as a polynomial
in $x_1, x_2, x'_1, x'_2$.
This expression holds under any specialization
of coefficients including the  non-coprime case.
(Alternatively, one can  prove the claim for the case $a=b=0$,
where the non-coprime case happens,
directly  from the result in Example  \ref{1ex:A1A11}.
See also Proposition \ref{1prop:A1A11}.)
\end{proof}

Now we prove Proposition \ref{1prop:LU1}.
\begin{proof}[Proof of Proposition \ref{1prop:LU1}]
The goal is to show that \eqref{1eq:upper2}
  reduces to \eqref{1eq:lower3}.

First, we show that
\begin{align}
\label{1eq:ZZZ1}
\bbZ\bbP[x_1^{\pm1}, x_2^{\pm1}]
 \cap
\bbZ\bbP[x'_1{}^{\pm1}, x_2^{\pm1}]
=
\bbZ\bbP[x_1, x'_1, x_2^{\pm1}].
\end{align}
The inclusion $\supset$ is clear by 
\eqref{1eq:relation1}. Let us show the inclusion $\subset$.
Let $L \in \bbZ\bbP[x_1^{\pm1}, x_2^{\pm1}]$.
Then, we have
\begin{align}
\begin{split}
L&=
\sum_{m\in \bbZ} x_1^m Q_m (x_2)
\qquad
(Q_m (x_2)\in
\bbZ\bbP[x_2^{\pm1}])
\\
&=
\sum_{m\in \bbZ} 
x'_1{}^{-m}
(p_1^+ x_2^a +p_1^-)^m
 Q_m (x_2)
,
\end{split}
\end{align}
where the sum is finite.
Imposing that $L$ also belongs to 
$\bbZ\bbP[x'_1{}^{\pm1}, x_2^{\pm1}]$,
we have, for any $m<0$,
\begin{align}
 (p_1^+ x_2^a +p_1^-)^m
 Q_m (x_2) \in \bbZ\bbP[x_2^{\pm1}].
 \end{align}
 Then,
 \begin{align}
L&=
\sum_{m\geq 0} x_1^m Q_m (x_2)
+
\sum_{m< 0} x'_1{}^{-m}
 (p_1^+ x_2^a +p_1^-)^{m}
 Q_{m} (x_2)
 \in
 \bbZ\bbP[x_1, x'_1, x_2^{\pm1}],
\end{align} 
which proves \eqref{1eq:ZZZ1}.
By \eqref{1eq:ZZZ1} and the one obtained by interchanging $x_1$ and $x_2$,
we obtain
the first step reduction
\begin{align}
\label{1eq:upper3}
\calU_{t_0}&=\bbZ\bbP[x_1, x'_1, x_2^{\pm1}]
 \cap
\bbZ\bbP[x_2, x'_2, x_1^{\pm1}].
\end{align}

Now our goal is to prove the following claim.
\par\noindent
{\bf Claim.} The following equality holds:
\begin{align}
\label{1eq:upper4}
\bbZ\bbP[x_1, x'_1, x_2^{\pm1}]
 \cap
\bbZ\bbP[x_2, x'_2, x_1^{\pm1}]
=
\bbZ\bbP[x_1, x_2, x'_1, x'_2].
\end{align}
The inclusion $\supset$ is clear.
We separate the proof
of the inclusion $\subset$
 into two cases.

{\bf Case 1: $a,b> 0$.}
First we prove the following equality.
\begin{align}
\label{1eq:upper5}
\bbZ\bbP[x_1, x'_1, x_2^{\pm1}]
=
\bbZ\bbP[x_1, x_2, x'_1, x'_2]
+
\bbZ\bbP[x_1, x_2^{\pm1}].
\end{align}
The inclusion $\supset$ is clear.
Let us show the inclusion $\subset$.
It is enough to prove that,
for $k, \ell > 0 $,
\begin{align}
\label{1eq:ZZ1}
 x'_1{}^k x_2^{-\ell} \in 
\bbZ\bbP[x_1, x'_1, x'_2]
+
\bbZ\bbP[x_1, x_2^{\pm1}].
\end{align}
Let $q:=-p_2^-/p_2^+$.
Then,  by \eqref{1eq:relation1}, we have
\begin{align}
x_2^{-1}\equiv qx_1^b x_2^{-1}
\mod \bbZ\bbP [x_2'].
\end{align}
Applying it repeatedly, we have, for any $k>0$,
\begin{align}
x_2^{-1}
\equiv
qx_1^b x_2^{-1}
\equiv
q^2x_1^{2b} x_2^{-1}
\equiv
\cdots
\equiv q^{k}x_1^{kb} x_2^{-1}
\mod \bbZ\bbP [x_1, x_2'].
\end{align}
Therefore, 
\begin{align}
x_2^{-\ell}\in 
\bbZ\bbP[x_1, x'_2]
+
x_1^{k}\bbZ\bbP[x_1, x_2^{\pm1}],
\end{align}
where we used the assumption $b>0$.
Noticing that $x'_1 x_1\in \bbZ\bbP[x_2]$,
we obtain \eqref{1eq:ZZ1}.
Now, applying the equality \eqref{1eq:upper5}
to the left hand side of \eqref{1eq:upper4}, we have 
\begin{align}
\label{1eq:upper6}
\begin{split}
&\quad\,
\bbZ\bbP[x_1, x'_1, x_2^{\pm1}]
 \cap
\bbZ\bbP[x_2, x'_2, x_1^{\pm1}]\\
&=
(\bbZ\bbP[x_1, x_2, x'_1, x'_2]
+
\bbZ\bbP[x_1, x_2^{\pm1}])
 \cap
\bbZ\bbP[x_2, x'_2, x_1^{\pm1}]\\
&=
\bbZ\bbP[x_1, x_2, x'_1, x'_2]
+
(\bbZ\bbP[x_1, x_2^{\pm1}]
 \cap
\bbZ\bbP[x_2, x'_2, x_1^{\pm1}]).
\end{split}
\end{align}
Then, the claim \eqref{1eq:upper4} follows from  the following equality:
\begin{align}
\label{1eq:upper7}
(\bbZ\bbP[x_1, x_2^{\pm1}]
 \cap
\bbZ\bbP[x_2, x'_2, x_1^{\pm1}])
=
\bbZ\bbP[x_1, x_2, x'_2].
\end{align}
The inclusion $\supset$ is clear.
We show the inclusion $\subset$.
Let $L\in \bbZ\bbP[x_2, x'_2, x_1^{\pm1}]$.
We have
\begin{align}
\begin{split}
L=\sum_{m,k,l}
c_{mk\ell}
x_1^m x_2^k x'_2{}^{\ell}
&=\sum_{m,k,l} 
c_{mk\ell}
x_1^m x_2^k x_2^{-\ell}(p_2^+ + p_2^- x_1^b)^{\ell},
\\
&\qquad
(m\in \bbZ,\ k, \ell \geq 0,\ c_{mk\ell}\in \bbZ\bbP).
\end{split}
\end{align}
Imposing that  $L$ also belong to $\bbZ\bbP[x_1, x_2^{\pm1}]$,
we see that
the powers $m$ in the above should be nonnegative.
This is the desired result.

{\bf Case 2: $a=b=0$.}
In this case, one can directly show
\eqref{1eq:upper4}.
We only need to show the inclusion $\subset$.
Let $L \in \bbZ\bbP[x_2, x'_2, x_1^{\pm1}]$.
Then, we have
\begin{align}
\begin{split}
L=\sum_{m,k,l}
c_{mk\ell}
x_1^m x_2^k x'_2{}^{\ell}
&=\sum_{m,k,l} 
c_{mk\ell}
x_1^m x_2^k x_2^{-\ell}(p_2^+ + p_2^-)^{\ell},
\\
&\qquad
(c_{mk\ell}\in \bbZ\bbP,
\ m\in \bbZ,\ k, \ell \geq 0).
\end{split}
\end{align}
A similar expression holds for $L'\in \bbZ\bbP[x_1, x'_1, x_2^{\pm1}]$.
It follows that
\begin{align}
L''=\sum_{m,k}
c_{mk}
x_1^m x_2^k \in \bbZ\bbP[x_1^{\pm 1}, x_2^{\pm1}],
\quad (c_{mk}\in \bbZ\bbP)
\end{align}
belongs to the left hand side of \eqref{1eq:upper4}
if the following condition is satisfied:
\begin{itemize}
\item[(i).]
$c_{mk}$ is divisible by $(p_1^+ + p_1^-)^{-m}$ if $m<0$,
\item[(ii).]
$c_{mk}$ is divisible by $(p_2^+ + p_2^-)^{-k}$ if $k<0$.
\end{itemize}
Here, we recall that $P_{1;t_0}=p_1^+ + p_1^-$ and $P_{2;t_0}=p_2^+ + p_2^-$
in \eqref{1eq:relation1}
are coprime in $\bbZ\bbP$ by assumption.
(This is the only place where the coprimeness condition
is concerned.)
Then, the conditions  (i) and (ii)  ensure the following condition,
\begin{itemize}
\item[(iii).]
$c_{mk}$ is divisible by $(p_1^+ + p_1^-)^{-m}(p_2^+ + p_2^-)^{-k}$
if $m,k<0$.
\end{itemize}
The conditions (i)--(iii) guarantee that $L''$ belongs to
$\bbZ\bbP[x_1, x_2, x'_1, x'_2]$
after replacing negative powers of $x_i$
with positive power of $x'_i/(p_i^+ + p_i^-)$.

This completes the proof of Proposition \ref{1prop:LU1}.
\end{proof}

In the non-coprime case, one can 
 directly verify the following fact,
  using  the result in  Example \ref{1ex:A1A11},
\begin{prop}
\label{1prop:A1A11}
Let $\bfSigma$ be any cluster pattern of rank 2 
such that the initial seed $\Sigma_{t_0}$ is not coprime.
Then,
we have
\begin{align}
\label{1eq:LAU2}
\calL_{t_0}=\calA=\overline{\calA}\subsetneqq\calU_{t_0}.
\end{align}
\end{prop}
\begin{proof}
Since $\Sigma_{t_0}$ is not coprime, 
we have $a=b=0$.
Then,
as shown in   Example \ref{1ex:A1A11},
there are four clusters
 $(x_1,x_2)$, $(x'_1,x_2)$, $(x_1,x'_2)$, $(x'_1,x'_2)$,
 where
\begin{align}
\label{1eq:A1A13}
x_1x'_1=P_{1;t_0}=p_1^++p_1^-,
\quad
x_2 x'_2=P_{2;t_0}=p_2^++p_2^-.
\end{align}
Therefore, $\calL_{t_0}=\bbZ\bbP[x_1,x_2,x'_1,x'_2]=\calA$ holds.
By assumption, there is
some  common factor $R\in \bbZ\bbP$ of $P_{1;t_0}$ and $P_{2;t_0}$
 with $R\not\in \bbZ\bbP^{\times}$ 
such that
\begin{align}
P_{1;t_0}
=Q_1R,
\quad
P_{2;t_0}
=Q_2R
\quad
(Q_1, Q_2\in \bbZ\bbP).
\end{align}
We have
\begin{align}
x_1^{-1} x_2^{-1}Q_1Q_2 R
=x_1'x_2^{-1}Q_2=x_1^{-1}x_2'Q_1
\in \calU_{t_0}.
\end{align}
On the other hand, this element does not belong to $\calA$,
because one cannot eliminate the negative powers of
$x_1$ and $x_2$ simultaneously.
Therefore, $\calA\neq \calU_{t_0}$.
On the other hand,
any element $L$ of $\overline{\calA}$ is explicitly written as
($c_{mn}\in \bbZ\bbP$)
\begin{align}
\begin{split}
L&=\sum_{m,n\geq 0} c_{mn} x_1^m x_2^n
+
\sum_{m\geq 0, n<0}  c_{mn}  P_{2;t_0}^{-n} x_1^m x_2^n
\\
&\qquad
+
\sum_{m<0, n\geq 0}  c_{mn}  P_{1;t_0}^{-m}  x_1^m x_2^n
+
\sum_{m,n<0} c_{mn} P_{1;t_0}^{-m} P_{2;t_0}^{-n} x_1^m x_2^n
\\
&=
\sum_{m,n\geq 0} c_{mn} x_1^m x_2^n
+
\sum_{m\geq 0, n<0}  c_{mn} x_1^m x'_2{}^{-n}
\\
&\qquad
+
\sum_{m<0, n\geq 0} c_{mn} x'_1{}^{-m} x_2^{n}
+
\sum_{m,n<0}   c_{mn} x'_1{}^{-m} x'_2{}^{-n}
\in \calA.
\end{split}
\end{align}
Therefore, $\calA=\overline{\calA}$.
\end{proof}

\subsection{Alternative proof of Laurent phenomenon}

As another application
of Proposition \ref{1prop:LU1} (and the results in the proof),
we present alternative proof of the Laurent phenomenon
in Section \ref{1subsec:Laurent1}.

First, we prove the following  fact.

\begin{prop}
\label{1prop:LU3}
Let $\bfSigma$ be any cluster pattern of rank 2 
such that any seed $\Sigma_t$ is coprime.
Then,
the following
equality holds for any $t\in \bbT_2$:
\begin{align}
\label{1eq:LU4}
 \calU_{t}=
\calL_{t}=
  \calU_{t_0}=
\calL_{t_0}.
\end{align}
In particular, the upper bound $ \calU_{t}$ is independent of $t$.
\end{prop}
\begin{proof}
Since any seed $\Sigma_t$ is coprime by assumption,
Proposition \ref{1prop:LU1} is applicable to any $t$,
so that we have $ \calL_{t}=
\calU_{t}$.
Therefore, it is enough to prove
\begin{align}
\label{1eq:LL1}
 \calL_{t_0}=
\calL_{t_1},
\end{align}
where
\begin{align}
\calL_{t_1}=\bbZ\bbP[x_1,x_2,x'_1,x''_2],
\end{align}
and $x''_2$ is the mutation of $x_2=x_{2;t_1}$ at $t_1$ in direction $2$.
Moreover,
by the symmetry of $x'_2$ and $x''_2$, it is enough to prove
that
\begin{align}
\label{1eq:x21}
x''_2\in  \calL_{t_0}=\bbZ\bbP[x_1,x_2,x'_1,x'_2].
\end{align}
We apply \eqref{1eq:sodd1} with $s=1$ therein to obtain
\begin{gather}
\label{1eq:relation2}
x_2 x''_2 =
q_2^+
x'_1{}^{b}
+
q_2^-,
\\
q_2^+=\frac{ y'_2}{ 1\oplus y'_2},
\quad
q_2^-=\frac{ 1}{ 1\oplus y'_2},
\quad
y'_2=y_2(1\oplus y_1)^b.
\end{gather}
The following relations hold:
\begin{align}
q_2^+/q_2^-&=y'_2=(p_2^+/p_2^- )(p_1^-)^{-b},
\\
x_1 x'_1&=  p_1^+ x_2^{a} + p_1^-,
\\
x_2 x'_2&=  p_2^+ + p_2^-x_1^{b}.
\end{align}
It follows that
\begin{align}
\begin{split}
x''_2
&=\frac{q_2^+x'_1{}^{b} + q_2^-}{x_2}
=\frac{q_2^+x'_1{}^{b}(x_2 x'_2 -p_2^-x_1^{b})}{p_2^+ x_2}
+\frac{ q_2^-}{x_2}\\
&
=
\frac{q_2^+x'_1{}^{b} x'_2 }{p_2^+ }
-
\frac{q_2^+p_2^-}{p_2^+ x_2}
\left(
(x'_1x_1)^{b}
- \frac{ q_2^- p_2^+ }{q_2^+p_2^-}
\right)
\\
&
=
\frac{q_2^+x'_1{}^{b} x'_2 }{p_2^+ }
-
\frac{q_2^+p_2^-}{p_2^+ x_2}
\left(
(p_1^+ x_2^a + p_1^-)^{b}
- (p_1^-)^b
\right).
\end{split}
\end{align}
The numerator of the second term in the last expression is
divisible by $x_2$.
Thus, we have  the property \eqref{1eq:x21}.
\end{proof}

We once again prove the Laurent phenomenon
in Theorem \ref{1thm:Laurent1}
based on Proposition
\ref{1prop:LU3}.
\index{Laurent phenomenon}

\begin{proof}[Alternative proof of Theorem \ref{1thm:Laurent1}]
Let $\bfSigma$ be any cluster pattern with rank $n$.
As in the previous proof in Section \ref{1subsec:Laurent1},
 it is enough to prove Theorem \ref{1thm:Laurent1} 
 assuming that $\bfSigma$ has free coefficients at $t_0$.
 Then, by Lemma \ref{1lem:coprime2},
any seed $\Sigma$ is coprime.
In particular,
the  equality \eqref{1eq:LU4}
is applicable
 for any rank 2 restriction of $\bfSigma$,
 where the frozen variables are regarded as a part of coefficients.
 For example, when we consider a cluster cluster pattern restricted in direction 1 and 2,
 we replace $\bbZ\bbP$ with $ \bbZ\bbP[x_3^{\pm1}, \dots, x_n^{\pm1}]$.
Then, in the same notation in the proof of  Proposition \ref{1prop:LU3},
 we have the following equality corresponding to  \eqref{1eq:LL1}:
   \begin{align}
( \bbZ\bbP[ x_3^{\pm1}, \dots, x_n^{\pm1}])
 [x_1,x_2,x'_1,x'_2]
 =
( \bbZ\bbP[ x_3^{\pm1}, \dots, x_n^{\pm1}])[x_1,x_2,x'_1,x''_2],
 \end{align}
which is also written as
  \begin{align}
 \label{1eq:LL2}
 \bbZ\bbP[x_1,x_2,x'_1,x'_2,x_3^{\pm1}, \dots, x_n^{\pm1}]
 =
 \bbZ\bbP[x_1,x_2,x'_1,x''_2,x_3^{\pm1}, \dots, x_n^{\pm1}].
 \end{align}

 Let us prove that  the upper bound $ \calU_{t}$ of $\bfSigma$ is independent of $t$.
 Let $t_i$ ($i=1,\dots,n$) be the vertex
that is $i$-adjacent to $t_0$.
Let $x'_i=x_{i;t_i}$.
It is enough to prove $\calU_{t_0}=\calU_{t_1}$.
Indeed, we have
\begin{align}
\begin{split}
\calU_{t_0}
&=\bbZ\bbP[\bfx^{\pm1}] \cap \bigcap_{i=1}^n \bbZ\bbP[\bfx_{t_i}^{\pm1}]
\\
&= \bigcap_{i=1}^n \left(\bbZ\bbP[\bfx^{\pm1}] \cap \bbZ\bbP[\bfx_{t_i}^{\pm1}]
\right)
\\
&=
\bigcap_{i=1}^n  \bbZ\bbP[ x_1^{\pm1}, \dots, x_{i-1}^{\pm1}, x_i, x'_i, x_{i+1}^{\pm1}, \dots, x_n^{\pm1}]
\quad
\text{(by \eqref{1eq:ZZZ1})}
\\
&=
\bigcap_{i=2}^n 
(
\bbZ\bbP[ x_1, x'_1, x_2^{\pm1} , \dots, x_n^{\pm1}]
\cap
 \bbZ\bbP[ x_1^{\pm1}, \dots, x_{i-1}^{\pm1}, x_i, x'_i, x_{i+1}^{\pm1}, \dots, x_n^{\pm1}]
 )
 \\
 &=
\bigcap_{i=2}^n 
\bbZ\bbP[ x_1, x'_1, x_2^{\pm1} , \dots,x_{i-1}^{\pm1},  x_i, x'_i, x_{i+1}^{\pm1}, \dots, x_n^{\pm1}]
 \\
 &=
\bigcap_{i=2}^n 
\bbZ\bbP[ x_1, x'_1, x_2^{\pm1} , \dots, x_{i-1}^{\pm1}, x_i, x''_i, x_{i+1}^{\pm1}, \dots, x_n^{\pm1}]
\quad
\text{(by \eqref{1eq:LL2})}
\\
&=
\calU_{t_1},
\end{split}
\end{align}
where the last equality is obtained by reversing the above procedure.
Then, the upper cluster algebra $\overline{\calA}$ is given by
\begin{align}
\overline{\calA}=\bigcap_{t\in \bbT_n} \bbZ\bbP[\bfx_{t}^{\pm1}]
=\bigcap_{t\in \bbT_n}\calU_{t}
=\calU_{t_0}.
\end{align}
On the other hand, for the cluster algebra $\calA$, we have
\begin{align}
\calA=\bbZ\bbP[\bfx_t \mid t\in \bbT_n]
\subset \bigcup_{t\in \bbT_n}\calU_{t}
=\calU_{t_0}.
\end{align}
Therefore, we obtain $\calA\subset \overline{\calA}$.
\end{proof}

\subsection{Further results}
Let us present
generalizations of 
Proposition \ref{1prop:LU1}
and 
Theorem \ref{1thm:generate1},
which are part of the main result in CA3, without proof.

\begin{defn}
An $n\times n$ skew-symmetrizable matrix $B$ is \emph{acyclic\/} \index{matrix!acyclic}
if the following condition holds:
\begin{itemize}
\item
For any $3\leq m \leq n$,
there is no cyclic sequence $i_1$, $i_2$, \dots, $i_m$, $i_{m+1}=i_1$
in $\{1,\dots,n\}$
such that $b_{i_s i_{s+1}}>0$ holds for any $s=1,\dots,m$.
\end{itemize}
\end{defn}

\begin{ex}
\label{1ex:acyclic1}
(a). Any $2\times 2$ skew-symmetrizable matrix $B$ is acyclic.
\par
(b).
Let $B$ be any skew-symmetrizable matrix such that 
the associated Cartan matrix $A(B)$  is of finite type
given in Definition \ref{1defn:Cartan2}.
Then, $B$ is acyclic because the corresponding Dynkin diagram
is a tree.
\end{ex}

By a similar technique using upper and lower bounds, a parallel result to
 Proposition \ref{1prop:LU1} and
 Theorem \ref{1thm:generate1} hold for any cluster pattern with
an acyclic initial exchange matrix.

\begin{thm}[\cite{Berenstein05}]
\label{1thm:generate2}
Let $\bfSigma$ be a cluster pattern 
 with coefficients
in any semifield $\bbP$ with a given initial vertex $t_0$.
Suppose that the initial exchange matrix $B_{t_0}$ is acyclic.
Then, the following facts hold.
\par
(a).  
We have
\begin{align}
\calL_{t_0}=\calA.
\end{align}
In particular,
the cluster algebra $\calA$ is finitely generated.
\par
(b). If the initial seed $\Sigma_{t_0}$ is coprime,
we have
\begin{align}
\label{1eq:LAU3}
\calL_{t_0}=\calA=\overline{\calA}=\calU_{t_0}.
\end{align}
\end{thm}

\newpage
\section{Generalized cluster algebras}

To conclude this introductory guide, we present a generalization of cluster algebras
introduced by Chekhov and Shapiro
as an extra material beyond basics.
Essentially all nice properties of cluster algebras are inherited to this generalization.
Therefore, they extensively widen the perspective of cluster algebra theory.

\subsection{Generalized cluster algebras}

Chekhov and Shapiro \cite{Chekhov11} introduced a  generalization of
cluster algebras  called  \emph{generalized cluster algebras\/} (GCA, for short),
motivated by
some examples naturally appeared in the study of Teichm\"uller space for
the Riemann surface with orbifold points.
Moreover, it turned out that
they are indeed a natural generalization such that
all essential properties of cluster algebras presented in this text are shown or conjectured to  hold
\cite{Chekhov11, Nakanishi14a}.

Let us explain the idea shortly.
The mutations of ordinary cluster algebras
are defined by  \emph{binomials\/}
$1+\haty_k$
and
$1\oplus y_k$ in   \eqref{1eq:xmut1} and \eqref{1eq:ymut1}.
We may ask how/why this binomial condition is essential.
It turns out that one can replace it  with  fairly general
 \emph{polynomials},
and they are still  as good as the ordinary one, \emph{if\/} we adequately
adjust the other parts of mutations.

\begin{defn}[Mutation data]
\index{generalized cluster algebra (GCA)!mutation data}
Let $\bbP$ be a given semifield.
We fix the following \emph{mutation data\/} $(\bfr,\bfz)$ in $\bbP$:
\begin{itemize}
\item an $n$-tuple of positive integers
$\bfr=(r_1,\dots,r_n)$, which are called \emph{mutation degrees}.
\item a collection
$\bfz=(z_{i,s})_{i=1,\dots,n; s=0,\dots, r_i}$
 of elements in $\bbP$ satisfying the following conditions:
\begin{align}
\label{1eq:gconst1}
\mbox{(unit constant term/monic property)}
\quad
z_{i,0}&=z_{i,r_i}=1,
%\hskip20pt
\\
\label{1eq:rec1}
\mbox{(reciprocity)}\quad z_{i,s}&=z_{i,r_i-s}.
\end{align}
\end{itemize}
They are equivalent to give an $n$-tuple of \emph{monic\/} and \emph{reciprocal\/}
polynomials in a single formal
variable $y$ with coefficients in $\bbP$,
\begin{align}
\label{1eq:P1}
P_{i,\bfz}(y)=
\sum_{s=0}^{r_i} z_{i,s} y^s,
\quad (i=1,\dots,n).
\end{align}
In the simplest case $\bfr=(1,\dots,1)$, they reduce to $P_{i,\bfz}(y)=1+y$, which is the
 binomial for the ordinary mutations of  seeds
 mentioned above.
\end{defn}

\begin{defn}[Seed mutation for GCA]
\index{generalized cluster algebra (GCA)!mutation}
Let $(\bfr,\bfz)$ be a mutation data in $\bbP$.
Let $\Sigma=(\bfx,\bfy,B)$ be an (ordinary) seed with coefficients in $\bbP$ in Definition
\ref{1defn:seed1}.
For the  polynomials in \eqref{1eq:P1}, let
\begin{align}
P_{k,\bfz}(\haty_k)=
\sum_{s=0}^{r_k} z_{k,s}\haty_k^s,
\quad
P_{k,\bfz}\vert_{\bbP}(y_k)=
\bigoplus_{s=0}^{r_k} z_{k,s}y_k^s,
\end{align}
where $\haty_k$ is the ordinary $\haty$-variable for $\Sigma$.
Then, the \emph{$(\bfr,\bfz)$-mutation of $\Sigma$
in direction $k$}, $\mu_k(\Sigma)=(\bfx',\bfy',B')$,
 is defined as follows:
\begin{align}
\label{1eq:gxmut1}
x'_i
&=
\begin{cases}
\displaystyle
x_k^{-1}
\Biggl(\, \prod_{j=1}^n x_j^{[-b_{jk}r_k]_+}
\Biggr)
\frac{ P_{k,\bfz}(\haty_k)}{ P_{k,\bfz}\vert_{\bbP}(y_k)}
& i=k,
\\
x_i
&i\neq k,
\end{cases}
\\
\label{1eq:gymut1}
y'_i
&=
\begin{cases}
\displaystyle
y_k^{-1}
& i=k,
\\
y_i y_k^{[r_k b_{ki}]_+} P_{k,\bfz}\vert_{\bbP}(y_k)^{-b_{ki}}
&i\neq k,
\end{cases}
\\
\label{1eq:gbmut2}
b'_{ij}&=
\begin{cases}
-b_{ij}
&
\text{$i=k$ or $j=k$,}
\\
b_{ij}
+
b_{ik} [r_kb_{kj}]_+
+
[-b_{ik}r_k]_+  b_{kj}
&
i,j\neq k.
\end{cases}
\end{align}
\end{defn}

It is easy to check the following properties.

\begin{prop}
(a). Any skew-symmetrizer $D$ of $B$ is also
a skew-sym\-metrizer of $B'$.
\par
(b). The mutation $\mu_k$ is an involution.
In particular, $(\bfx',\bfy',B')$ is a seed.
\par
(c). The $\haty$-variables mutate as
\begin{align}
\label{1eq:gyhatmut1}
\haty'_i
&=
\begin{cases}
\displaystyle
\haty_k^{-1}
& i=k,
\\
\haty_i \haty_k^{[r_k b_{ki}]_+} P_{k,\bfz}(\haty_k)^{-b_{ki}}
&i\neq k.
\end{cases}
\end{align}
\par
(d). The $\varepsilon$-expressions, which are parallel to 
the ones in Proposition \ref{1prop:epsilon1}, hold,
where we put $\varepsilon$ in every $[\, \cdot \, ]_+$
and also replace $P_{k,\bfz}(\haty_k)$ and $P_{k,\bfz}\vert_{\bbP}(y_k)$
with $P_{k,\bfz}(\haty_k^{\varepsilon})$ and $P_{k,\bfz}\vert_{\bbP}(y_k^{\varepsilon})$
in \eqref{1eq:gxmut1}--\eqref{1eq:gbmut2}.

\end{prop}
\begin{proof}
(a). This is proved in the same way as
Proposition \ref{1prop:Bmut1} (a).
\par
(b). 
By  the reciprocity condition \eqref{1eq:rec1},
we have a parallel formula to \eqref{1eq:yy1},
\begin{align}
\label{1eq:gyy1}
\frac{ P_{k,\bfz}(\haty_k^{-1})}{ P_{k,\bfz}\vert_{\bbP}(y_k^{-1})}
=
\frac{P_{k,\bfz}(\haty_k)}{  P_{k,\bfz}\vert_{\bbP}(y_k)}
\prod_{j=1}^n  x_j^{-b_{jk}r_k}.
\end{align}
Then, one can repeat the proofs of Propositions \ref{1prop:Bmut1} (c)
and  \ref{1prop:smut1} (a).
\par
The  properties (c) and (d) can be proved in a similar way as
Propositions \ref{1prop:yhat1} and \ref{1prop:epsilon1}.
\end{proof}

\begin{rem}
One can lift the reciprocity condition \eqref{1eq:rec1} by
putting the data $\bfz$ in a seed and introducing its
mutation \cite{Nakanishi15}.
\end{rem}

\begin{defn}
For a given mutation data $(\bfr,\bfz)$ in $\bbP$,
a \emph{generalized cluster pattern\/}   \index{generalized cluster algebra (GCA)!cluster pattern}
 $\bfSigma$ with coefficients in $\bbP$
are defined in the same way as  the  ordinary one
by replacing  mutations with  $(\bfr,\bfz)$-mutations.
The \emph{generalized cluster algebra\/} $\calA$ associated with $\bfSigma$ is
 \index{generalized cluster algebra (GCA)}\index{GCA|see{generalized cluster algebra}}
defined in the same way as  the ordinary one.
\end{defn}

The first serious test for this generalization is the Laurent phenomenon,
which is the \emph{raison d'etre\/} of cluster algebras.
This indeed holds literally in the same way as the ordinary ones.

\begin{thm}[Laurent phenomenon \cite{Chekhov11}]
\index{generalized cluster algebra (GCA)!Laurent phenomenon}
\label{1thm:gLaurent1}
Let $\bfSigma$ be any generalized cluster pattern with coefficients in any semifield $\bbP$.
Let $t_0,t\in \bbT_n$ be any vertices.
Then, any cluster variable $x_{i;t}$ is expressed as a Laurent polynomial in 
$\bfx_{t_0}$ with coefficients in $\bbZ\bbP$.
\end{thm}
\begin{proof}
The proof is also the same as the proof of Theorem \ref{1thm:Laurent1},
just by replacing the binomial $1+\haty_k$ therein with 
the  polynomial $P_{k,\bfz}(\haty_k)$.
\end{proof}

\begin{ex}
Let us consider a rank 2 example.
As the simplest nontrivial mutation data $(\bfr,\bfz)$,
we take $\bfr=(2,1)$, where the only nontrivial data in $\bfz$ is $z_{1,1}$.
We choose $z_{1,1}=z$  to be arbitrary in $\bbP$.
We also choose the simplest nontrivial initial exchange matrix
\begin{align}
B_0=
\begin{pmatrix}
0 & -1 \\
1 & 0\\
\end{pmatrix}.
\end{align}
Accordingly,
\begin{align}
\hat{y}_1=y_1x_2,
\quad
\hat{y}_2=y_2x_1^{-1}.
\end{align}
We note that, for the diagonal matrix $R=\mathrm{diag}(r_1,r_2)$,
\begin{align}
B_0R=
\begin{pmatrix}
0 & -1 \\
2 & 0\\
\end{pmatrix},
\end{align}
which is the exchange matrix for an
 ordinary cluster algebra of type $B_2$
  we studied in Section \ref{1sec:rank2}.
We repeat a similar calculation as in Section
\ref{1sec:rank2}.
The result is given below up to $t=6$.
\begin{align}
&
\begin{cases}
x_{1;0}=x_1,\\
x_{2;0}=x_2,\\
\end{cases}
\qquad
\begin{cases}
y_{1;0}=y_1,\\
y_{2;0}=y_2,\\
\end{cases}
\\
\allowbreak
&
\begin{cases}
\displaystyle
x_{1;1}=x_1^{-1}\frac{1+z\hat{y}_1+\hat{y}_1^2}{1\oplus z y_1\oplus y_1^2},\\
x_{2;1}=x_2,\\
\end{cases}
\qquad
\begin{cases}
y_{1;1}=y_1^{-1},\\
y_{2;1}=y_2(1\oplus z y_1\oplus y_1^2),\\
\end{cases}
\\
    \begin{split}
&
\begin{cases}
\displaystyle
x_{1;2}=x_1^{-1}\frac{1+z\hat{y}_1+\hat{y}_1^2}{1\oplus z y_1\oplus y_1^2},\\
\displaystyle
x_{2;2}=x_2^{-1}
\frac{1 + \hat{y}_2 +z \hat{y}_1\hat{y}_2+
\hat{y}_1^2 \hat{y}_2}
{1\oplus y_2\oplus z y_1y_2
\oplus y_1^2y_2},
    \rule{0pt}{18pt}\hskip-1pt
    \\
\end{cases}
\\
&
\begin{cases}
y_{1;2}=y_1^{-1}(1\oplus y_2 \oplus z y_1y_2\oplus y_1^2 y_2),\\
y_{2;2}=y_2^{-1}(1\oplus z y_1\oplus y_1^2)^{-1},\\
\end{cases}
\end{split}
\\
\allowbreak
    \begin{split}
&
\begin{cases}
\displaystyle
x_{1;3}=x_1x_2^{-2}
\frac{
1+ 2 \hat{y}_2+ \hat{y}_2^2
 + z\hat{y}_1\hat{y}_2+ z\hat{y}_1\hat{y}_2^2 + \hat{y}_1^2\hat{y}_2^2}
{1\oplus 2 y_2\oplus y_2^2
 \oplus zy_1y_2\oplus zy_1y_2^2 \oplus y_1^2y_2^2},
\\
\displaystyle
x_{2;3}=x_2^{-1}
\frac{1 + \hat{y}_2 +z \hat{y}_1\hat{y}_2+
\hat{y}_1^2 \hat{y}_2}
{1\oplus y_2\oplus z y_1y_2
\oplus y_1^2y_2},
    \rule{0pt}{18pt}\hskip-1pt
\\
\end{cases}
\\
\allowbreak
&
\begin{cases}
y_{1;3}=y_1(1\oplus y_2 \oplus z y_1y_2\oplus y_1^2 y_2)^{-1},\\
y_{2;3}=y_1^{-2}y_2^{-1}(1\oplus 2 y_2\oplus y_2^2
 \oplus zy_1y_2\oplus zy_1y_2^2 \oplus y_1^2y_2^2
),
\\
\end{cases}
\end{split}
\\
    \begin{split}
&
\begin{cases}
\displaystyle
x_{1;4}=x_1x_2^{-2}
\frac{
1+ 2 \hat{y}_2+ \hat{y}_2^2
 + z\hat{y}_1\hat{y}_2+ z\hat{y}_1\hat{y}_2^2 + \hat{y}_1^2\hat{y}_2^2}
{1\oplus 2 y_2\oplus y_2^2
 \oplus zy_1y_2\oplus zy_1y_2^2 \oplus y_1^2y_2^2},
\\
\displaystyle
x_{2;4}=
x_1x_2^{-1}
\frac{1 + \hat{y}_2 }{1\oplus y_2},
    \rule{0pt}{18pt}\hskip-1pt
\\
\end{cases}
\\
&
\begin{cases}
y_{1;4}=y_1^{-1}y_2^{-1}(1\oplus y_2),\\
y_{2;4}=y_1^{2}y_2(1\oplus 2 y_2\oplus y_2^2
\oplus zy_1y_2\oplus zy_1y_2^2 \oplus y_1^2y_2^2
)^{-1},\\
\end{cases}
\end{split}
\\
\allowbreak
&
\begin{cases}
\displaystyle
x_{1;5}=x_1,\\
\displaystyle
x_{1;5}=x_1x_2^{-1}\frac{1+\hat{y}_2}{1\oplus y_2},\\
\end{cases}
\qquad
\begin{cases}
\displaystyle
y_{1;5}=y_1y_2(1\oplus y_2)^{-1},\\
y_{1;5}=y_2^{-1},\\
\end{cases}
\\
\allowbreak
&
\begin{cases}
x_{1;6}=x_1,\\
x_{2;6}=x_2,\\
\end{cases}
\qquad
\begin{cases}
y_{1;6}=y_1,\\
y_{2;6}=y_2.\\
\end{cases}
\end{align}
We observe that all characteristic feature of
ordinary cluster  patterns are preserved.
In particular, we see the same periodicity of
the ordinary cluster pattern of type $B_2$
 in Section
\ref{1sec:rank2}.
In fact, if we formally set $z=0$ in the above
(though such a specialization is prohibited in $\bbP$)
we recover the result for the ordinary cluster pattern of type $B_2$ therein.
This is an example of  the general  fact 
that any generalized cluster pattern reduces to 
some (ordinary) cluster pattern called the \emph{right companion cluster pattern\/}
under such a specialization
\cite{Nakanishi15}.

\end{ex}

\subsection{Separation formulas for GCA}
As one more confirmation of the rightness of GCA,
we present the separation formulas for generalized cluster patterns,
which are parallel to the ordinary ones.

\begin{defn}
\label{1defn:gprincipal1}
We say that a generalized cluster pattern $\bfSigma$
with mutation data $(\bfr,\bfz)$
is \emph{with principal coefficients at $t_0\in \bbT_n$\/}
if the following conditions are satisfied:
\begin{itemize}
\item
The coefficient semifield of $\bfSigma$ is a tropical semifield
$\mathrm{Trop}(\bfy,\bfz)$ with generators $\bfy=(y_1,\dots,y_n)$,
$\bfz=(z_{i,s})_{i=1,\dots,n;\, s=1,\dots, r_i-1}$
with $z_{i,s}=z_{i,r_i-s}$ as formal variables.
\item
The coefficient tuple $\bfy_{t_0}$ at $t_0$ coincides with $\bfy$.
\item
The mutation data $\bfz$ coincides with $\bfz$.
\end{itemize}
\end{defn}

We  need a stronger version 
of the Laurent phenomenon,
which is parallel to Theorem \ref{1thm:Laurent2}.

\begin{thm}[{\cite{Nakanishi14a}}]
\label{1thm:gLaurent2}
For any generalized cluster pattern with principal coefficients at $t_0$, we have
\begin{align}
\label{1eq:gLaurent1}
x_{i;t} \in \bbZ[\bfx_{t_0}^{\pm1}, \bfy,\bfz].
\end{align}
\end{thm}
\begin{proof}
This is proved by showing a similar claim
in the proof of Theorem \ref{1thm:Laurent2}
for both variables $\bfy$ and $\bfz$,
where the  property  \eqref{1eq:gconst1} of $\bfz$ is  essential.
\end{proof}

For a generalized cluster pattern
one can define its (generalized) $C$- and $G$-matrices and $F$-polynomials in two ways:
\begin{itemize}
\item
Define them through a generalized cluster pattern
with principal coefficients at $t_0$
\item
Define them by the underlying (generalized) $B$-pattern $\bfB$ and $t_0$.
\end{itemize}

Below we skip the first definition, which is parallel to the ordinary case,
 and only give the \emph{second\/} one.
One notable difference to the ordinary case is that the $F$-polynomials are
now polynomials in \emph{both $\bfy$ and $\bfz$},
where $\bfz$ are the formal variables $\bfz$
corresponding to the mutation data $\bfz$
(under our usual abuse of notations).

\begin{defn}
For a given generalized cluster pattern $\bfSigma$ with a given
mutation data $(\bfr,\bfz)$ and a given initial vertex $t_0$,
the (generalized) $C$-, $G$, $F$-patterns 
$\bfC^{t_0}$, $\bfG^{t_0}$, $\bfF^{t_0}$
of $\bfSigma$ are uniquely determined by the following
initial conditions and the mutation formulas,
where $t$ and $t'$ are $k$-adjacent:
\begin{align}
\label{1eq:gCmat1}
C_{t_0}&=I,
\\
 \label{1eq:gcmut1}
c_{ij;t'}&=
\begin{cases}
-c_{ik;t} & j = k,\\
c_{ij;t} 
+ c_{ik;t}[r_kb_{kj;t}]_+
+ [-c_{ik;t}r_k]_+b_{kj;t} 
& j \neq k,
\end{cases}
\\
\label{1eq:gGmat1}
G_{t_0}&=I,
\\
 \label{1eq:ggmut1}
 g_{ij;t'}&=
 \begin{cases}
 \displaystyle
 -g_{ik;t}
 + \sum_{\ell=1}^n g_{i\ell;t} [-b_{\ell k;t}r_k]_+
 -  \sum_{\ell=1}^nb_{i\ell;t_0}  [-c_{\ell k;t}r_k]_+ 
 &
 j= k,
 \\
 g_{ij;t}  &
 j \neq k,
 \end{cases}
 \end{align}
 \begin{align}
\label{1eq:gFinit1}
F_{i;t_0}(\bfy,\bfz)&=1,
\\
 \label{1eq:gFmut1}
 F_{i;t'}(\bfy,\bfz)&=
 \begin{cases}
\frac
 {
  \displaystyle
 M_{k;t}(\bfy,\bfz)
 }
{ \displaystyle
 F_{k;t}(\bfy,\bfz)
 }
   &
 i= k,
 \\
 F_{i;t}(\bfy,\bfz)  &
 i \neq k,
 \end{cases}
 \end{align}
where
\begin{align}
\label{1eq:gM1}
\begin{split}
 M_{k;t}(\bfy)
   &=
\Biggl(\,
    \prod_{j=1}^{n}
  y_j^{[-c_{jk;t}r_k]_+}
    \prod_{j=1}^{n}
  F_{j;t}(\bfy,\bfz)^{[-b_{jk;t}r_k]_+}
\Biggr)
\\
&\qquad
\times
 \sum_{s=0}^{r_k}
  z_{k,s}
  \Biggl(\,
 \prod_{j=1}^{n}
  y_j^{c_{jk;t}}
    \prod_{j=1}^{n}
  F_{j;t}(\bfy,\bfz)^{b_{jk;t}}
  \Biggr)^s.
  \end{split}
\end{align}
 \end{defn}

In the first definition of $F$-polynomials, 
$ F_{i;t}(\bfy,\bfz)$  is defined from
the cluster variable
$x_{i;t}$ with principal coefficients in $t_0$
by the specialization $x_1=\cdots = x_n=1$. 
Therefore, 
thanks to Theorem \ref{1thm:gLaurent2},
$ F_{i;t}(\bfy,\bfz)$  is a polynomial
in $\bfy$ and $\bfz$.
.

As expected, the separation formulas for
generalized cluster patterns are given exactly in
the same form as  the ordinary ones
in Theorem \ref{1thm:gsep1}
just by adding the variables $\bfz$ for $F$-polynomials.

\begin{thm}[{Separation Formulas \cite{Nakanishi14a}}]
\index{generalized cluster algebra (GCA)!separation formula}
\label{1thm:gsep1}
Let $\bfSigma$ be any generalized cluster pattern
with coefficients in $\bbP$, mutation data $(\bfr,\bfz)$,
and a given initial vertex $t_0$.
Let 
\begin{align}
\bfx_{t_0}=\bfx,
\quad
\bfy_{t_0}=\bfy,
\quad
\hat\bfy_{t_0}=\hat\bfy
\end{align}
be the initial cluster variables,
coefficients, and $\haty$-variables.
Then,
the following formulas hold.
\begin{align}
\label{1eq:gsep1}
x_{i;t}&=
\Biggl(\,
\prod_{j=1}^n
x_j^{g_{ji;t}}
\Biggr)
\frac{F_{i;t}(\hat{\bfy},\bfz)}{F_{i;t}\vert_{\bbP}(\bfy,\bfz)},
\\
\label{1eq:gsep2}
y_{i;t}&=
\Biggl(\,
\prod_{j=1}^n
y_j^{c_{ji;t}}
\Biggr)
\prod_{j=1}^n
F_{j;t}\vert_{\bbP}(\bfy,\bfz)^{b_{ji;t}},
\\
\label{1eq:gsep3}
\haty_{i;t}&={}
\Biggr(
\prod_{j=1}^n
\haty_j^{c_{ji;t}}
\Biggr)
\prod_{j=1}^n
F_{j;t}(\hat\bfy,\bfz)^{b_{ji;t}}.
\end{align}
\end{thm}
\begin{proof}
One can repeat the proof of
Theorem \ref{1thm:sep1}  
taking care of the modification by the mutation degree $\bfr$.
\end{proof}

There is one notable feature of GCA compared with the ordinary one.
For any generalized cluster pattern 
with mutation data $(\bfr, \bfz)$,
let $\bfB=\{ B_t\}_{t\in \bbT_n}$ be the (generalized) $B$-pattern of $\bfSigma$.
Also, let $\bfC^{t_0}=\{ C_t\}_{t\in \bbT_n}$
and $\bfG^{t_0}=\{ G_t\}_{t\in \bbT_n}$
be the (generalized) $C$- and $G$-patterns of $\bfSigma$.

Let $R=(r_i\delta_{ij})_{i,j=1}^n$ be the diagonal matrix whose diagonal entries are
given by the mutation degrees.
Let $R\bfB=\{ RB_t\}_{t\in \bbT_n}$ and
$\bfB R=\{ B_t R\}_{t\in \bbT_n}$.
Then, one can verify from \eqref{1eq:gbmut2} by inspection
that 
both $R\bfB$ and $\bfB R$ are $B$-pattern by the ordinary matrix mutations.

Let $^L\bfC^{t_0}=\{ ^LC_t\}_{t\in \bbT_n}$ and  $^L\bfG^{t_0}
=\{ ^LG_t\}_{t\in \bbT_n}$
be  the (ordinary) $C$- and $G$-patterns associated with the (ordinary) $B$-pattern $R \bfB $.
Similarly,
let $^R\bfC^{t_0}=\{ ^RC_t\}_{t\in \bbT_n}$ and  $^R\bfG^{t_0}
=\{ ^RG_t\}_{t\in \bbT_n}$
be  the (ordinary)  $C$- and $G$-patterns associated with the (ordinary) $B$-pattern $\bfB R$.

\begin{prop}[\cite{Nakanishi14a}]
The following relations holds:
\begin{align}
\label{1eq:gCLR1}
C_t&={} ^LC_t=R (^RC_t)R^{-1},\\
\label{1eq:gGLR1}
G_t&={} ^RG_t=R^{-1} (^LG_t)R.
\end{align}
\begin{proof}
Let us prove \eqref{1eq:gCLR1}.
The first equality follows from the inspection of \eqref{1eq:gcmut1}.
To obtain the second equality, multiply $r_i^{-1} r_j$ to \eqref{1eq:gcmut1}.
It yields the equality $R^{-1}C_t R={}^RC_t$.
The proof of \eqref{1eq:gGLR1} is similar, 
where we also use \eqref{1eq:gCLR1}.
\end{proof}
\end{prop}

Finally, we present a parallel result to Theorem \ref{1thm:dual2}.

\begin{thm}[Duality \cite{Nakanishi14a}] 
\index{duality}
\label{1thm:dual3}
(a).
For  a common
skew-symmetrizer $D$ of the  (ordinary) $B$-pattern $R\bfB$, we have the equality
\begin{align}
\label{1eq:dual3}
D^{-1} R^{-1} G_t^T RD C_t=I.
\end{align}
(b).
For a common skew-symmetrizer $D$ of 
the  (ordinary) $B$-pattern $\bfB R$,
we have the equality
\begin{align}
\label{1eq:dual4}
D^{-1} R G_t^T R^{-1}D C_t=I.
\end{align}
(c).
For a common skew-symmetrizer $D$ of 
the generalized $B$-pattern $\bfB$,
we have the equality
\begin{align}
\label{1eq:dual5}
D^{-1}  G_t^T D C_t=I.
\end{align}
\end{thm}
\begin{proof}
(a).
By applying Theorem \ref{1thm:dual2} to ${} ^LC_t$ and ${} ^LG_t$,
we obtain the equality
\begin{align}
D^{-1}  ({}^L G_t)^T D ({}^L C_t)=I.
\end{align}
Then, by \eqref{1eq:gCLR1} and \eqref{1eq:gGLR1}, it is written as
\eqref{1eq:dual3}.

(b). This is obtained  in the same way as (a).
Alternatively, note that, if $D$ is a common
skew-symmetrizer  of the  $B$-pattern $R\bfB$,
$D R^2 $ is a common
skew-symmetrizer  of the  $B$-pattern $\bfB R$.
Then,  \eqref{1eq:dual4} follows from \eqref{1eq:dual3}.

(c).  If $D$ is a common
skew-symmetrizer  of the  $B$-pattern $R\bfB$,
 $DR$ is a common
skew-symmetrizer  of the  generalized $B$-pattern $\bfB $.
Then, \eqref{1eq:dual5} follows from \eqref{1eq:dual3}.

\end{proof}

More results on generalized cluster patterns are found in \cite{Nakanishi19}.

%%%%%%%%%%%%%%%%% part1.tex 
%%%%%%%%%%%%%%%%%%%%%%%%%%%%%%%%%

%%%%%%%%%%%%%%%%% part1.tex 
%%%%%%%%%%%%%%%%%%%%%%%%%%%%%%%%%
\newpage
\pagestyle{myheadings}
\bibliography{../../biblist/biblist.bib}

\fontsize{10pt}{10.5pt}\selectfont
\printindex
\end{document}